\documentclass{amsart}
\usepackage{amssymb,amsmath}
\usepackage[mathscr]{euscript}
\input{epsf}

\newcounter{sec}

\newcounter{punct}[sec]

\def\punct{\refstepcounter{punct}{\arabic{sec}.\arabic{punct}.  }}

\newtheorem{theorem}{Theorem}[sec]
\newtheorem{proposition}[theorem]{Proposition}

\newtheorem{lemma}[theorem]{Lemma}

\newtheorem{corollary}[theorem]{Corollary}

\def\COUNTERS{\addtocounter{sec}{1}
              \setcounter{punct}{0}
          \setcounter{equation}{0}
          \setcounter{theorem}{0}
          }
          
          \def\sm{\smallskip}

\begin{document}

\newcommand{\supp}{\mathop {\mathrm {supp}}\nolimits}
\newcommand{\rk}{\mathop {\mathrm {rk}}\nolimits}
\newcommand{\Aut}{\mathop {\mathrm {Aut}}\nolimits}
\newcommand{\Out}{\mathop {\mathrm {Out}}\nolimits}
\renewcommand{\Re}{\mathop {\mathrm {Re}}\nolimits}

\def\ov{\overline}
\def\wh{\widehat}
\def\wt{\widetilde}

\renewcommand{\rk}{\mathop {\mathrm {rk}}\nolimits}
\renewcommand{\Aut}{\mathop {\mathrm {Aut}}\nolimits}
\renewcommand{\Re}{\mathop {\mathrm {Re}}\nolimits}
\renewcommand{\Im}{\mathop {\mathrm {Im}}\nolimits}
\newcommand{\sgn}{\mathop {\mathrm {sgn}}\nolimits}

\def\bfa{\mathbf a}
\def\bfb{\mathbf b}
\def\bfc{\mathbf c}
\def\bfd{\mathbf d}
\def\bfe{\mathbf e}
\def\bff{\mathbf f}
\def\bfg{\mathbf g}
\def\bfh{\mathbf h}
\def\bfi{\mathbf i}
\def\bfj{\mathbf j}
\def\bfk{\mathbf k}
\def\bfl{\mathbf l}
\def\bfm{\mathbf m}
\def\bfn{\mathbf n}
\def\bfo{\mathbf o}
\def\bfp{\mathbf p}
\def\bfq{\mathbf q}
\def\bfr{\mathbf r}
\def\bfs{\mathbf s}
\def\bft{\mathbf t}
\def\bfu{\mathbf u}
\def\bfv{\mathbf v}
\def\bfw{\mathbf w}
\def\bfx{\mathbf x}
\def\bfy{\mathbf y}
\def\bfz{\mathbf z}

\def\bfA{\mathbf A}
\def\bfB{\mathbf B}
\def\bfC{\mathbf C}
\def\bfD{\mathbf D}
\def\bfE{\mathbf E}
\def\bfF{\mathbf F}
\def\bfG{\mathbf G}
\def\bfH{\mathbf H}
\def\bfI{\mathbf I}
\def\bfJ{\mathbf J}
\def\bfK{\mathbf K}
\def\bfL{\mathbf L}
\def\bfM{\mathbf M}
\def\bfN{\mathbf N}
\def\bfO{\mathbf O}
\def\bfP{\mathbf P}
\def\bfQ{\mathbf Q}
\def\bfR{\mathbf R}
\def\bfS{\mathbf S}
\def\bfT{\mathbf T}
\def\bfU{\mathbf U}
\def\bfV{\mathbf V}
\def\bfW{\mathbf W}
\def\bfX{\mathbf X}
\def\bfY{\mathbf Y}
\def\bfZ{\mathbf Z}

\def\frD{\mathfrak D}
\def\frL{\mathfrak L}

\def\bfw{\mathbf w}

\def\R {{\mathbb R }}
 \def\C {{\mathbb C }}
  \def\Z{{\mathbb Z}}
  \def\H{{\mathbb H}}
\def\K{{\mathbb K}}
\def\N{{\mathbb N}}
\def\Q{{\mathbb Q}}
\def\A{{\mathbb A}}

\def\T{\mathbb T}
\def\P{\mathbb P}

\def\G{\mathbb G}

\def\cD{\EuScript D}
\def\cL{\mathscr L}
\def\cK{\EuScript K}
\def\cM{\EuScript M}
\def\cN{\EuScript N}
\def\cR{\EuScript R}
\def\cW{\EuScript W}
\def\cY{\EuScript Y}
\def\cF{\EuScript F}

\def\bbA{\mathbb A}
\def\bbB{\mathbb B}
\def\bbD{\mathbb D}
\def\bbE{\mathbb E}
\def\bbF{\mathbb F}
\def\bbG{\mathbb G}
\def\bbI{\mathbb I}
\def\bbJ{\mathbb J}
\def\bbL{\mathbb L}
\def\bbM{\mathbb M}
\def\bbN{\mathbb N}
\def\bbO{\mathbb O}
\def\bbP{\mathbb P}
\def\bbQ{\mathbb Q}
\def\bbS{\mathbb S}
\def\bbT{\mathbb T}
\def\bbU{\mathbb U}
\def\bbV{\mathbb V}
\def\bbW{\mathbb W}
\def\bbX{\mathbb X}
\def\bbY{\mathbb Y}

\def\kappa{\varkappa}
\def\epsilon{\varepsilon}
\def\phi{\varphi}
\def\le{\leqslant}
\def\ge{\geqslant}

\def\B{\mathrm B}

\def\la{\langle}
\def\ra{\rangle}

\def\F{{}_2F_1}
\def\FF{{}^{\vphantom{\C}}_2F_1^\C}

\newcommand{\dd}[1]{\,d\,{\overline{\overline{#1}}} }

\def\lambdA{{\boldsymbol{\lambda}}}
\def\alphA{{\boldsymbol{\alpha}}}
\def\betA{{\boldsymbol{\beta}}}
\def\mU{{\boldsymbol{\mu}}}

\def\1{\boldsymbol{1}}
\def\2{\boldsymbol{2}}

\def\const{\mathrm{const}}
\def\rem{\mathrm{rem}}
\def\even{\mathrm{even}}
\def\SO{\mathrm{SO}}
\def\SL{\mathrm{SL}}

\def\Pii{\Pi_{\mathrm {cont}}}

\def\aphi{\acute {\phi}}
\def\apsi{\acute {\psi}}

\begin{center}

{\bf\Large A pair of commuting hypergeometric operators

\sm

on the complex plane  and bispectrality}

\bigskip

\sc \large
Vladimir F. Molchanov%
\footnote{Supported by the State order of the
Ministry of Education and Science of the Russian Federation  3.8515.2017.},
Yury A. Neretin%
\footnote{Supported by the grants FWF Der Wissenschaftsfonds, P28421, P31591.}

\end{center}

{\small We consider the standard hypergeometric differential operator
$\frD$ regarded as an operator on the complex plane $\C$ and the 
complex conjugate operator $\ov\frD$. These operators formally commute
and are formally adjoint one to another with respect to an appropriate weight.
We find conditions when they commute in the Nelson sense and
write explicitly their joint spectral decomposition. It is determined by a two-dimensional
counterpart of the  Jacobi transform  (synonyms: 
generalized Mehler--Fock transform, Olevskii  transform).
We also show that the inverse transform is an operator of spectral decomposition
for a pair of commuting difference operators defined in terms of shifts in imaginary
direction.}

\section{Introduction}

\COUNTERS

{\bf\punct Spectral problem.%
\label{ss:1.1}}
Denote by $\dot \C$ the complex plane without the points $0$ and $1$,
by $\cD(\dot \C)$ the space of smooth compactly supported functions
on $\dot\C$. Denote by $\dd z$ the standard Lebesgue measure on $\C$. 

Fix real $a$ and $b$. Consider the following measure on $\dot\C$
\begin{equation}
\mu_{a,b}(z)\,\dd z:=
|z|^{2a+2b-2}|1-z|^{2a-2b}\,\dd{z}
\label{eq:mu}
\end{equation}
 and the corresponding space $L^2(\C,\mu_{a,b})$, 
$$
\la f, g\ra=
\int_\C f(z)\ov{g(z)}\,\mu_{a,b}(z)\,\dd z.
$$
Consider the following pair of differential operators in the space $L^2(\C,\mu_{a,b})$:
\begin{align}
\frD:=&z(1-z)\frac{\partial^2}{\partial z^2}+
\bigl(a+b-(2a+1)z\bigr)\frac\partial{\partial z}-a^2;
\label{eq:frD}
\\
\ov \frD:=&\ov z(1-\ov z)\frac{\partial^2}{\partial \ov z^2}+
\bigl(a+b-(2a+1)\ov z\bigr)\frac\partial{\partial \ov z}
-a^2.
\label{eq:ovfrD}
\end{align}

These operators formally commute, i.e.,
$$
\frD \ov\frD f=\ov\frD \frD f,\qquad \text{where $f\in \cD(\dot\C)$.}
$$
A straightforward calculation shows that they are formally adjoint,
$$
\la \frD f, g\ra=\la f,\ov \frD g\ra, \qquad \text{where $f$, $g\in \cD(\dot\C)$.}
$$
Therefore the operators 
$
\tfrac 12(\frD+\ov\frD)$,
$ \frac1{2i}(\frD-\ov\frD)
$
are  symmetric on the domain $\cD(\dot\C)$.

The purpose of this paper is to construct an explicit  spectral decomposition
of this pair, i.e., a unitary operator $U$, which diagonalizes both  operators  $\frD$, $\ov\frD$.

As we know after the famous work of  Edward Nelson  \cite{Nel}, 1959, (see, also \cite{RS}, Sect. VIII.5)
 a question about commutativity of two unbounded self-adjoint operators 
 can be highly nontrivial%
 \footnote{The topic of the Nelson paper was finite-dimensional Lie algebras of unbounded operators in  Hilbert spaces.
 However, his results remain to be non-trivial for  pairs of commuting operators and even for one operator.}. Recall that two
self-adjoint operators $A$, $B$ commute if  they can be simultaneously realized as operators
of multiplication by functions in some $L^2$.
Equivalently,
the corresponding one-parametric groups commute:
$$
e^{is A} e^{itB}=  e^{itB} e^{is A}, \qquad \text{where $s$, $t$ in $\R$.}
$$
Equivalently, resolvents $(A-\lambda)^{-1}$ and $(B-\mu)^{-1}$, 
commute.
However these properties do not follow from the identity $AB=BA$ and are difficult
for a verification. There are some useful  sufficient conditions and necessary conditions for commutativity
(for necessity we use the result of Kostyuchenko and Mityagin \cite{KM1}-\cite{KM2}),
but quite often a question remains to be heavy%
\footnote{A famous example is a problem, see \cite{Fug}, which was raised by Irving Segal in 1958
and which was discussed  during almost 30 years: Let $\Omega$ be an open connected domain in
$\R^n$. Assume that the operators $i\partial/\partial x_k$ in $\cD(\Omega)$
admit commuting self-adjoint extensions. Is it correct that
$\Omega$ is essentially a fundamental domain of $\R^n$ with respect to a certain discrete group?
The answer is affirmative.}.


\sm

Define two domains $\Pi\supset \Pii$ of the parameters $(a,b)$:
\begin{align}
\Pi:&\quad
0< a+b< 2,\qquad -1<a-b< 1
.
\label{eq:restrictions1}
\\
\Pii:&\quad
0\le a \le 1, \qquad 0\le b\le 1, \quad \text{and $(a,b)\ne (\pm 1,\pm 1)$, $(\pm 1,\mp 1)$.} 
 \label{eq:restrictions2}
\end{align}

\begin{figure}
\epsfbox{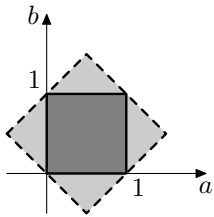}

\caption{To Theorem \ref{th:conditions1}.
The domain $\Pi$ of commutativity, and the domain
$\Pii\subset \Pi$, where the spectrum is purely continuous.}
\end{figure}

\begin{theorem}
	\label{th:conditions1}
%
%
 The operators
$
\tfrac 12(\frD+\ov\frD)$,
 $
\tfrac 1{2i}(\frD-\ov\frD)
$
admit  extensions to a pair of commuting
self-adjoint operators if and only if 
$(a,b)\in \Pi$.
\end{theorem}


Next, we define a natural domain for our operators.
Consider the subspace $\cR_{a,b}(\dot\C)\subset L^2(\dot\C,\mu_{a,b})$
consisting of smooth functions $f$ on $\dot\C$
satisfying the following 
conditions%
\footnote{If $(a,b)\notin\Pi$,
then $\cR_{a,b}(\dot\C)$ is not contained in $L^2(\C,\mu_{a,b})$.}: 

\sm

$1^\circ$. In a neighborhood of $z=0$ a function $f$ has an expansion of the form%
\footnote{Boundary conditions in this spirit sometimes arise in  spectral theory
of ordinary differential operators $D$ for operators with deficiency indices $(1,1)$ or $(2,2)$,
see, e.g., \cite{Ner-jacobi}, Section 1.}
\begin{equation}
f(z)=
\begin{cases}
\alpha(z)+ \beta(z) |z|^{2-2a-2b},\qquad & \text{if $a+b\ne 1$;}
\\
\alpha(z)+ \beta(z) \ln|z|,\qquad & \text{if $a+b=1$,}
\end{cases}
\label{eq:alpha-beta}
\end{equation}
where $\alpha(z)$, $\beta(z)$ are smooth functions.

\sm

$2^\circ$.  In a neighborhood of $z=1$ a function $f$ has an expansion of the form
\begin{equation}
\label{eq:gamma-delta}
f(z)=
\begin{cases}
\gamma(z)+ \delta(z) |z-1|^{2b-2a},\qquad& \text{if $a-b\ne 0$;}
\\
\gamma(z)+ \delta(z) \ln|z-1|,\qquad &\text{if $a-b=0$,}
\end{cases}
\end{equation}
where $\gamma(z)$, $\delta(z)$ are smooth.

\sm

$3^\circ$.  For each $p$, $q$, $N$ we have 
\begin{equation}
\frac{\partial^{p+q} f}{\partial z^p \partial \ov z^q}
=O\Bigl(|z|^{-2a-p-q}\,(\ln|z|)^{-N}\Bigr) \qquad \text{as $z\to\infty$.}
\label{eq:alpha-beta-N}
\end{equation}

\begin{theorem}
\label{th:conditions2}
{\rm a)} For $(a,b)\in \Pi$ 
the
operators $\tfrac 12(\frD+\ov\frD)$,
$\tfrac 1{2i}(\frD-\ov\frD)$ are essentially
self-adjoint on $\cR_{a,b}(\dot\C)$ and
commute in the Nelson sense.

\sm 

{\rm b)}
If $(a,b)\in \Pii$,
then the spectrum of the problem
\begin{equation}
\frD f=\zeta f,\qquad \ov\frD f=\ov\zeta f
\label{eq:sigma-ov-sigma}
\end{equation}
is multiplicity free and
consists of $\zeta$ having the form
\begin{equation}
\zeta=\bigl(\tfrac{k+is}2\bigr)^2,
\qquad
\text{where $k\in \Z$, $s\in\R$.}
\label{eq:sigma=}
\end{equation}
If  $(a,b)\in \Pi\setminus \Pii$, then the spectrum consists on the same
set plus one eigenvalue $\zeta_0>0$.
\end{theorem}

Let us explain the obstacle for  commutativity. Consider a second order
 differential operator $D$
on an interval. For each $\zeta\in\R$ the differential equation $Df=\zeta f$
has two solutions, and we can select generalized eigenfunctions of $D$
as solutions that have $L^2$- or almost $L^2$-asymptotics at the ends of the interval.
In our case the system
(\ref{eq:sigma-ov-sigma})
locally has 4 solutions. Furthermore, $\dot \C$ is not simply connected, solutions
are ramified at $0$, $1$,  $\infty$. As a result there are few  single valued solutions
and we have no freedom for selection of asymptotics. Such  considerations
(see Section \ref{s:differential}) allow to establish necessity of the conditions of Theorem \ref{th:conditions1}.

Unfortunately, we do not know an a priori proof of sufficiency
and obtain it as a byproduct  of the explicit joint spectral decomposition of the
operators $\frD$, $\ov \frD$.
Such detour makes our work long and requires numerous
explicit calculations and estimates.

\sm

{\bf \punct The index hypergeometric transform.}
Our work is a counterpart of the following classical topic.
Consider the hypergeometric differential operator 
$$
D:=x(x+1)\frac {d^2}{dx^2} +\bigl( (a+b)+ (2a+1) x\bigr)\frac d{dx}+a^2
$$
on the half-line $\R_+$, i.e., $x>0$ . Consider the integral operator
\begin{equation}
I_{a,b} f(s):=\frac1{\Gamma(a+b)}\int_0^\infty f(x)\, _2F_1\left[\begin{matrix}
                                                                 a+is,a-is\\a+b
                                                                \end{matrix}
;x\right] \, x^{a+b-1}(1+x)^{a-b}\,dx.
\label{eq:index}
\end{equation}
Then  $I_{a,b}$ is a unitary operator
\begin{equation}
L^2\Bigl(\R_+,  x^{a+b+1}(1+x)^{a-b}\,dx\Bigr)
\to L^2\Bigl(\R_+, 
\pi^{-1}\Bigl|\frac{\Gamma(a+is)\Gamma(b+is)}{\Gamma(2is)}\Bigr|^2 \,ds  \Bigr).
\label{eq:L2even}
\end{equation}
The operator $I_{a,b}$ sends $D$ to the multiplication by $s^2$, see \cite{Wey},
\cite{Tit}, \cite{Ole}, \cite{Koo-jacobi},
\cite{Koo-paley}, \cite{Ner-index}, \cite{Ner-add}.
This transform%
\footnote{A special case $a=1/2$, $b=1$ of this transform
was discovered by Gustav Mehler in 1881, the general transform
was obtained by Weyl \cite{Wey} in 1910.
The $I_{a,b}$ is a representative of a large family of {\it index integral transforms},
which involve
indices of hypergeometric functions, see numerous examples in \cite{Yak},
\cite{Gro}, \cite{Ner-shift}.
} is known as {\it'the generalized Mehler-Fock transform', 'the Olevskii transform', or
'the Jacobi transform'}.

Such operators arise in a natural way in the analysis on rank one Riemannian symmetric spaces,
on the other hand they  are  special cases of multi-dimensional Harish-Chandra spherical transforms
and 
more general Heckman--Opdam  \cite{HO} integral  transforms,
which arise as spectral decompositions of certain families of commuting partial differential
operators.

Next, consider the following difference operator in the space of even functions
depending on the variable $s$:
\begin{multline}
Lg(s):=\frac{(a-is)(b-is)}{(-2is)(-2is+1)} \bigl(g(s+i)-g(s)\bigr)
+\\+
\frac{(a+is)(b+is)}{(2is)(2is+1)} \bigl(g(s-i)-g(s)\bigr),
\label{eq:L-real}
\end{multline}
where $i^2=-1$. A domain of this operator is a space of even functions
holomorphic in the strip $|\Im s|<1+\epsilon$ with some condition of decreasing at infinity.
It turns out that $L$ is  essentially self-adjoint in the space of even functions $L^2_\even\bigl(\R, 
\bigl|\frac{\Gamma(a+is)\Gamma(b+is)}{\Gamma(2is)}\bigr|^2 \,ds  \bigr)$
and the operator $I_{a,b}^{-1}$ sends it to the operator
of multiplication by $x$.

So we have a bispectrality in the spirit 
of Gr\"unbaum \cite{Gru}, \cite{DuG}. Notice
that simpler index integral transforms as the Kontorovich--Lebedev transform
and the $_1F_1$-Wimp transforms also are bispectral, see \cite{Ner-shift}.

Cherednik showed \cite{Che} that inverse Heckman--Opdam transforms
provide spectral decompositions of families of commuting
difference operators, see also van Diejen,  Emsiz \cite{Die}.

\sm

{\bf\punct Radial parts of Laplace operators.%
\label{ss:radial}} Recall one more classical topic. Consider the usual sphere $S^2_\R$:
$$x^2+y^2+z^2=1,$$
the orthogonal  group $\mathrm{SO}(3)$ acts in $L^2(S^2_\R)$ by rotations. Recall one of possible
ways to decompose
 this unitary representation into irreducible components. Consider the Beltrami--Laplace operator $\Delta$ on
 the sphere and restrict it to the space of functions depending on the  height $z$.
 We get a differential operator
 $$
 L_z:=(1-z^2)\frac{\partial^2}{\partial z^2}-2z \frac{\partial}{\partial z}
 $$
 in $L^2[-1,1]$. Eigenfunctions of $L_z$ are the Legendre polynomials.
 Simple arguments show that
 the spectral decomposition of $\Delta$ is a priori equivalent to the spectral decomposition
 of $L_z$ (the reason of this equivalence is compactness of the group $\SO(2)$ of rotations of 
 $S^2_\R$ about the vertical axis).
 
 Now consider the complex manifold $S^2_\C\subset \C^3$ defined by the same equation
 $x^2+y^2+z^2=1$. The complex orthogonal group $\SO(3,\C)$ (the Lorentz group)  acts on the quadric $S^2_\C$, the 
 action
  admits an $\SO(3,\C)$-invariant measure, and again we come to a problem%
  \footnote{This problem was solved by Naimark in \cite{Nai1} in a completely different way.}
  of decomposition
 of the unitary representation of $\SO(3,\C)$ in $L^2$ on $S^2_\C$. Now we have
 two Beltrami-Laplace operators, a holomorphic operator $\Delta$ and
 an antiholomorphic operator $\ov \Delta$. They commute in the Nelson sense.
 Restricting them to functions depending on the coordinate $z\in\C$ we get two operators%
 \footnote{This pair corresponds to $a=b=1/2$ in our parameters.}:
 $$
 L_z:=(1-z^2)\frac{\partial^2}{\partial z^2}-2z \frac{\partial}{\partial z},
 \qquad
   L_{\ov z}:=(1-{\ov z}^2)\frac{\partial^2}{\partial {\ov z}^2}-2{\ov z} \frac{\partial}{\partial {\ov z}}.
 $$
 However, now the stabilizer of the point $(x,y,z)=(0,0,1)$ is a {\it noncompact} subgroup $\SO(2,\C)$,
 and this breaks the a priori argumentation. A joint spectral decomposition of $\Delta$, $\ov \Delta$
 can be reformulated as a certain problem%
 \footnote{Such reductions for families of Laplace operators were widely used by Harish-Chandra
 (in his famous works on the Plancherel formula for real semisimple Lie groups)
 and by his successors.
 The problem for $L_z$, $\ov L_z$ is more similar to decompositions of $L^2$ on real rank one pseudo-Riemannian
 symmetric spaces, which was solved
 by one of the authors of the present paper \cite{Mol1}--\cite{Mol3}.}
 for $L_z$, $\ov L_z$, but this is not
 precisely a problem of a joint spectral decomposition of $L_z$, $\ov L_z$.

 \sm
 
 Notice that a similar separation of variables can be done for $L^2$ on an arbitrary
 rank one complex symmetric space $G_\C/H_\C$ (and, more generally,
 for spaces of $L^2$-sections of line bundles on $G_\C/H_\C$).
 In all the cases we get pairs of hypergeometric operators of our type.
 We hope
 that our spectral decomposition allows to write the explicit Plancherel formula
 for such spaces and to give another proof of old Naimarks's results
 \cite{Nai1}--\cite{Nai3} on tensor products 
 of representations of the Lorentz group. However, the present paper does not have
 such purposes.

\sm

{\bf \punct Homographic transformations of the operators $\frD$, $\ov\frD$.}
Our next purpose is to present the explicit joint spectral decomposition of the pair
 $\frD$, $\ov\frD$. We need some preparations.
 
 Consider the following 
8 transformations of functions
on $\dot\C$:
$$
f(z)\mapsto \gamma_j(z) f(z),\qquad f(z)\mapsto \gamma_j(z) f(1-z), 
$$
where
$$\gamma_j(z)=1,\quad
|1-z|^{2(b-a)},\quad |z|^{2(1-a-b)},\quad |z|^{2(1-a-b)}|1-z|^{2(b-a)},
$$
cf. Erd\'elyi etc., \cite{HTF1}, Subsect. 2.6.1.
It can be readily checked that these transformations send the operators $\frD$, $\ov\frD$
to operators of the same type with other values of the parameters $(a,b)$,
as 
$$(a,b)\mapsto (b,a),\qquad (a,b)\mapsto (1-a,b),\qquad\text{etc.}$$
Thus we get all isometries of the square $\Pi$.
In particular, such transformations send spectral problems to equivalent spectral
problems.

\sm

{\bf\punct Notation. Generalized  powers.}
Denote by $\C^\times$ the multiplicative group  of $\C$.
We need  a notation for  characters 
of $\C^\times$. Let  $z\in \C^\times$  and $a$, $a'\in\C$ satisfy 
$a-a'\in\Z$.
We define a {\it generalized power} of $z$ by
 $$
z^{\bfa}=z^{a|a'}:=z^a \ov z^{a'}= e^{a\ln z+a'\ov{\ln z}}=|z|^{2a}\, \ov z^{\,a'-a},
 $$
Denote by $\Lambda^\C$ the set of all pairs $a|a'$ such that $a-a'\in\Z$.
Denote by $\Lambda\subset \Lambda_\C$ the set of all pairs
 \begin{equation}
 a\bigl|a'
 =\tfrac 12(k+is)\Bigl|\tfrac 12(-k+is), \qquad\text{where $k\in\Z$, $s\in \R$}.
 \label{eq:ks}
 \end{equation}
 Equivalently, $a|a'\in \Lambda$ if $a-a'\in\Z$, $a+a'\in i\R$.
 We also will use the notation
 \begin{equation}
 \label{eq:[]}
 [\bfa]= [a|a']:=\tfrac 12 \Re(a+a').
 \end{equation}
We have
$$
\bigl|z^{a|a'}\bigr|=|z|^{2[a|a']},
$$
in particular, for $\bfa\in\Lambda$ we have $\bigl|z^{a|a'}\bigr|=1$.

\sm

  We fix  the standard Lebesgue measure $\wt d\lambda$ on  the set $\Lambda$:
 $$
 \int_\Lambda \phi(\lambda)\,\wt d \lambda:=\sum_k \int_\R \phi\bigl( \tfrac{k+is}2\bigr)\, ds .
 $$

{\bf\punct Hypergeometric function of the complex field.}
Following Gelfand, Graev, and Retakh \cite{GGR}, we define the gamma function $\Gamma^\C$,
the beta function $B^\C$, and the hypergeometric
function $\FF$ of the complex field (see, also,
Gelfand, Graev, Vilenkin,
\cite{GGV}, Subsect. II.3.7, and Mimachi \cite{Mim}). The {\it gamma function} $\Gamma^\C$ is
\begin{multline}
\Gamma^\C(\bfa)=\Gamma^\C(a|a'):= \frac 1\pi
\int_\C z^{\bfa-\1} e^{2i\Re z} \dd{z}:=
\\:=
\frac 1\pi
\int_\C z^{a-1|a'-1} e^{2i\Re z} \dd{z}= \\=
i^{a-a'}
\frac{\Gamma(a)}{\Gamma(1-a')}=
i^{a'-a} \frac{\Gamma(a')}{\Gamma(1-a)}=
\frac{i^{a'-a}}\pi \Gamma(a)\Gamma(a')\sin \pi a'
.
\label{eq:gamma}
\end{multline}

%
%
%


The {\it beta function $B^\C$} is%
\footnote{This integral has a multi-dimensional counterpart 
of the Selberg type, see \cite{DF}.}
\begin{multline}
\B^\C(\bfa,\bfb):=\frac 1\pi \int_\C t^{\bfa-1}(1-t)^{\bfb-1} \,\dd{t}
=\frac{\Gamma^\C(\bfa)\Gamma^\C(\bfb)}
          {\Gamma^\C(\bfa+\bfb)}
          =\frac{\Gamma(a)\Gamma(b)\Gamma(1-a'-b')}
                {\Gamma(a+b)\Gamma(a')\Gamma(b')}.
\label{eq:beta}
\end{multline}

The {\it hypergeometric function of the complex field} is defined by
\begin{multline}
\FF[\bfa,\bfb;\bfc;z]=
\FF\Bigl[\begin{matrix}\bfa;\bfb\\ \bfc\end{matrix};z\Bigr]
=
\FF\Bigl[\begin{matrix}a|a',b|b'\\ c|c'\end{matrix};z\Bigr]
:=\\:=
\frac1{\pi B^\C(\bfb,\bfc-\bfb)}
\int_\C
t^{\bfb\bf-1}(1-t)^{\bfc-\bfb\bf-1}(1-zt)^{-\bfa}
\dd{t}.
\label{eq:def-FF}
\end{multline}
Recall that the Gauss hypergeometric functions are defined by
$$
\F[a,b;c;z]:=
\frac1{ B(b,c-b)}
\int_0^1
t^{b-1}(1-t)^{c-b-1}(1-zt)^{-a}
dt=
\sum_{p=0}^\infty\frac{(a)_p(b)_p}{ p!\,(c)_p} \, z^p,
$$
where $(c)_p:=c(c+1)\dots (c+p-1)$ is the Pochhammer symbol.
The functions $\FF[\bfa,\bfb;\bfc;z]$ admit expressions in the terms of
$_2F_1$, see Theorem \ref{th:long}.

\sm

{\bf\punct Spectral decomposition.}
For  $(a, b)\in\Pi$ we define the kernel
$\cK_{a,b}(z,\lambda)$ on $\C\times \Lambda$ by
\begin{equation}
 \cK_{a,b}(z,\lambda)=
 \frac1{\Gamma^\C(a+b|a+b)}
 \,
 \FF\left[\begin{matrix}
  a+\lambda|a-\ov\lambda,\,
 a-\lambda|a+\ov\lambda
  \\ a+b|a+b
                               \end{matrix}
                               ;z
\right].
\label{eq:cK}
\end{equation}

\begin{theorem}
	\label{th:spectral1}
 Let $(a,b)\in\Pii$.
 Then the operator
 $$
 J_{a,b} f(\lambda) := \int_\C \cK_{a,b}(z,\lambda)\,f(z) \,\mu_{a,b}(z)\,\dd z
 $$
 is a unitary operator from $L^2(\C,\mu_{a,b})$ to
 $L^2_\even(\Lambda, \kappa_{a,b})$ of even functions on $\Lambda$ with respect to the
 Plancherel measure  
 \begin{equation}
dK_{a,b}(\lambda)=\kappa_{a,b} (\lambda)\,\wt d\lambda
 =\frac 1{4\pi^2} \Bigl|
 \lambda\,
 \Gamma^\C(a-\lambda|a+\ov\lambda)\,
 \Gamma^\C(b+\lambda|b-\ov\lambda)\Bigr|^2
 \,\wt d\lambda.
 \label{eq:plancherel}
 \end{equation}
\end{theorem}
 
 Next, we  modify the definition of the measure for $(a,b)\in \Pi\setminus\Pii$.
  Due to the homographic transformations%
  \footnote{Changing of kernels  $\cK_{a,b}$ by the  homographic transformation
  can be observed from Proposition \ref{pr:kummer}.}
   it is sufficient to consider the case $a<0$.
  We define the Plancherel measure $d K_{a,b}(\lambda)$ on $\Lambda_\C$
  that is the sum of $\kappa_{a,b}\,\wt d\lambda$ and two
  $\delta$-measures located at the points
  $\pm a|\pm a\in\Lambda_\C$,
  \begin{equation}
   \Gamma^\C(a+b|a+b)\,\Gamma^\C(b-a|b-a)\,\Gamma^{\C}(2a|2a) \cdot \bigl(\delta_{a|a}+ \delta_{-a|-a}\bigr).
  \end{equation}
Define a constant function $v(z)$ on $\dot\C$
by
$$
v(z)=\Gamma^\C(a+b|a+b)^{-1}.
$$
For $f\in \cD(\dot\C)$ we define an even function $J_{a,b}(\lambda)$ on the support of $d K_{a,b}(\lambda)$
given by the same formula (\ref{eq:plancherel}) on $\Lambda$,
its value at $(\pm a|\pm a)$ is
$$
J_{a,b} f(\pm a|\pm a):=\la f,v\ra_{L^2(\C,\mu_{a,b})}.
$$

\begin{theorem}
\label{th:spectral1add}
If $(a,b)\in \Pi$  and $a<0$, then
the operator $J_{a,b}$ is  unitary as an operator $L^2(\C,\mu_{a,b})$ to
$L^2_\even(\Lambda_\C,d K_{a,b})$.
\end{theorem}

Our operator really determines the spectral decomposition:

\begin{theorem}
 For each $(a,b)\in \Pi$
  for any 
$f\in \cD_{\even}(\dot \C)$ we have 
$$
 J_{a,b} \frD f(\lambda)= \lambda^2 J_{a,b} f(\lambda),
 \qquad J_{a,b} \ov\frD f(\lambda)= \ov \lambda^2 J_{a,b} f(\lambda).
$$
\end{theorem}

This means that
$$
\frD \cK(z,\lambda)=\lambda^2 \cK(z,\lambda), \qquad \ov\frD \cK(z,\lambda)=\ov\lambda^{\,2} \cK(z,\lambda).
$$

 Next, we consider the space $\cD_\even(\dot\Lambda)$, which consists of
 even
 smooth compactly supported functions on $\Lambda$
 that are zero on a neighborhood of the point $0|0$.
The following statement explains the appearance of the space $\cR_{a,b}$
and also is one of the arguments for the proofs of our main statements.

\begin{theorem}
\label{th:image-J*}	
 If $F\in \cD_\even(\dot \Lambda)$, then $J^*_{a,b}F\in \cR_{a,b}$.
\end{theorem}

The images of $\delta$-functions also are contained in $\cR_{a,b}$.

\sm

{\bf \punct The transformation $J_{a,b}$ in the complex domain.%
\label{ss:complex-domain}}
Let us extend our kernel $\cK$ to the complex domain.
For 
$$
\{\lambda|\lambda'\}=\bigl\{\tfrac{k+\sigma}2\bigl|\tfrac{-k+\sigma}2\bigr\}\in \Lambda_\C
$$
we set
\begin{multline}
\cK_{a,b}(z; \lambda|\lambda')=
\cK_{a,b}(z; k, \sigma):= \\:=
\frac 1{\Gamma^\C(a+b|a+b)}\,
\FF\left[\begin{matrix}
       a+\frac{k+\sigma}2\bigl| a+\frac{-k+\sigma}2,\,\, a+\frac{-k-\sigma}2\bigl| a+\frac{k-\sigma}2 
       \\
       a+b|a+b
                          \end{matrix}; z
                          \right],
                          \label{eq:cK-complex}
\end{multline}
 where $k$ ranges in $\Z$, $\sigma$ ranges in $\C$. The previous expression (\ref{eq:cK})
 corresponds to a pure imaginary $\sigma$.

 For $f\in\cD(\dot\C)$ we define a meromorphic function on $\Lambda_\C$
 by
 $$
 J_{a,b} f(k,\sigma):=\int_{\dot\C} f(z)\,\cK(z;k,\sigma)\,d\mu_{a,b}\dd z.
 $$

 \begin{theorem}
 \label{th:paley}
For  $f\in\cD(\dot\C)$ the function  $J_{a,b} f$
is contained in the space $\cW_{a,b}$ defined as follows.  
 \end{theorem}
 
We define a space $\cW_{a,b}$ as the space of all meromorphic functions%
\footnote{We say that a function $F(k,\sigma)$
is meromorphic if it is meromorphic as a function in $\sigma$ for any fixed $k$.}
$F(k,\sigma)$ on $\Lambda_\C$ satisfying the conditions a)--d):
 
{\rm  a)}  $F$ is even, i.e., $F(-k,-\sigma)=F(k,\sigma)$.

\sm
 
{\rm  b)}
 Possible poles of $F(k,\sigma)$  
are located at the points
\begin{equation}
\sigma=\pm(-2a+|k|+2j), \qquad \pm(-2b+|k|+2j),\qquad \text{where $j=1$, $2$, $3$, \dots}
 \label{eq:possible-poles-1}
\end{equation}
A maximal possible order of a pole at a point $(l,c)$
is a multiplicity of $(l,c)$ in the collection%
\footnote{For $(a,b)\in \Pi$ orders of poles $\le 2$. Poles of order 2
arise only if $a=b$, $a=1$, $b=1$.} (\ref{eq:possible-poles-1})

\sm

{\rm  c)} For each $A>0$ for each $N>0$ in the union of  strips $|\Re \sigma|<A$  we have an estimate
\begin{equation}
 F(k,\sigma)=O\bigl(k^2+(\Im \sigma)^2\bigr)^{-N} \qquad \text{as  $k^2+(\Im \sigma)^2\to\infty$}.
\end{equation}

\sm

{\rm d)} For each $p$, $q\in \Z$
\begin{equation}
  F(p,q)=  F(q,p).
  \label{eq:kk}
\end{equation} 

\sm

 Next, we extend the spectral density $\kappa_{a,b}$ to the complex domain.
 \begin{multline}
  \kappa_{a,b} (\lambda|\lambda')= \kappa_{a,b} (k,\sigma):=\\:=
  \frac 1{4\pi^2} 
   (k+\sigma)(k-\sigma)\,
  \Gamma^\C  \bigl(a+\tfrac{k+\sigma}2\bigl| a+\tfrac{-k+\sigma}2\bigr)
 \,
 \Gamma^\C  \bigl(a+\tfrac{-k-\sigma}2\bigl| a+\tfrac{-k-\sigma}2\bigr)
 \times\\
 \Gamma^\C  \bigl(b+\tfrac{k+\sigma}2\bigl| b+\tfrac{-k+\sigma}2\bigr)
 \,
 \Gamma^\C  \bigl(b+\tfrac{-k-\sigma}2\bigl| b+\tfrac{-k-\sigma}2\bigr).
 \label{eq:spectral-density}
 \end{multline}
 
 In the case $a<0$ discussed above, $ \kappa_{a,b}$ has a pole at $k=0$, $\sigma=a$ and the inner product in 
 $L^2_\even(\Lambda, dK_{a,b})$ can be written
 as
 \begin{multline*}
 \la F,G\ra=
 \frac 1 i\sum_k \int_{-i\infty}^{i\infty}
 F(k,\sigma)\,\ov{G(k,-\ov \sigma)} \kappa_{a,b}(k,\sigma)\,d\sigma
 +\\+
2\, \mathrm{res}_{s=a} \Bigl( F(k,\sigma)\,\ov{G(k,-\ov \sigma)} \kappa_{a,b}(0,\sigma)\Bigr).
 \end{multline*}
 If $a>1$, then the spectral density has a zero at $k=0$, $\sigma=a-1$ but
  both  functions
 $F(k,\sigma)$, $\ov{G(k,-\ov \sigma)}$ admit simple poles at this point, and we have 
 a similar formula.

\sm

{\bf\punct Difference spectral problem.}
It turns out that our problem is bispectral, and the bispectrality is a crucial argument 
of our proof.
We define analogs of the difference operator (\ref{eq:L-real}).
%
Consider meromorphic functions $\Phi$ depending on 
$$\lambda|\lambda'=\tfrac12(k+is)\Bigl|\tfrac 12(-k+is)\in\Lambda_\C$$
and the operators in the space of meromorphic functions
defined by
$$
T \Phi(k,s)=\Phi(k+1,s-i),\qquad \wt T \Phi(k,s)=\Phi(k+1,s+i),
$$
or, equivalently,
\begin{equation}
T \Phi(\lambda|\lambda')= \Phi(\lambda+1|\lambda'), \quad \wt T \Phi(\lambda|\lambda')= \Phi(\lambda|\lambda'+1).
\label{eq:shifts}
\end{equation}
We define the following difference operators
\begin{align}
 \frL&:=& \frac{(a+\lambda)(b+\lambda)}{2\lambda(1+2\lambda)}\, (T-1)
 +  \frac{(a-\lambda)(b-\lambda)}{-2\lambda(1-2\lambda)} \,(T^{-1}-1);
 \label{eq:frL}
 \\
  \ov\frL&:=& \frac{(a+\lambda')(b+\lambda')}{2\lambda'(1+2\lambda')}
  \, (\wt T^{-1}-1)
 +  \frac{(a-\lambda')(b-\lambda')}{-2\lambda(1-2\lambda')}\, (\wt T-1).
 \label{eq:ovfrL}
\end{align}

Formally,
$$
\frL \ov\frL =\ov\frL \frL.
$$

\begin{theorem}
	\label{th:spectral2}

 {\rm a)} The
 operators $\frac 12 (\frL+\ov \frL)$, $\frac 1{2i} (\frL-\ov \frL)$
 defined on  the space $\cW_{a,b}$
 are essentially self-adjoint and commute in the Nelson sense.

 \sm
	
{\rm b)}
 For $\Phi\in J_{a,b} \cD(\dot\C)$ we have
 \begin{equation}
  J_{a,b}^{-1} \,\frL\, \Phi(z)= z\, J_{a,b}^{-1}\, \Phi(z), \qquad   J_{a,b}^{-1}\, \ov\frL\, \Phi(z)= \ov z \,J_{a,b}^{-1}\, \Phi(z).
 \end{equation}
\end{theorem}

Thus the operator $J_{a,b}^{-1}$  determines a  joint spectral decomposition
of $\frac 12 (\frL+\ov \frL)$ and $\frac 1{2i} (\frL-\ov \frL)$.

\sm

{\bf \punct The structure of the proofs.}
We derive asymptotics of the kernel $\cK(z,\lambda)$ as $z\to 0$, $1$, $\infty$ for fixed
$\lambda$ (Theorem \ref{th:long}) and  as $|\lambda|\to\infty$ for fixed $z$ (Theorem \ref{th:asymptotic}).
Next, we prove inclusions
$$
J_{a,b}^* \,\cD_\even(\dot \Lambda) \subset \cR_{a,b},\qquad
J_{a,b}^{\phantom{*}}\, \cD(\dot \C)\subset \cW_{a,b} 
$$
(Proposition \ref{l:JRab} and Corollary \ref{cor:JCL2}) and symmetries 
\begin{align}
\la \frD f,g\ra_{L^2(\C,\mu_{a,b})}=\la f,\ov \frD g\ra_{L^2(\C,\mu_{a,b})}
,\qquad \text{where $f$, $g\in \cR_{a,b}$};
\\
\la \frL F,G\ra_{L^2(\dot\Lambda,dK_{a,b})}=\la  F,\ov\frL G\ra_{L^2(\Lambda,dK_{a,b})},\quad
\text{where $F$, $G\in \cW_{a,b}$},
\end{align}
see Proposition \ref{pr:D-symmetry} and Theorem \ref{th:symmetry}.
This  implies a generalized orthogonality, i.e., 
$$
\la J^*_{a,b} F,  J^*_{a,b} G\ra_{L^2(\C,\mu_{a,b})}=0 \quad \text{
if supports of $F$, $G\in \cD_\even(\dot\Lambda)$ are disjoint,}
$$
and a similar statement for $J_{a,b}$, see Lemmas \ref{l:packet2},
\ref{l:orthogonal2}.
Next, we show that for any $F$, $G\in \cD_\even(\dot\Lambda)$ 
the inner products of their preimages can be written as:
$$
\la J^*_{a,b} F,  J^*_{a,b} G\ra_{L^2(\C,\mu_{a,b})}=
\la F,G\ra_{L^2(\Lambda,dK_{a,b})}+
\int\limits_\Lambda \int\limits_\Lambda H(\lambda_1,\lambda_2) F(\lambda_1) G(\lambda_2)\,\wt d\lambda_1\,\wt d\lambda_2,
$$
where $H$ is a locally integrable function, see Lemma \ref{l:orthogonal2}. We also prove a similar statement
for $J_{a,b}$, see Lemma \ref{l:zero}. Then generalized orthogonality implies $H(\cdot,\cdot)=0$.
Thus we get 
\begin{equation}
J^*_{a,b}J_{a,b}^{\phantom{*}}=1,\qquad J_{a,b}^{\phantom{*}} J^*_{a,b}=1,
\end{equation}
and this is our main statement.

Some steps of this double way are straightforward, some points require long calculations and estimates, and 
we meet some points of good luck (the proofs of Theorem \ref{th:symmetry}
and  Lemma \ref{l:zero}). We also need a lot of information
about functions $\FF$ (in particular, to cover the cases $a+b\in \Z$ and $a-b\in\Z$ 
we need  a tedious examination of possible degenerations
of functions $\FF$).

 The bispectrality   allows to avoid a direct proof of completeness of the
system of generalized eigenfunctions of $\frD$, $\ov\frD$.

To prove the necessary conditions of self-adjointness in Theorem \ref{th:conditions1}  we analyze  
common generalized  eigenfunctions
of $\frD$, $\ov\frD$ for $(a,b)\notin \Pi$ and after a natural 
selection we reduce a set of possible candidates to a finite family. 
 This is done in Section \ref{s:differential}. 

\sm

This text is focused to a proof of unitarity of $J_{a,b}$. 
An introduction to functions $\FF$ in Section \ref{s:hypergeometric}
can be a point of 
an
independent interest.
Also, we get two relatively pleasant statements 
about asymptotic behavior of  integrals
$$
M(\epsilon)=
\int_{\C}
t^{\alpha-1|\alpha'-1}(\epsilon-t)^{\beta-1|\beta'-1}
\psi(t)\,\dd{t}\qquad\text{as $\epsilon\to0$}
$$
and 
$$
I(\lambda)=
\int_\C |f(t)|^2 e^{i\Re \bigl(\lambda\phi(t)\bigr)}\,\dd t\qquad \text{as $|\lambda|\to\infty$,}
$$
where $f$, $\phi$ are holomorphic and
$\lambda\in\C$ (Theorems \ref{th:weak-singularities}
and \ref{th:parametric}).

\sm

{\bf\punct Final remarks.} The index hypergeometric transform
(\ref{eq:index})
can be applied as a heavy tool of theory of special functions, see
\cite{Koo-jacobi}, \cite{Ner-w}, \cite{Ner-add}.  In 
 \cite{Ner-dou} we use our operators $J_{a,b}$ to obtain a beta integral
 over $\Lambda$, which is a counterpart of the Dougall $_5H_5$-summation 
 formula and of the de Branges--Wilson integral.

 Also, we notice that functions, which can be regarded 
 as higher hypergeometric functions  $^{\vphantom \C}_4F^\C_3$ of the complex field, arise in a natural way in 
 the work of Ismagilov \cite{Ism} as analogs of the Racah coefficients
 for {\it unitary} representations of the Lorentz group
 $\SL(2,\C)$ (see, also a continuation in \cite{DS}).
 
  \sm
 
It seems that
our problem can be a representative of some family of spectral  problems, but
now it is too early to claim  something certainly.

\section{Preliminaries. Gamma function, the Mellin transform, weak singularities \label{s:preliminaries}}

\COUNTERS

This section is a union of 3 disjoint topics:

\sm

--- some properties of the function $\Gamma^\C$, which are intensively used below:

\sm

--- some properties of the Mellin transform on $\C$, they are used in a proof of Proposition \ref{pr:FF-meromorphic}
 and in Sections \ref{s:asymptotic}--\ref{s:isometry2}:
 
 \sm
 
 --- a lemma from asymptotic analysis, which is used only in a proof of Theorem \ref{th:long}
 (the last statement can be independently  established by a straightforward tiresome way):

 \sm


{\bf \punct Some properties of the gamma function.}
The usual functional equations for the $\Gamma$-function can be easily rewritten 
for $\Gamma^\C$ (recall that $a-a'\in \Z$!):
\begin{align}
 &\Gamma^\C(a|a')=  \Gamma^\C(a'|a);
 \\
  &\Gamma^\C(a+1|a')= i\,a\, \Gamma^\C(a|a');
 \\
 &\Gamma^\C(a|a')\,\Gamma^\C(1-a|1-a')=(-1)^{a-a'};
  \label{eq:complement}
 \\
 &\ov {\Gamma^\C(a|a')}= (-1)^{a-a'}\Gamma^\C(\ov a|\ov a')
 .
\label{eq:ov-gamma}
\end{align}
Also,
$$
\prod_{p=0}^{m-1}
\Gamma^\C\Bigl(a+\frac {p-1}m\Bigl|a'+\frac {p-1}m \Bigr)=
m^{1-m(a+a')}\,
\Gamma^\C(ma|ma').
$$

The identity (\ref{eq:ov-gamma}) implies
\begin{equation}
 \ov{\phantom{a} B^\C[a|a',b|b']\phantom{a}}=
 B^\C[\ov a|\ov a',\ov b|\ov b'].
 \label{eq:ov-beta}
\end{equation}

Let
$k_1$, $k_2\in\Z$. Then
\begin{equation}
\Gamma^\C(k_1|k_2)=\begin{cases}
                    \infty, &\text{if $k_1$, $k_2\in \Z_-$},
                    \\
                    0,  &\text{if $k_1$, $k_2\in \N$,}
                    \\
                    i^{k_1-k_2} \frac{(k_1-1)!}{(-k_2)!}, &
                    \text{if $k_1\in \N$, $k_2\in \Z_-$},
                    \\
                    i^{k_2-k_1} \frac{(k_2-1)!}{(-k_1)!}, &
                    \text{if $k_2\in \N$, $k_1\in \Z_-$},
                   \end{cases}
                   \label{eq:Gamma-poles}
\end{equation}
where $\Z_-$ denotes the set of integers $\le 0$.

\sm

The following lemma gives  us the asymptotics of the Plancherel density
(\ref{eq:plancherel}).

\begin{lemma}
\label{l:beta-as}
We have the following asymptotics in $\lambda\in\Lambda$: 
\begin{equation}
\Gamma^\C(a-\lambda|a+\ov\lambda)\,\, \Gamma^\C( b+\lambda|b-\ov \lambda)
\sim \lambda^{a+b-1|a+b-1} \qquad
\text{\rm as $|\lambda|\to \infty$}.
\end{equation}
The asymptotics is uniform in $a$, $b$ if they range in a bounded domain.
\end{lemma}

{\sc Proof.}
Denote $\Re \lambda=k/2$.
Let $|\arg\lambda|<\pi-\epsilon$.  Then we write our expression as
$$
\frac{i^k\Gamma(a+\ov\lambda)}{\Gamma(1-a+\lambda)}
\cdot\frac{i^{-k} \Gamma(b+\lambda)}{\Gamma(1-b+\ov \lambda)}
$$
and apply the standard asymptotic formula $\Gamma(z+\alpha)/\Gamma(z+\beta)\sim z^{\alpha-\beta}$
in the sector $|\arg z|<\pi-\epsilon$, see Erd\'elyi etc., \cite{HTF1}, formula (1.18.4).
If $|\arg(-\lambda)|<\pi-\epsilon$, we write
$$
\frac{i^{-k}\Gamma(a-\lambda)}{\Gamma(1-a-\ov\lambda)}\cdot\frac{i^k \Gamma(b-\ov\lambda)}{\Gamma(1-b-\lambda)}
$$
and come to the same asymptotics.
\hfill $\square$


\sm

{\bf\punct The Mellin transform.%
\label{ss:mellin}}
Denote by $\C^\times:=\C\setminus 0$ the multiplicative group of $\C$.
The {\it Mellin transform} (see, e.g., \cite{GGR}) on $\C^\times$ is defined by 
\begin{equation}
g(\mU)=\cM f(\mU)=\frac 1{2\pi}\int_\C f(z) z^{\mU-\1}\dd z, 
\label{eq:mellin}
\end{equation}
where $\mU=\{\mu|\mu'\}=\bigl\{\frac{k+is}2|\frac{-k+is}2\bigr\}\in\Lambda_\C$
(here we allow complex $s$).
This operator is the Fourier transform on the  group $\C^\times\simeq (\R/2\pi\Z)\times \R$,
so it is reduced to the Fourier transforms on $(\R/2\pi\Z)$ and on $\R$.
Indeed,
changing variables
$$z=e^\rho e^{i\phi}$$
we come to
$$
g(k,s)=\frac 1{2\pi} \int_0^{2\pi}\int_{-\infty}^\infty f\bigl(e^\rho e^{i\phi}\bigr) e^{ik\phi+is\rho}
\,d\rho\,d\phi.
$$

The inversion formula is given by 
$$
f(z)=\cM^{-1} g(z)  = \frac 1 {2\pi}\int_\Lambda g\bigl(\mu\bigl|-\ov\mu\bigr)  z^{-\mu|\ov \mu}\, \wt d \mu.
$$
Equivalently,
$\cM$ is a unitary operator 
$
L^2(\C^\times,|z|^{-2})\to L^2(\Lambda)
$.

{\bf\punct The Mellin transform of even functions.\label{ss:Mellin-even}}
We say that a function $f$ on $\C^\times$ is {\it $\times$-even} if $f(z^{-1})=f(z)$.
Denote by $L^2_+(\C^\times, |z|^{-2})$ the subspace of $L^2(\C^\times, |z|^{-2})$
consisting of $\times$-even
functions. Obviously, the Mellin transform sends $\times$-even functions in $z$
to even functions in $\mu$. Also, for a $\times$-even function
$f$ we have 
\begin{equation}
\cM f(\mU)=\frac12\int_\C f(z) (z^{\mU-\1}+z^{-\mU-\1})\dd z, \qquad\text{where $f$ is $\times$-even.}
\label{eq:even-Mellin}
\end{equation}

{\bf \punct The Mellin transform of smooth compactly supported functions.}

\begin{theorem}
\label{th:mellin-meromorphic}
{\rm a)} Let $f$ be a compactly supported smooth function on $\C$.
 Then $\cM f(\mu|\mu')$ extends to a meromorphic function in the variable $\mu$
 with possible poles located at the points
 $\mu|\mu'\in \Z_-\times \Z_-$. 
  Moreover, for any $p$, $p'\in\Z_+$
 for $\Re(\mu+\mu')>-p-p'$
 we have
 \begin{equation}
  I(\mu|\mu')=\frac{(-1)^{p+p'}}{2\pi(-\mu)_p \,(-\mu')_{p'}}\int_\C z^{\mu-p-1|\mu'-p'-1} 
  \frac{\partial^{p+p'}}{\partial z^p\,\partial \ov z^{\,p'}}
 f(z)\,\dd{z}.
  \label{eq:meromorphic-continuation}
 \end{equation}
  The residues
 at the poles are
 \begin{equation}
 \mathrm{res}_{\mu|\mu'=-p|-p'} I(\mu|\mu')=\frac {1}{(p-1)!\,(p'-1)!}
 \frac{\partial^{p+p'}f(0,0)}{\partial z^p\,\partial \ov z^{\,p'}}
 .
 \label{eq:residues-M}
\end{equation}
 
 {\rm b)} For each $N$ for each $A$ for all pairs $(k,s)$ satisfying $|\Im s|<A$ we have
 $$
 \cM f\bigl(\tfrac {k+is}2|\tfrac {-k+is}2\bigr)=O(k^2+|s|^2)^{-N}\qquad\text{as $|k^2|+|s^2|\to\infty$.}
 $$
 
 \end{theorem}

For a proof of statement a),
see Gelfand, Shilov \cite{GSh}, Sect. B.1.3, or equivalently 
Russian edition of Gelfand, Graev, Vilenkin \cite{GGV}, Addendum 1.3
(the term 'Mellin transform' in that place is absent, but the statement is proved).
Formula (\ref{eq:meromorphic-continuation}) is obtained from
(\ref{eq:mellin}) by integration by parts. The statements about location of poles and about residues
require more careful  arguments.

\sm

{\sc Proof of statement} b).
We pass to polar coordinates, $z= e^{i\theta}$
and get 
$$
 \cM f\bigl(\tfrac {k+is}2|\tfrac {-k+is}2\bigr)=\frac 1{2\pi}\int_0^\infty \int_0^{2\pi} H(\theta, r)\, r^{-1+is} 
e^{i\theta k} \,d\theta\,dr,
$$
where
$H(\theta,r):=\Phi(r e^{i\theta})$ is a smooth function $2\pi$-periodic in $\theta$, 
the  $H(\theta,0)$ does not depend on $\theta$, also
$H(\theta+\pi, -r)=H(\theta,r)$. Integrating by parts in $r$, we  get
$$
 \cM f\bigl(\tfrac {k+is}2|\tfrac {-k+is}2\bigr)= 
 \frac{(-1)^l}{2\pi\, (is)_l}
 \int_0^{2\pi} \int_0^\infty  \frac {\partial^l}{\partial r^l} H(\theta,r)\,
r^{-1+is+l}\,dr\,\, e^{i\theta k} \,d\theta.
$$
For   $l>A$  the integral
absolutely converges.
Integrating by parts in $\theta$, we get 
$$
\frac{(-1)^{k+l}}{2\pi\, (i k)^m (is)_{l}}
\int_0^{2\pi}
\int_0^\infty  \frac {\partial^{l+m}}{\partial r^{l}\,\partial \theta^m} H(\theta,r) 
r^{-1+is+l}dr  e^{i\theta k} \,d\theta,
$$
and 
$$
| \cM f\bigl(\tfrac {k+is}2|\tfrac {-k+is}2\bigr)|\le \frac{\mathrm{const}}{|(2\pi i k)^m (is)_{l}|}
.
$$
If $|s|>|k|$ we take $m=0$ and large $l$, if $|k|>|s|$, we take $l>|\Im s|$ and large $m$.
\hfill $\square$

{\bf\punct Weak singularities.}
Here we imitate one standard trick of asymptotic
analysis, see, e.g., \cite{Fed}, Sect. I.4.
Fix $R$ and a smooth function $\psi(t)$ on $\C$.
Consider the integrals of the following type
$$
M(\epsilon)=M_{\alphA,\betA}(\epsilon):=
\int_{|t|\le R}  t^{\alpha-1|\alpha'-1}(\epsilon-t)^{\beta-1|\beta'-1}
\psi(t)\,\dd{t}.
$$
Clearly, $M_{\alphA,\betA}(\epsilon)$
 is holomorphic in $\alphA$, $\betA$ in the domain of convergence
  and admits a meromorphic continuation%
\footnote{For instance, see the proof of Proposition \ref{pr:FF-meromorphic} below.}
to 
$(\alpha,\beta)\in\Lambda^2$.

\begin{theorem}
\label{th:weak-singularities}
Let $\alphA$, $\betA$ satisfy
the  condition
\begin{equation}
\alphA,\, \betA,\, \alphA+\betA-\1\notin \Z_-\times \Z_-. 
\label{eq:weak-condition}
\end{equation}
 Then $M(\epsilon)$ {\rm(}defined in the sense of analytic 
continuation{\rm )}
admits the following asymptotic expansion
at $0$:
\begin{multline}
M(\epsilon)\sim
\sum_{i,i'\ge 0}
B^\C(\alpha+i|\alpha'+i',\,\beta|\beta')
\cdot
\frac 1{i!\,i'!}\frac{\partial^{i+i'}\psi(0,0)}
{\partial t^i \partial \ov t^{i'}}
\cdot
\epsilon^{\alpha+\beta+i-1|\alpha'+\beta'+i'-1}
+\\+
\sum_{j,j'\ge 0} r_{j|j'} \epsilon^{j|j'}.
\label{eq:weak-singularities}
\end{multline}
The coefficients of the expansion are meromorphic in $\alpha|\alpha'$, $\beta|\beta'$.

\sm

If\,%
\footnote{Recall notation (\ref{eq:[]}).}
$[\alpha|\alpha']>0$, $[\beta|\beta']>0$, $[\alpha|\alpha']+[\beta|\beta']>1$,
 then
\begin{equation}
r_{00}=\int_{|t|<R} t^{\alpha+\beta-2|\alpha'+\beta'-2}\psi(t)
\dd{t}.
\label{eq:r00}
\end{equation}
\end{theorem}

First, we prove the following lemma

\begin{lemma}
\label{l:cut-beta}
 Let $\alpha$, $\alpha'$, $\beta$, $\beta'$ satisfy
the  condition {\rm(\ref{eq:weak-condition})}.
 Then
the following integral {\rm(}defined in the sense of analytic continuation{\rm)}
\begin{equation}
R(\epsilon)=
\int_{|t|<R}  t^{\alpha-1|\alpha'-1}(\epsilon-t)^{\beta-1|\beta'-1}
\,\dd{t}
\label{eq:tge0}
\end{equation}
 admits an asymptotic expansion
of the form
\begin{equation}
R(\epsilon)=B^\C(\alpha|\alpha',\beta|\beta')\cdot
 \epsilon^{\alpha+\beta-1|\alpha'+\beta'-1}+
 \sum_{j\ge 0, j'\ge 0} p_{j|j'}\epsilon^{j|j'}.
\label{eq:incomplete-beta}
\end{equation}
Moreover, the series 
\begin{equation}
\sum_{j\ge 0, j'\ge 0} p_{j|j'}\epsilon^{j|j'}
\label{eq:pp}
\end{equation}
  converges
in the circle $|\epsilon|<1/R$, and the coefficients $p_{j|j'}(\bfa,\bfb)$ are holomorphic
in the domain {\rm(\ref{eq:weak-condition})}.
\end{lemma}

{\sc Proof.}
Set
$$
Q_{\alphA,\betA}(\epsilon_1,\epsilon_2):=
(-1)^{\beta'-\beta}
\int\limits_{|t|>R}
t^{\alpha+\beta-2|\alpha'+\beta'-2}
\Bigl(1-\frac{\epsilon_1} t\Bigr)^{\beta-1} 
\Bigl(1-\frac{\epsilon_2}t\Bigr)^{\beta'-1}
\,\dd t.
$$

This function is meromorphic in $\alpha$, $\beta$, and in $\epsilon_1$,
$\epsilon_2$ in the bidisk $|\epsilon_1|<1/R$, $|\epsilon_2|<1/R$.
Let
\begin{equation}
[\alpha|\alpha']>0,\quad  [\beta|\beta']>0,\quad [\alpha|\alpha']+[\beta|\beta']<1.
\label{eq:restriction}
\end{equation}
Under these conditions the integral $R(\epsilon)$ converges, and
$$
R(\epsilon)=
\int_\C -\int_{|t|>R}
=
B^\C(\alpha|\alpha',\beta|\beta')\cdot
\epsilon^{\alpha+\beta-1|\alpha'+\beta'-1}- Q_{\alphA,\betA}(\epsilon,\ov\epsilon)
.
$$
Expanding the integrand in $Q_{\alphA,\betA}$ in a series in $\epsilon_1$,
$\epsilon_2$ and 
integrating termwise we come to
\begin{multline}
Q_{\alphA,\betA}(\epsilon_1,\epsilon_2)=
\\
=
\sum_{j\ge 0,j'\ge 0:\, \alpha+\beta-j=\alpha'+\beta'-j'}
\frac {(-\beta+1)_j\,(-\beta'+1)_j\, R^{\alpha+\alpha'+\beta+\beta'-j-j'}}
{(j+j'-\alpha-\alpha'-\beta-\beta')\, j!\,j'! } \epsilon_1^j\epsilon_2^{j'}
.
\label{eq:Q}
\end{multline}
Now we can omit restrictions (\ref{eq:restriction}).
Indeed, under  conditions (\ref{eq:weak-condition})
the series (\ref{eq:Q}) converges in the bidisk $|\epsilon_1|<1$, $|\epsilon_2|<1$
and therefore its sum coincides with the meromorphic continuation.
\hfill $\square$

\sm

{\sc Proof of Theorem \ref{th:weak-singularities}.}
We expand the function $\psi$ as a sum
$$
\psi(t,\ov t)=\sum_{j+j'\le N} \frac 1{j!\,j'!} \frac{\partial^{j+j'}\psi(0)}{\partial t^j\,\partial t^{j'}}
\,
t^{j|j'} + H_N(t),
$$
where $H_N(t)$ is a smooth function and
$$
H_N(t)=O(|t|^{N+1})\qquad \text{as $t\to 0$}.
$$
Substituting this to the initial integral we get a sum of integrals of the form
(\ref{eq:tge0}),  we apply Lemma \ref{l:cut-beta} to each summand. Also we get a summand
$$
I(\epsilon)=
\int_{|t|\le R} t^{\alpha-1|\alpha'-1} (\epsilon+t)^{\beta-1|\beta'-1} H_N(t) \dd{t}.
$$
We wish to show that $T(\epsilon)$ has partial derivatives at 0 up to order $N-k$, where $k$ is constant
depending only on $\alphA$ and $\betA$. Consider a partition of unity,
$1=\phi_1+\phi_2$ such that $\phi_2$ is zero at  some smaller  circle $|t|<R'$.
According to this, we split $I=I_1+I_2$.
Obviously, $I_2$ has an expansion of the form
$$
I_2\sim \sum_{j,j'\ge 0} c_{j|j'} \epsilon^{j|j'}
$$
with coefficients meromorphic in $\alphA$, $\betA$. Next, we integrate $I_1$ by parts several times,
$$
I_1(\epsilon)= \frac 1{(\beta)_m(\beta')_m}\int_{|t|\le R} (\epsilon+t)^{\beta-1+m\bigl|\beta'+m-1} 
\frac{\partial t^{2m}} {\partial^m \partial \ov t^m}
\Bigl( t^{\alpha-1|\alpha'-1}  H_N(t)\phi_1(t)\Bigr) \dd{t}.
$$
Choosing $m$ we can make $\beta+m-1$, $\beta'+m-1$ as large, as we want, say $> q$. Next, we choose a large $N$,
such that $\frac{\partial^{2m}}{\partial t^m\partial \ov t^m}\bigl(\dots\bigr)$ is continuous at 0. Now we can 
differentiate
$I_1(\epsilon)$ with respect to $\epsilon$, $\ov\epsilon$ \, $q$ times  at 0 and  consider its Taylor expansion.
This finishes a derivation of the asymptotic expansion for $R(\epsilon)$.

If the integral $R(0)$
converges, we substitute $\epsilon=0$ to the expansion and get the expression
for $r_{00}$.
\hfill $\square$

\section{The hypergeometric function of the complex field\label{s:hypergeometric}}

\COUNTERS

Here we discuss basic properties of the functions $\FF[\cdot]$.

\sm

{\bf\punct Domain of convergence and  analytic continuation.}
The hypergeometric function $\FF[\bfa,\bfb;\bfc;z]$
of the complex field
is defined by the Euler type integral (\ref{eq:def-FF}):
\begin{equation}
\FF[\bfa,\bfb;\bfc;z]
=
\frac1{\pi B^\C(\bfb,\bfc-\bfb)}
\int_\C
t^{\bfb\bf-1}(1-t)^{\bfc-\bfb\bf-1}(1-zt)^{-\bfa}
\dd{t}.
\label{eq:def}
\end{equation}

For $z\ne 0$, $1$, the integral absolutely converges (see notation (\ref{eq:[]}) if 
$\bfa$, $\bfb$, $\bfc$ is contained in the following tube $\wh\Xi$,
\begin{equation}
\wh\Xi:\quad
[\bfb]>0,\quad [\bfc]-[\bfb]>0, \quad  [\bfa]<1, \quad  [\bfc]-[\bfa]<2.
 \label{eq:convergence-hyper}
\end{equation}
%
%
In other words, the integral absolutely converges
if and only if the point $\bigl([\bfa], [\bfb], [\bfc]\bigr)$ is contained in the simplex $\Xi$  in $\R^3$ with vertices
\begin{equation}
 (1,0,0),\qquad (-1,0,0), \qquad (1,0,1), \qquad (1,2,2).
 \label{eq:simplex}
\end{equation}
 We have
$ \Lambda^\C\simeq \C\times \Z
 ,$
therefore triples ($\bfa$, $\bfb$, $\bfc$) depend on 3 integers and 3 complex parameters.
Clearly, each component of the set $\Z^3\times \C^3$ has an open intersection with
 the domain of convergence%
 \footnote{The map $(\bfa,\bfb,\bfc)\to \bigl([\bfa], [\bfb], [\bfc]\bigr)$ from $\Lambda^3\to \R^3$
 is surjective on all components.}.

\begin{proposition}
\label{pr:FF-meromorphic}
For $z\in \dot \C$, the expression  $\FF[\bfa,\bfb;\bfc;z]$ as a function of 
$\bfa$, $\bfb$, $\bfc$ admits a meromorphic extension to arbitrary values of 
$\bfa$, $\bfb$, $\bfc$
with  poles at a countable union of surfaces
\begin{gather}
\bfa\in \N\times\N,\,\,
\bfb\in \N\times\N,\,\,
\bfc-\bfa\in \N\times\N,\,\,
\bfc-\bfb\in \N\times\N,\,\,
\label{eq:poles}
\\
\bfc\in \Z_-\times \Z_-
\label{eq:poles-1}
\end{gather}
and vanishes for all $z\in\dot\C$ at
\begin{equation}
\bfc\in \N\times \N.
\label{eq:zeros-F}
\end{equation}
\end{proposition}

\sm

{\sc Proof.}
Consider a partition of unity 
$1=\phi_0(t)+ \phi_1(t)+\phi_{1/z}(t)+ \phi_\infty(t)+\phi_{\varnothing} (z)$,
where all summands are smooth and nonnegative, $\phi_0$, $\phi_1$, $\phi_{1/z}$,  $\phi_\infty$ are zero outside small neighborhoods of
of $0$, $1$, $1/z$, $\infty$ respectively, and $\phi_{\varnothing}=0$ in small neighborhoods of these points.
Denote $P(t,\ov t)$ the integrand in the integral representation of
$\FF[\bfa,\bfb;\bfc;z]$. Then
\begin{multline*}
\pi \B^\C(\bfb, \bfc-\bfb)\,\,
\FF[\bfa,\bfb;\bfc;z]=\\=\int \phi_0 P\,\dd{t}+ \int \phi_1 P\,\dd{t}+\int \phi_{1/z} P\,\dd{t} + 
\int \phi_\infty P\,\dd{t} +\int \phi_{\varnothing} P\,\dd{t}.
\end{multline*}
The last summand is an entire function  
in $\bfa,$ $\bfb$, $\bfc$.
By Theorem \ref{th:mellin-meromorphic} other
 summands are meromorphic and can have poles 
at
$$
\bfb\in\Z_-\times \Z_-,\,\, \bfc-\bfb\in \Z_-\times \Z_-,\,\, \bfa\in \N\times\N,\,\, \bfc-\bfa\in  \N\times\N.
$$
However, $B^\C$-factor in the front of the integral (\ref{eq:def}) kills 
the first and the second families of poles and  produces new poles and also zeros.
This gives us (\ref{eq:poles})--(\ref{eq:zeros-F}), in particular the factor $\Gamma^\C(\bfc)$
produces poles (\ref{eq:poles-1}) and zeros (\ref{eq:zeros-F}).

All these possible poles really are poles, the simplest way to observe this is to look at formulas
(\ref{eq:A})--(\ref{eq:C11}) derived below. Formulas (\ref{eq:A})--(\ref{eq:A11}) show that
  (\ref{eq:poles}) are poles. To check a presence of poles (\ref{eq:poles-1}) we apply
(\ref{eq:C})--(\ref{eq:C11}).
\hfill $\square$

\sm

{\bf\punct Kummer symmetries.}  This section contains a collection of elementary formulas,
they partially depend on Theorem \ref{th:long} proved below. However, our
proof of this theorem is based on differential equations and asymptotic analysis and is independent of our formulas.

First we notice two trivial identities
\begin{align}
\FF\left[\begin{matrix}
         a|a',b|b'\\
         c|c'
        \end{matrix};\, \ov z
 \right]
 =
 \FF\left[\begin{matrix}
         a'|a,b'|b\\
         c'|c
        \end{matrix};\,  z
 \right];
 \label{eq:z-conjugate}
 \\
 \ov{
 \phantom{aa}
  \FF\left[\begin{matrix}
         a'|a,b'|b\\
         c'|c
        \end{matrix};\,  z\right]
\phantom{aa}
}
= \FF\left[\begin{matrix}
         \ov a'|\ov a,\ov b'|\ov b\\
         \ov c'|\ov c
        \end{matrix};\,  z
 \right]^{\vphantom{\bigl|}}.
 \label{eq:f-conjugate}
\end{align}

To verify (\ref{eq:z-conjugate}) we substitute $t\mapsto \ov t$ to the integral
(\ref{eq:def}).

\begin{proposition} {\rm  a) (Gauss identity)}
	Let 
	$
	[\bfc]-[\bfa]-[\bfb]>0
	$.
	Then
\begin{equation}
\FF\left[\begin{matrix}
         \bfa,\bfb\\
         \bfc
        \end{matrix};\, 1
 \right]:=\lim_{z\to 1}
 \FF\left[\begin{matrix}
         \bfa,\bfb\\
         \bfc
        \end{matrix};\, z
 \right]
 =
 \frac{\Gamma^\C(\bfc)\, \Gamma^\C(\bfc-\bfa-\bfb) }{\Gamma^\C(\bfc-\bfa)
 	\, \Gamma^\C(\bfc-\bfb)}.
 \label{eq:Gauss-identity}
\end{equation}

{\rm b)} Let%
\footnote{If $[\bfc]\ge 1$, then $\lim_{z\to 0} \bigl|\FF[\bfa,\bfb;\bfc;z]\bigr|=\infty$.}
$
[\bfc]<1.
$
Then
$$
\FF[\bfa,\bfb;\bfc;0]:=\lim_{z\to 0} \FF[\bfa,\bfb;\bfc;z]=1.
$$

\end{proposition}

{\sc Proof.} a)
We substitute $z=1$ to (\ref{eq:def}) and come to a beta function,
$$\pi B^\C[\bfb, \bfc-\bfa]/\pi B^\C[\bfb,\bfc-\bfc].$$
However, this argument is valid only if the beta integral
$B^\C[\bfb, \bfc-\bfa]$ converges.
The general  statement follows from Theorem
\ref{th:long}.b proved below.

\sm

b) also is reduced to a beta-function with the same problem
with the domain of convergence. The general statement follows from Theorem
\ref{th:long}.a.
\hfill $\square$  


\begin{proposition}
 \begin{equation}
\FF[\bfa,\bfb;\bfc;z]=\FF[\bfb,\bfa;\bfc;z].
\label{eq:ab-ba}
\end{equation}
\end{proposition}

This will become obvious after Theorem \ref{th:long}.
We use this symmetry in the next two proofs.

\begin{proposition} {\rm(}Euler and Pfaff transformations{\rm)},
\begin{align}
\FF\left[\begin{matrix}
         \bfa,\bfb\\
         \bfc
        \end{matrix};\, z
 \right]&=
 (1-z)^{-\bfa}\,
 \FF\left[\begin{matrix}
         \bfa,\bfc-\bfb\\
         \bfc
        \end{matrix};\, \frac z{z-1}
 \right]
 \label{eq:pf1}
 \\
 &=
  (1-z)^{-\bfb}\,
 \FF\left[\begin{matrix}
         \bfc-\bfa,\bfb\\
         \bfc
        \end{matrix};\, \frac z{z-1}
 \right]
  \label{eq:pf2}
 \\
 &=
   (1-z)^{\bfc-\bfa-\bfb}\,
 \FF\left[\begin{matrix}
         \bfc-\bfa,\bfc-\bfb\\
         \bfc
        \end{matrix};\, z
 \right]
  \label{eq:pf3}
\end{align}
\end{proposition}

{\sc Proof.}
We substitute $t=1-s$ to (\ref{eq:def}) and get (\ref{eq:pf1}).
Applying (\ref{eq:ab-ba}) we get  (\ref{eq:pf2}). Applying (\ref{eq:pf1}) and (\ref{eq:pf2}),
we get (\ref{eq:pf3}).
\hfill $\square$

\newpage

\begin{proposition}
\label{pr:kummer}
\rm (Kummer symmetries)
{\it The following functions $u_j^\C(z)$ are equal}%
\footnote{The meaning of subscripts $j$ in $u_j^\C$, references,
	and comments
are explained in a remark after the proof.}:
\begin{equation}
u^\C_1(z)= \FF\left[\begin{matrix}
         \bfa,\bfb\\
         \bfc
        \end{matrix};\, z
 \right]
\qquad\qquad\qquad \text{(compare with \rm \cite{HTF1}, (2.2.9.1))};
\label{eq:u1}
\end{equation}
\begin{multline}
u^\C_5(z)= 
\\ =\frac{(-1)^{\bfc-\bfa-\bfb}\Gamma^\C(\bfc)\,\Gamma^\C(\bfc-1) }{\Gamma^\C(\bfa)\, \Gamma^\C(\bfb)\,\Gamma^\C(\bfc-\bfa)\, \Gamma^\C(\bfc-\bfb)}
z^{1-c}
\FF\left[ \begin{matrix}
\bfb-\bfc+\1;\, \bfa-\bfc+\1
\\ \2-\bfc
\end{matrix};\, z
\right]
\label{eq:u5}
\end{multline}
\begin{equation*}
\,\,\,\,\text{(see  \cite{HTF1}, (2.2.9.17)) and ratio of coefficients
at $u_1$, $u_5$ in (2.2.10.35));}
\end{equation*}
\begin{align}
u^\C_3(z)&= \frac{\Gamma^\C(\bfc)\,\Gamma^\C(\bfb-\bfa)}{\Gamma^\C(\bfb)\,\Gamma^\C(\bfc-\bfa)}\,\,
(-z)^{-\bfa} \,\, \FF\left[\begin{matrix}
\bfa, \bfa-\bfc+\1
\\
\bfa-\bfb+\1
\end{matrix};\, z^{-1}
\right]
\label{eq:u3}
\\
 & \qquad\qquad \qquad\qquad\qquad\qquad \text{(see  \cite{HTF1}, (2.2.9.9) and $B_1$ in (2.2.10.5));}
\nonumber
\\
u^\C_4(z)&=
 \frac{\Gamma^\C(\bfc)\, \Gamma^\C(\bfa-\bfb)}{\Gamma^\C(\bfa)\,\Gamma^\C(\bfc-\bfb)}
(-z)^{-\bfb} \FF\left[\begin{matrix}
\bfb, \bfb-\bfc+\1\\
\bfb-\bfa+\1
\end{matrix}; \, z^{-1}
\right]
\label{eq:u4}
\\
& \qquad\qquad \qquad\qquad\qquad\qquad  \text{(see  \cite{HTF1}, (2.2.9.10) and $B_2$ in (2.2.10.5));}
\nonumber
\\
u^\C_2(z)&=
\frac{\Gamma^\C(\bfc)\, \Gamma^\C(\bfc-\bfa-\bfb) }{\Gamma^\C(\bfc-\bfa) \Gamma^\C(\bfc-\bfb)}\,
\FF\left[\begin{matrix}
\bfa,\bfb\\
\bfa+\bfb+\1-\bfc
\end{matrix};\, 1-z
\right]
\label{eq:u2}
\\
& \qquad\qquad \qquad\qquad\qquad\qquad  \text{(see  \cite{HTF1}, (2.2.9.5) and $A_1$ in (2.2.10.5));}
\nonumber
\\
u^\C_6(z)&= 
\frac{\Gamma^\C(\bfc)\,\Gamma^\C(\bfa+\bfb-\bfc)}{\Gamma^\C(\bfa)\Gamma^\C(\bfb)} \,
(1-z)^{\bfc-\bfa-\bfb} 
\FF\left[ \begin{matrix}
\bfc-\bfa,\bfc-\bfb\\
\bfc -\bfa-\bfb+\1
\end{matrix}; \,1-z
\right]
\label{eq:u6}
\\ 
&
\qquad\qquad \qquad\qquad\qquad\qquad  \text{(see  \cite{HTF1}, (2.2.9.5) and $A_1$ in (2.2.10.5)).}
\nonumber
\end{align}
\end{proposition}

{\sc Remark.}
For each expression (\ref{eq:u1})--(\ref{eq:u6})
we can apply one of the transformations (\ref{eq:pf1})-(\ref{eq:pf3}).
In this way we get 24 expressions of this type.
\hfill $\boxtimes$

\sm

{\sc Proof.} The formula for $u_3$. Changing a variable $t=1/s$
 in (\ref{eq:def}) we come to
 \begin{multline*}
 \frac{(-1)^{\bfc-\bfa-\bfb-\1} z^{-a}}{\pi B^\C(\bfb,\bfc-\bfb)}
 \int_\C s^{\bfa-\bfc}  (1-s)^{\bfc -\bfb-\1} (1-s/z)^{-\bfa}\dd s
 =\\
 \frac{(-1)^{\bfc-\bfb-\1} \pi B^\C (\bfa-\bfc+\1, \bfc-\bfb) }
 {\pi B^\C(\bfb,\bfc-\bfb)} (-z)^{-a}
 \,\,\FF \left[\begin{matrix}
 \bfa, \bfa-\bfc+\1
 \\
 \bfa-\bfb+\1
 \end{matrix};\, z^{-1}
 \right].
 \end{multline*}
 We cancel $\Gamma^\C(\bfc-\bfb)$ and apply
 (\ref{eq:complement}) two times.
 
 The formula for $u_4$. We transpose $\bfa$ and $\bfb$ in the formula for
 $u_3$.
 
 The formula for $u_5$. We combine the transformations (\ref{eq:u3}) and (\ref{eq:u4}).

 The formula for $u_2$. We combine the transformations (\ref{eq:u3}), (\ref{eq:pf1}),
  and again (\ref{eq:u3}).
  
 We combine transformations (\ref{eq:u3}), (\ref{eq:pf2}),
 and again (\ref{eq:u3}). 
 \hfill $\square$
 
 \sm

{\sc Remark.}  Proposition \ref{pr:kummer} is a self-closed
collection of identities. However, they are reflections of the 
Kummer table of solutions of the hypergeometric equation
$$
\Bigl(z(1-z)\frac{\partial^2}{\partial z^2}+
\bigl(c-(a+b+1)z\bigr)\frac\partial{\partial z}-ab\Bigr)\, u(z)=0,
$$
see Erd\'elyi, et al., \cite{HTF1}, Section 2.2.9, formulas (1)--(24).
The Kummer table
contains 6 solutions, each of them is defined in a neighborhood of one of the  singular
points
$0$, 1, $\infty$.
\begin{align*}
u_1(z)&= \alpha_1(z),  &u_5(z)&= z^{1-c} \alpha_5(z),
\\
u_3(z)&=(-z)^{-a} \alpha_3(z^{-1}), &
u_4(z)&=(-z)^{-b} \alpha_4(z^{-1}),\\
u_2(z)&= \alpha_2(1-z), &
u_6(z)&= (1-z)^{c-a-b}\alpha_6(1-z),
\end{align*}
where $\alpha_j(x)$ are power series, $\alpha_j(0)=1$.
Generally, these solutions are ramified at the points $0$, $1$, $\infty$.
Each solution is represented in 4 forms, which can be obtained one from another by the Pfaff
transformations, see  Erd\'elyi, et al., \cite{HTF1}, Sect. 2.1, (22)--(23). In the table above
we present the corresponding expressions for $\FF[\bfa,\bfb; \bfc; z]$,  they 
correspond to Kummer's expressions with change $(a,b,c)\mapsto (\bfa,\bfb, \bfc)$. The resulting
functions $u_j^\C$ are non-ramified (by definition) and
differ by factors independent of $z$, we normalize them to make them equal
one to another. Counterparts
of these  factors (except one formula)
are present in the Kummer formulas  as coefficients of transfer-matrices
$(u_1,u_5)$ to $(u_3, u_4)$ and $(u_2,u_6)$, with the same replacement  $(a,b,c)\mapsto (\bfa,\bfb, \bfc)$,
see Erd\'elyi, et al., \cite{HTF1}, display (2.2.10.5) and the coefficients $A_1$, $A_2$, $B_1$, $B_2$.
So, in each line of  Proposition \ref{pr:kummer} we give a reference to the corresponding formula
in Erd\'elyi, et al., \cite{HTF1}, (2.2.9.1)--(2.2.9.24) and to the corresponding coefficient in \cite{HTF1}, display (2.2.10.5). 
\hfill $\boxtimes$

\sm

{\bf\punct Differential equations.}

\begin{lemma}
\label{l:differentsiruem}
\begin{align*}
\frac\partial{\partial z}\,\,
\FF\Bigl[\begin{matrix}a|a';b|b'\\ c|c'\end{matrix};z\Bigr]
=
\frac{ab}{c}\,\,
\FF\Bigl[\begin{matrix}a+1|a';b+1|b'\\ c+1|c'\end{matrix};z\Bigr];
\\
\frac\partial{\partial \ov z}\,\,
\FF\Bigl[\begin{matrix}a|a';b|b'\\ c|c'\end{matrix};z\Bigr]
=
\frac{a'b'}{c'}\,\,
\FF\Bigl[\begin{matrix}a|a'+1;b|b'+1\\ c|c'+1\end{matrix};z\Bigr].
\end{align*}
\end{lemma}

{\sc Proof.} We differentiate the integral with respect
to the parameter $z$, and get an integral of the same form.
The calculation is valid if
$\Xi\cap\bigl( \Xi+(\frac 12, \frac 12,\frac 12)\bigr)\ne \varnothing$,
where $\Xi$  is the simplex defined by (\ref{eq:convergence-hyper})--(\ref{eq:simplex}).
This intersection is open and nonempty. It remains to refer to the meromorphic
continuation.
\hfill $\square$

\smallskip

Denote
\begin{align}
D=D[a,b,c]:=z(1-z)\frac{\partial^2}{\partial z^2}+
\bigl(c-(a+b+1)z\bigr)\frac\partial{\partial z}-ab;
\label{eq:D[]}
\\
D'=D'[a',b',c']:=\ov z(1-\ov z)\frac{\partial^2}{\partial \ov z^2}+
\bigl(c'-(a'+b'+1)\ov z\bigr)\frac\partial{\partial \ov z}-a'b'.
\label{eq:D-prime[]}
\end{align}

\begin{proposition}
\label{pr:system-hypergeometric}
The complex hypergeometric function $\cF(z)=\FF[\bfa,\bfb;\bfc;z]$  satisfies
the following system of partial differential equations
\begin{equation}
D[a,b,c]\,\,\cF
=0;
\qquad\quad
D'[a',b',c']\,\,
\cF=0.
\label{eq:two-equations}
\end{equation}
\end{proposition}

We call these equation by {\it complex hypergeometric system}.

\sm

{\sc Proof.} This follows from the identity
\begin{equation}
D[a,b,c] \Bigl(t^{b-1}(1-t)^{c-b-1}(1-tz)^{-a}\Bigr)
=-a\frac\partial{\partial t}
\Bigl( t^{b}(1-t)^{c-b}(1-tz)^{-a-1}\Bigr)
\label{eq:identity}
\end{equation}
(cf. \cite{HTF1}, (2.1.3.11)).
Consider sufficiently small positive $\epsilon$, $\delta$, $\kappa$
and take $\bfa$, $\bfb$, $\bfc$ such that
$$
[b|b']=\epsilon, \quad [c|c']=\epsilon+\delta,\qquad [a|a']=-\tfrac 12 +\epsilon+\delta+\kappa.
$$
We multiply both parts 
 of (\ref{eq:identity})
 by $\ov t^{\,b'-1}(1-\ov t)^{c'-b'-1}(1-\ov t\ov z)^{-a'}$
 and integrate 
 over $\C$.
In the left hand side for such values of the parameter we can permute  integration and differentiation in $z$.
In the right hand side the integrand is an integrable  derivative of an integrable function. Therefore the right hand side is zero.
\hfill $\square$

\smallskip

\begin{proposition}
\label{pr:apriory-hyper}
{\rm a)} Any solution of  system {\rm (\ref{eq:D[]})--(\ref{eq:D-prime[]})} is real analytic in $z$.

\sm

{\rm b)}
Let $z_0\ne 0$, $1$, $\infty$.
Denote by $\phi_1(z)$, $\phi_2(z)$ a pair of independent
holomorphic solutions of the ordinary differential equation
$D[a,b,c]\,f(z)=0$ at a neighborhood of $z_0$. Denote
by $\psi_1(\ov z)$, $\psi_2(\ov z)$ a pair of antiholomorphic
solutions of the ordinary differential equation
$D'[a',b',c']\,f(\ov z)=0$. Then 
any solution of the system {\rm (\ref{eq:D[]})--(\ref{eq:D-prime[]})}
can be represented as
\begin{equation}
\label{eq:Dij}
\sum_{i,j=1,2} \tau_{ij}(\bfa,\bfb,\bfc) \phi_i(z)\psi_j(\ov z).
\end{equation}

{\rm c)} If we choose $\phi_i$, $\psi_j$ meromorphic in the parameters
$\bfa$, $\bfb$, $\bfc$
in some domain in $\Lambda^3_\C$,
 then the coefficients $\tau_{ij}$ also are meromorphic
in the parameters $\bfa$, $\bfb$, $\bfc$.
\end{proposition}

{\sc Proof.}
a)  Indeed, $D[a,b,c]$ is an elliptic differential
operator, therefore solutions of the equation
$D \cF=0$ are analytic functions, see, e.g., \cite{Hor},Theorem 9.5.1.

\sm

 b) Consider a solution
$$
\FF[\bfa,\bfb;\bfc;z]= h_{00}+h_{10}(z-z_0)+h_{01}(\ov z-\ov z_0)+
h_{11}( z- z_0)(\ov z-\ov z_0)+\dots
$$
 of the system of partial differential equations
(\ref{eq:two-equations}). These equations determine
recurrence relations for the Taylor coefficients $h_{ij}$ of
$\FF[\dots]$ at $z_0$. It  can be easily checked that
all the coefficients $h_{ij}$
admit linear expressions in terms of
 $h_{00}$, $h_{01}$, $h_{10}$, $h_{11}$.
  On the other hand, for given
$h_{00}$, $h_{01}$, $h_{10}$, $h_{11}$,
   we can find
a local solution of the complex hypergeometric  system (\ref{eq:two-equations})
in the form
$\sum C_{ij} \phi_i(z)\psi_j(\ov z)$.

\smallskip

c) By Lemma \ref{l:differentsiruem}, the coefficients
$h_{00}$, $h_{10}$, $h_{01}$, $h_{11}$
depend on $\bfa$, $\bfb$, $\bfc$ meromorphically.
If $\phi_i(z_0)$, $\phi_i'(z_0)$,
$\psi_j(z_0)$, $\psi_j'(z_0)$ are meromorphic in the
parameters, then the coefficients $C_{ij}$ also are meromorphic.
\hfill $\square$

\smallskip

{\bf \punct Expressions for $\FF$.} Let us write expansions of $\FF[\dots;z]$ near the singular points
$z=0$, 1, $\infty$.
Explicit formulas for fundamental systems of solutions
of the hypergeometric differential equation
are well-known, see Erd\'elyi, et al., \cite{HTF1}, 2.9 (the Kummer series).
For definiteness, consider $z_0=0$. If $c\notin\Z$, then
for generic values of the parameters
the hypergeometric equation
$D[a,b,c]f(z)=0$ has two holomorphic (ramified)  solutions 
on
a punctured neighborhood
of 0,
$$
\phi_1(z)=\F[a,b;c;z],\qquad \phi_2(z)=z^{1-c}
F[a+1-c,b+1-c; 2-c;z].
$$
The equation $D'[a',b',c']f(\ov z)=0$
has two antiholomorphic solutions
$$
\psi_1(z)=\F[a',b';c';\ov z],\qquad \psi_2(z)=\ov z^{1-c'}
F[a'+1-c',b'+1-c'; 2-c';\ov z].
$$
Therefore near $z=0$ we have 
solutions of  system (\ref{eq:D[]})--(\ref{eq:D-prime[]})
of the same  form (\ref{eq:Dij}) with new $\phi$, $\psi$.
We get a family of functions depending of 4 parameters $\tau_{ij}$,
therefore for generic $\bfa$, $\bfb$, $\bfc$ this formula gives all
multivalued solutions near $z=0$.

Solutions (\ref{eq:Dij}) that are single valued 
in a neighborhood of 0 have the form
\begin{equation}
A_1 \phi_1(z) \psi_1(\ov z) + A_2 \phi_2(z) \psi_2(\ov z).
\label{eq:AB}
\end{equation}

\begin{theorem}
\label{th:long}
{\rm a)} In the disk $|z|<1$ we have the following expansion:
\begin{multline}
\FF\Bigl[\begin{matrix}\bfa;\bfb\\ \bfc\end{matrix};z\Bigr]
=A_{0}\cdot
\F\Bigl[\begin{matrix}a,b\\c  \end{matrix}; z  \Bigr]
\,
\F\Bigl[\begin{matrix}a',b'\\c'  \end{matrix}; \ov z  \Bigr]
+\\+
A_{1}\cdot z^{1-c|1-c'}
\F\Bigl[\begin{matrix}a+1-c,b+1-c\\2-c  \end{matrix}; z  \Bigr]
\,
\F\Bigl[\begin{matrix}a'+1-c',b'+1-c'\\2-c'  \end{matrix};
 \ov z  \Bigr],
 \label{eq:A}
\end{multline}
where
\begin{align}
A_{0}&=1,
\label{eq:A00}
\\
A_{1}&=
(-1)^{\bfc-\bfa-\bfb} \frac{\Gamma^\C(\bfc)\,\Gamma^\C(\bfc-\1)}
{\Gamma^\C(\bfa)\, \Gamma^\C(\bfb)\, \Gamma^\C(\bfc-\bfa)\,\Gamma^\C(\bfc-\bfb)}.
\label{eq:A11}
 \end{align}

{\rm b)} In the disk $|z-1|<1$ the following expansion holds:
\begin{multline}
\FF\Bigl[\begin{matrix}\bfa;\bfb\\ \bfc\end{matrix};z\Bigr]=
B_{0}\cdot
\F\Bigl[\begin{matrix}a,b\\a+b+1-c\end{matrix}; 1-z  \Bigr]
\F\Bigl[\begin{matrix}a',b'\\a'+b'+1-c'\end{matrix};1-\ov z  \Bigr]
+\\+
B_{1} \cdot
(1-z)^{c-a-b|c'-a'-b'}
\F\Bigl[\begin{matrix}c-a,c-b\\c+1-a-b\end{matrix};1-z  \Bigr]
\F\Bigl[\begin{matrix}c'-a',c'-b'\\c'+1-a'-b'\end{matrix};
1-\ov z  \Bigr],
\label{eq:B}
\end{multline}
where
\begin{align}
B_{0}=
\frac{\Gamma^\C(\bfc)\, \Gamma^\C(\bfc-\bfa-\bfb)}{\Gamma^\C(\bfc-\bfa)\Gamma^\C(\bfc-\bfb)}
,
\label{eq:B00}
\\
B_{1}=\frac{\Gamma^\C(\bfc)\Gamma^\C(\bfa+\bfb-\bfc)}{\Gamma^\C(\bfa)\, \Gamma^\C(\bfb)}.
\label{eq:B11}
\end{align}

{\rm c)} In the disk $|z|>1$
the following expansion holds:
\begin{multline}
\FF\Bigl[\begin{matrix}\bfa;\bfb\\ \bfc\end{matrix};z\Bigr]=
C_{0}\cdot (-z)^{-a|-a'}
\F\Bigl[\begin{matrix}a,a+1-c\\a+1-b\end{matrix}; z^{-1}  \Bigr]
\F\Bigl[\begin{matrix}a',a'+1-c'\\a'+1-b'\end{matrix};
 \ov z^{\,\,-1}  \Bigr]
 +\\+
 C_{1}\cdot
(-z)^{-b|-b'}
\F\Bigl[\begin{matrix}b,b+1-c\\b+1-a\end{matrix}; z^{-1}  \Bigr]
\F\Bigl[\begin{matrix}b',b'+1-c'\\b'+1-a'\end{matrix};
 \ov z^{\,\,-1}  \Bigr],
\label{eq:C}
\end{multline}
where
\begin{align}
C_{0}:=\frac{\Gamma^\C(\bfc)\,\Gamma^\C(\bfb-\bfa) }{\Gamma^\C(\bfc-\bfa)\,\Gamma^\C(\bfb)}
,
\label{eq:C00}
\\
C_{1}:=\frac{\Gamma^\C(\bfc)\,\Gamma^\C(\bfa-\bfb)}{\Gamma^\C(\bfc-\bfb)\,\Gamma^\C(\bfa)}
\label{eq:C11}
.
\end{align}
\end{theorem}

\sm


 \smallskip

  {\bf\punct Proof of Theorem \ref{th:long}.%
 \label{ss:asy}}
Forms (\ref{eq:A}), (\ref{eq:B}), (\ref{eq:C}) 
for the desired expressions follow from the preceding considerations.
Also we know that the coefficients $A_{0}$, $A_{1}$, $B_{0}$, $B_{1}$,
$C_{0}$, $C_{1}$ are meromorphic in $\bfa$, $\bfb$, $\bfc$.
Now we apply asymptotic expansions from Theorem \ref{th:weak-singularities}.

\sm

 1. {\it Asymptotic of $\FF[\bfa,\bfb;\bfc;z]$
 as $z\to\infty$.}
 Assume that the defining integral for $\FF[\bfa,\bfb;\bfc;z]$
 converges, and also 
 \begin{equation}
 [\bfb]-[\bfa]>0
 \label{eq:condition}
 \end{equation}
 Then
\begin{multline*}
\frac1{\pi B^\C(\bfb,\bfc-\bfb)}
\int_\C t^{\bfb\bf-1} (1-t)^{\bfc\bf-\bfb-1}
(1-zt)^{-\bfa}\,\dd t
=\\=
\frac{z^{\bf-\bfa}}{\pi B^\C(\bfb,\bfc-\bfb)}
\int_\C t^{\bfb\bf-1} (1-t)^{\bfc\bf-\bfb-1}
(z^{-1}-t)^{-\bfa}\,\dd t
\sim \\ \sim
\frac{B^\C(\bfb-\bfa,\bfc-\bfb)}
{B^\C(\bfb,\bfc-\bfb)}\cdot (-z)^{-\bfa}\cdot
\Bigl(1+\sum_{(i,i')\ne (0,0)} p_{ii'}z^{-i|-i'}\Bigr)
+\\+
\frac{B^\C(\bfb,\bf 1-\bfa)}
{B^\C(\bfb,\bfc-\bfb)}\cdot z^{-\bfb}\cdot
\Bigl(1+\sum_{(i,i')\ne (0,0)} q_{ii'}z^{-i|-i'}\Bigr).
\end{multline*}

 Precisely, denote $z^{-1}$ by $\epsilon$, and
 denote the integrand in the last integral by
 $H(\cdot)$. Let $\phi(t)\ge 0$, $\psi(t)\ge 0$ be smooth functions
 such that $\phi(t)+\psi(t)=1$, $\phi(t)=1$ near 0,
 and $\psi(t)=1$ near $\infty$. A straightforward differentiation with respect to the parameter $\epsilon$
 shows that
 $$
\int H(t;\epsilon) \psi(t)
\dd{t}
 $$
is smooth  near $\epsilon=0$. For
$$
\int_\C H(t;z) \phi(t;\epsilon)
\dd{t}
$$
we apply Theorem \ref{eq:weak-singularities},
due to the restriction (\ref{eq:condition}) we can also apply (\ref{eq:r00}).
Thus we get explicit coefficients
$C_{0}$, $C_{1}$ in the expansion (\ref{eq:C}).
To remove restrictions for the parameters, we refer to the analytic
continuation.

Finally, we transform $B^\C(\bfb,\bf 1-\bfa)$ with formula (\ref{eq:complement}),
$$
B^\C(\bfb,{\bf 1}-\bfa)=\frac{\Gamma^\C(\bfb)\,\Gamma^\C(\1-\bfa)}{\Gamma^\C(\1+\bfb-\bfa)}
=(-1)^{\bfb}\frac{\Gamma^\C(\bfb)\Gamma^\C(\bfa-\bfb)}{\Gamma^\C(\bfa)}.
$$

2. {\it  Asymptotic of $\FF[\bfa,\bfb;\bfc;z]$ as $z\to 0$.}
Substituting $t=1/s$ to the definition
(\ref{eq:def-FF}) of $\FF$, we
get
\begin{multline*}
\FF[\bfa,\bfb;\bfc;z]
=
\frac{(-1)^{\bfc-\bfa-\bfb}}{\pi B^\C(\bfb,\bfc-\bfb)}
\int_\C
s^{-\bfc+\bfa}(z-s)^{-\bfa}(1-s)^{\bf\bfc-\bfb-1}
\,\dd s
\sim\\ \sim
\frac{(-1)^{\bfc-\bfa-\bfb} B^\C(\bfa-\bfc\bf+1, \bf 1-\bfa)}
{B^\C(\bfb,\bfc-\bfb)}\cdot
z^{\1-\bfc}
\cdot \Bigl(1+\sum_{(i,i')\ne (0,0)}
p_{ii'} z^{i|i'}\Bigr)
+\\+
\frac{(-1)^{\bfc-\bfb} B^\C(\bf 1-\bfc, \bfc-\bfb)}
{B^\C(\bfb,\bfc-\bfb)}
\cdot \Bigl(1+\sum_{(i,i')\ne (0,0)}
q_{ii'} z^{i|i'}\Bigr).
\end{multline*}

\smallskip

3. {\it
Asymptotic of $\FF[\bfa,\bfb;\bfc;z]$ as $z\to 1$.}
We substitute $t=\frac 1{1-s}$ to (\ref{eq:def-FF})
and get
\begin{multline*}
\FF[\bfa,\bfb;\bfc;z]=
\frac {(-1)^{\bfc-\bfb}}
{\pi B^\C(\bfb,\bfc-\bfb)}
\int_\C
s^{\bf c-b-1}
(1-s)^{\bfa-\bfc}
(1-z-s)^{-\bfa} \,\dd s
\sim\\ \sim
\frac {(-1)^{\bfc-\bfb}B^\C(\bfc-\bfb, \bf 1-\bfa)}
{ B^\C(\bfb,\bfc-\bfb)}
\cdot (1-z)^{\bf c-b-a}\cdot
\Bigl(1+\sum_{(i,i')\ne (0,0)}
p_{ii'} (1-z)^{i|i'} \Bigr)
+\\+
\frac {(-1)^{\bfc-\bfb-\bfa}B^\C(\bfc-\bfb-\bfa, \bf 1+\bfa-\bfc)}
{ B^\C(\bfb,\bfc-\bfb)}
\cdot
\Bigl(1+\sum_{(j,j')\ne (0,0)}
q_{ii'} z^{j|j'} \Bigr).
\end{multline*}


{\sc Remark. Another way of a proof of Theorem \ref{th:long}.}
Applying  the Kummer formulas,
	Erd\'elyi, et al., \cite{HTF1}, Section 2.9, we can
	write the analytic continuation of (\ref{eq:AB}) to a neighborhood
	of this point. 
	The resulting expression for 
	$\FF$ 
	must be non-ramified at $z=1$.
	This gives us the
	 coefficients in (\ref{eq:AB}) 
	up to a common  factor. In fact this calculation is done 
	below in the proof of Proposition \ref{pr:uniqueness}. 
	The scalar factor can be evaluated
	using (\ref{eq:Gauss-identity}). It remains to apply 
	the Kummer formulas (\cite{HTF1}, Section 2.9) for the analytic continuation again  and to get
	an expansion at $\infty$.
	\hfill $\boxtimes$

\smallskip

{\bf\punct
Additional symmetry.}

\begin{proposition}
 Let $a-b\in\Z$. Then
 \begin{equation}
  \FF\left[\begin{matrix}a|a',b|b'\\
  c|c'\end{matrix};z\right]=
    \FF\left[\begin{matrix}a|b',b|a'\\
  c|c'\end{matrix};z\right].
  \label{eq:additional}
 \end{equation}
\end{proposition}

{\sc Proof.} The expansions  (\ref{eq:A})--(\ref{eq:A11}) at 0 for  both  functions
are identical. We only must verify the equality of the denominators in
(\ref{eq:A11}):
\begin{multline}
 \Gamma^\C(a|a')\,\Gamma^\C(b|b')\,\Gamma^\C(c-a|c'-a')\,\Gamma^\C(c-b|c'-b')=\\
  \Gamma^\C(a|b')\,\Gamma^\C(b|a')\,\Gamma^\C(c-a|c'-b')\,\Gamma^\C(c-b|c'-a').
  \label{eq:R-additional}
\end{multline}
The both sides are equal to
$$
\frac{(-1)^{c'-c} \pi^4 \Gamma(a)\Gamma(a')\Gamma(b)\Gamma(b')\Gamma(c-a)\Gamma(c-a')\Gamma(c-b)\Gamma(c-b')}
{\sin\pi a' \sin\pi b' \sin\pi (c'-a') \sin\pi (c'-b') }.
 \qquad\qquad \square
$$

\sm

{\bf\punct Degenerations and logarithmic expressions.}

\sm

a) {\sc Residues and zeros.}
Notice that poles and zeros of $\FF[\bfa,\bfb;\bfc;z]$ 
as  function of $\bfa$, $\bfb$, $\bfc$
depend on a choice of a normalizing factor in the front
of the integral (\ref{eq:def}).

It is easy to see that residues at poles also are solutions of the complex hypergeometric system
(\ref{eq:two-equations}). 
The expressions for the residues can be obtained from our expansions.

For obtaining the residues at $\{a|a'\}\in \N\times \N$ we can use the expansion of $\FF$ at $z=0$,
see (\ref{eq:A})--(\ref{eq:A11}). We get 
$$
z^{1-c|1-c'}
\F\Bigl[\begin{matrix}a+1-c,b+1-c\\2-c  \end{matrix}; z  \Bigr]\,
\F\Bigl[\begin{matrix}a'+1-c',b'+1-c'\\2-c'  \end{matrix};
 \ov z  \Bigr]
$$
with an obvious $\Gamma^\C$-factor. Applying the Pfaff transformations 
of $\F$, we  observe that these expressions  are elementary functions.
Formulas (\ref{eq:A})--(\ref{eq:A11}) allow to calculate residues at the poles of all the types (\ref{eq:poles}).

\sm

Next, consider another  normalization%
\footnote{In fact, in the main part of our work we use this normalization of the kernel,
see (\ref{eq:cK}). Due to this we do not lose the case of $L^2$
on the complex quadric discussed in Subsect. \ref{ss:radial}.} of the functions $\FF$:
\begin{equation}
_2^{\vphantom{\C}}\wt F^\C_1[\bfa,\bfb;\bfc;z]:=\frac 1{\Gamma^\C(\bfc)}\,\,\FF[\bfa,\bfb;\bfc;z].
\label{eq:wt}
\end{equation}
This operation cancels the factor $\Gamma^\C(\bfc)$
in expansion of $\FF[z]$ at $\infty$, see (\ref{eq:C})--(\ref{eq:C11}).
So we get a finite expression at the poles (\ref{eq:poles-1}) and non-zero
function at the zeros (\ref{eq:zeros-F}).

Thus, at all exceptional planes (\ref{eq:poles})--(\ref{eq:zeros-F})
we get explicit nonzero expressions. Such expressions  also
depend on normalization of $\FF[\dots,z]$, but for a point $(\bfa_0,\bfb_0,\bfc_0)$
being in a general position on an exceptional plane such nonzero expression 
is canonically defined up to a constant factor.

\sm

{\sc b) Further degenerations.} Classical hypergeometric differential equation has a sophisticated
list of degenerations, see  \cite{HTF1}, Sect. 2.2. In our case new difficulties arise
if at least two of the parameters $\bfa$, $\bfb$, $\bfc-\bfa$, $\bfc-\bfb$ are contained in $\Z\times \Z$. 
We stop here further analysis and only notice that {\it for exceptional values of the parameters a solution of
the complex hypergeometric system {\rm (\ref{eq:two-equations})} can be non-unique}.

For instance, if $\bfa\in\Z_-\times \Z_-$, $\bfc-\bfb\in \N\times \N$, then both  summands in
(\ref{eq:AB}) are single-valued (since all hypergeometric series are terminating).

\sm

{\sc c) Logarithmic expressions.} For definiteness we discuss the case 
$$\bfc\in \N\times \N$$
(which is interesting for our further purposes). Consider the function $_2^{\vphantom{\C}}\wt F^\C_1$
defined by (\ref{eq:wt}). It has a removable singularity at our $\bfc$.
Recall that for $c=n\in\N$ the usual hypergeometric differential equation
$D[a,b,n] f=0$ has two solutions. The first is $\F[a,b;c;z]$
and the second has the form
\begin{equation}
\Psi[a,b;n;z]=\sum_{j=-n+1}^\infty p_j z^j + \ln z\cdot \F[a,b;n;z],
\label{eq:logarithmic-solution}
\end{equation}
where $p_j$ are explicit coefficients, $p_{-n+1}\ne 0$, and this form does not depend on further degenerations,
see \cite{AAR}, Section 2.3. Passing around 0 we get
$$
\Psi[a,b;n;e^{i\phi}z]\Bigr|_{\phi=2\pi}-\Psi[a,b;n;z]=\F[a,b;n;z].
$$
Thus the system
$$D[a,b;n] \cF=0, \qquad D[a',b';n']\cF=0$$
has two solutions that are single-valued near zero,
the first is obvious
$$
\F[a,b,n;z]\, \F[a',b';n';\ov z],
$$
and the second is
\begin{equation}
\F[a,b;n;z]\, \Psi[a',b';n';\ov z]+ \Psi[a,b;n;z]\,  \F[a',b';n';\ov z].
\label{eq:log}
\end{equation}
Our function $_2^{\vphantom{\C}}\wt F^\C_1[\bfa,\bfb;\bfn;z]$  is certain linear combination of these solutions. 

\smallskip

{\sc d) On uniqueness of a solution of the hypergeometric system.}

\begin{proposition}
	\label{pr:uniqueness}
	Let 
	\begin{align*}
\text{$a$, $b$, $c$, $c-a-b$, $c-a$, $c-b \notin \Z$},
\\
\text{$a'$, $b'$, $c'$, $c'-a'-b'$, $c'-a'$, $c'-b' \notin \Z$}.
	\end{align*}
	Let the system $D[a,b,c]\,\cF=0$, $D'[a',b', c']\,\cF=0$ 
	have a non-ramified non-zero solution. Then $c-c'\in \Z$
	and
	\begin{equation}
	\text{$a-a'$, $b-b'\in \Z$ or $a-b'$, $b-a'\in \Z$}
	\label{eq:abprime}
	\end{equation}
	Such solution is unique up to a scalar factor and therefore is 
	$\FF[a|a',b|b';c|c';z]$ or $\FF[a|b',b|a';c|c';z]$
\end{proposition}


{\sc Proof.} 
First, we examine the behavior of a solution near $z=0$.
Let
$$
\phi(z):=\F\left[\begin{matrix}
a,b\\c
\end{matrix};z\right], \qquad
\psi(z):= \F\left[\begin{matrix}
a+1-c,b+1-c\\2-c
\end{matrix};z\right],
$$
i.e., $\phi$, $z^{1-c}\psi$
are the Kummer solutions of the hypergeometric equation $D[a,b,c]f=0$
 at $0$, see \cite{HTF1}, (2.9.1), (2.9.17). Denote by $\aphi(\ov z)$,
 $\apsi(\ov z)$ the similar functions obtained be the change
 $a\mapsto a'$, $b\mapsto b'$, $c\mapsto c'$, $z\mapsto \ov z$.
A solution of our system near 0 has the form
$$
G(z)=
\sigma\, \phi(z)\aphi(\ov z)+ \mu\, \ov z^{1-\ov c'} \phi(z)\apsi(\ov z) + \nu\, z^{1-c} \psi(z)\aphi(\ov z)
+\tau\, z^{1-c|1-\ov c^{\,\prime}} \psi(z)\apsi(\ov z).
$$
Passing $m$ times around 0 we come to
\begin{multline*}
G(z)=
\sigma\, \phi(z)\acute{\phi}(\ov z)+ \mu e^{2\pi mc'i}\, \ov z^{1-\ov c'} \phi(z)\apsi(\ov z)
+\\ + 
\nu e^{-2\pi mc i}\,z^{1-c} \psi(z)\aphi(\ov z)
+\tau e^{2\pi m(c'-c)i}\, z^{1-c|1-\ov c^{\,\prime}} \psi(z)\apsi(\ov z).
\end{multline*}
Since $c$, $c'\notin \Z$, we have $ e^{2\pi mc'i}$,  $e^{-2\pi mci}\ne 1$,
on the other hand they are   $\ne e^{2\pi m(c'-c)i}$. If $G(z)$ is single-valued,  
then
 $\mu=\nu=0$. Also, we need $\tau=0$ or $c-c'\in \Z$.

\sm

 To examine the behavior of $G$ near  $z=1$ 
we apply a formula for analytic continuation,
see \cite{HTF1}, Subsect. 2.10. Near $z=1$ we have
\begin{multline}
F\left[\begin{matrix}
a,b\\c
\end{matrix};z\right]=
A_1(a,b,c)\,
 \F\left[\begin{matrix}
a,b\\ a+b-c+1
\end{matrix}; 1-z
\right]+
\\
+ A_2(a,b,c)\,(1-z)^{c-a-b} 
\F\left[\begin{matrix}
c-a,c-b\\ c-a-b+1
\end{matrix}; 1-z
\right]
\label{eq:analytic-cont}
,
\end{multline}
where
\begin{equation}
A_1(a,b,c):=\frac{ \Gamma(c) \Gamma(c-a-b)}{\Gamma(c-a)\Gamma(c-b)},
\qquad
A_2(a,b,c):=
\frac{\Gamma(c)\Gamma(a+b-c)}{\Gamma(a)\Gamma(b)}.
\label{eq:analytic-cont1}
\end{equation}

Since $c-a-b$, $c'-a'-b'\notin \Z$, the expression $\phi(z)\aphi(\ov z)$
is not single valued. Thus $\tau\ne 0$, $c-c'\in\Z$, and
$$
G(z)=\sigma\, \phi(z)\aphi(\ov z)+ 
\tau\, z^{1-c|1-\ov c^{\,\prime}} \psi(z)\apsi(\ov z).
$$
Applying for $\phi$, $\aphi$, $\psi$, $\apsi$
formula (\ref{eq:analytic-cont}) and the identity
$$
\F(\alpha,\beta;\gamma;z)=(1-z)^{\gamma-\alpha-\beta}
\F(\gamma-\alpha,\gamma-\beta;\gamma;z),
$$
 we come to
 {\small
\begin{multline*}
G(z)=\sigma A(a,b,c)A(a',b',c')\,\F\left[\begin{matrix}
a,b\\c
\end{matrix};1-z\right] \F\left[\begin{matrix}
a',b'\\c'
\end{matrix};1-\ov z\right]+\\+
\Bigl\{\sigma A_1(a,b,c) A_2 (a',b',c')+\tau
A_1 (a+1-c,b+1-c,2-c)A_2(a'+1-c',b'+1-c',2-c')
\Bigr\}\times\\
\times (1-\ov z)^{c'-a'-b'} \F\left[\begin{matrix}
a,b\\c
\end{matrix};1-z\right] 
\F\left[\begin{matrix}
c'-a',c'-b'\\
c'-a'-b'+1
\end{matrix};1-\ov z\right]
+
\\+
\Bigl\{\sigma A_2(a,b,c) A_1 (a',b',c')+\tau
A_2 (a+1-c,b+1-c,2-c)A_1(a'+1-c',b'+1-c',2-c')
\Bigr\}\times\\
\times (1-z)^{c-a-b}
\F\left[\begin{matrix}
c-a,c-b\\
c-a-b+1
\end{matrix};1- z\right]
\F\left[\begin{matrix}
a',b'\\c'
\end{matrix};1-\ov z\right]
+\\+
A_2(a+1-c,b+1-c,2-c)A_2(a'+1-c',b'+1-c',2-c')
\times\\\times
(1-z)^{c-a-b|c'-a'-b'}
\F\left[\begin{matrix}
c-a,c-b\\
c-a-b+1
\end{matrix};1- z\right]
\F\left[\begin{matrix}
c'-a',c'-b'\\
c'-a'-b'+1
\end{matrix};1- \ov z\right]
\end{multline*} 
}
The coefficients $A_1(\cdot)$, $A_2(\cdot)$ have no zeros and no
poles under our restrictions.
The expression is single-valued if and if two curly 
brackets are zero and $(c-a-b)-(c'-a'-b')\in \Z$.
This implies 
$$
(a+b)-(a'+b')\in \Z.
$$

Two curly brackets give a system of linear equations 
for $\sigma$, $\tau$. It has a nonzero solution if and only if its determinant
$\Delta$ is zero. Straightforward calculations
give
\begin{multline*}
\!\!\!
\Delta=
\pi^{-4}\Gamma(c)\Gamma(c')\Gamma(2-c)\Gamma(2-c')\Gamma(c-a-b)\Gamma(c'-a'-b')\Gamma(a+b-c)
\Gamma(a'+b'-c')\cdot \Xi,
\end{multline*}
where 
$$
\Xi=\sin \pi(c-a)\sin\pi(c-b)\sin \pi a'\sin\pi b'-
\sin \pi(c'-a')\sin\pi(c'-b')\sin \pi a\sin\pi b.
$$
Clearly, the set $\Xi(a,b,c,a',b',c')=0$ is invariant with respect to the shifts
$a\mapsto a+1$, $b\mapsto b+1$, $c\mapsto c+1$.
Therefore to examine the set of zeros we can assume
$c'=c$, $b'=a+b-a'$. Under these conditions $\Xi$ can be reduced to the following
form:
$$
\Xi(a,b,c,a',b',c')=\sin\pi(a-a')\sin\pi (a'-b)\sin\pi c\sin\pi(c-a-b)
$$
(this non-obvious identity can be verified by decompositions of both sides into sums
if exponentials).
This implies (\ref{eq:abprime}). 

If $\Delta=0$ then $\sigma$, $\tau$ are defined up to a common
scalar factor, this proves the uniqueness
(and gives an expression for $\sigma/\tau$).
\hfill $\square$

\sm

{\sc e) Non-interesting solutions.} However, we have seen that the complex  hypergeometric  system
for some values of the parameters has two single-valued solutions. Also, there are solutions that do not seem reasonable.
For instance, we have
$$
D[0,b_1,c_1] \cdot 1=0,\qquad D'[0,b_2,c_2] \cdot 1=0
$$
for arbitrary $b_1$, $c_1$, $b_2$, $c_2\in\C$.

\sm

{\bf\punct Differential-difference equations for $\FF$.}
We can regard 
$\F[a,b;c;z]$ as a family of functions of a complex variable
$z$ depending on 3 parameters $a$, $b$, $c$. But we also can regard 
$\F[a,b;c;z]$ as {\it one function} of the four complex variables $a$, $b$, $c$, $z$.
Then $\F[a,b;c;z]$ satisfy a non-obvious system of linear
differential-difference equations, some examples of such equations are 
in Erd\'elyi, et al., \cite{HTF1}, (2.8.20-45). Below we show that such equations
can be automatically transformed to differential-difference equations 
for the function $\FF[a|a',b|b';c|c';z]$ of 7 complex variables.

Consider a space of functions in the variables $a$, $b$, $c$, $z$. 
Define operators
\begin{align*}
T_a f(a,b,c,z)&= f(a+1,b,c,z) ,\qquad T_b f(a,b,c,z)=f(a,b+1,c,z),
\\
T_c f(a,b,c,z)&=f(a,b,c+1,z).
\end{align*}
Consider {\it finite} sums of the form
\begin{equation}
\cL=
\sum_{j\ge 0} \sum_{k,l,m\in \Z} U_{j,k,l,m}(a,b,c,z)\, T_a^k\, T_b^l \, T_c^m \frac{\partial^j}{\partial z^j},
\label{eq:cL}
\end{equation}
where $U_{j,k,l,m}(a,b,c,z)$ are polynomial expressions in $z$ with 
 coefficients rationally depending on $a$, $b$, $c$.

Assume that
$$
\cL \,\, \F[a,b;c;z]=0.
$$
We can regard an operator (\ref{eq:cL}) as an operator on functions
$f(a|a',b|b',c|c',z)$ on $\Lambda^3\times \dot\C$. We also define operators
\begin{align*}
T_{a'} f(a|a',\,b|b',\,c|c',\,z)= f(a|a'+1,\,b|b',\,c|c',\,z) ,
\\T_{b'} f(a|a',\,b|b',\,c|c',\,z)=f(a|a',\,b|b'+1,\,c|c',\,z),
\\
T_{c'} f(a|a',\,b|b',\,c|c',\,z)=f(a|a',\,b|b',\,c|c'+1,\,z).
\end{align*}
For such an operator $\cL$ we define the operator $\cL'$ by 
\begin{equation*}
\cL'=
\sum_{j\ge 0} \sum_{k,l,m\in \Z} U_{j,k,l,m}(a',b',c',\ov z) \,
T_{a'}^k\, T_{b'}^l\, T_{c'}^m \frac{\partial^j}{\partial \ov z^j}.
\end{equation*}
From  the definition it follows that
$$
\cL \cL'= \cL' \cL.
$$

\begin{proposition}
 \label{pr:contiguous}
 Let the function $Q(a,b,c,z)=\F[a,b;c;z]$ satisfy an equation $\cL\,Q=0$.
 Then the function 
 $$R(a|a',\,b|b',\,c|c',\,z):=\FF[a|a',\,b|b',\,c|c',\,z]$$
 satisfies the system of equations
 \begin{equation}
 \cL\, R(a|a',\,b|b',\,c|c',\,z)=0,\qquad \cL'\, R(a|a',\,b|b',\,c|c',\,z)=0.
 \label{eq:cLcL}
 \end{equation}
\end{proposition}

\begin{lemma}
\label{l:for-difference}
Let $Q=\F[a,b;c;z]$ satisfy an equation
$
\cL Q=0
$. Then 
\begin{equation}
\frac{e^{\pi i(c-a-b)}\Gamma(c)\Gamma(c-1)}{\Gamma(a)\Gamma(b)\Gamma(c-a)\Gamma(c-b)}
z^{1-c}\F[a+1-c,b+1-c;2-c,z]
\label{eq:difference-analytic}
\end{equation}
satisfies the same equation.
\end{lemma}

{\sc Remark.} The same statement holds for the
functions
\begin{align}
u_1&=\frac{ \Gamma(c-a) \Gamma(c-a-b)}{\Gamma(c)\Gamma(c-b)} \F\left[\begin{matrix}
                                                               a,b\\ a+b-c+1
                                                              \end{matrix}; 1-z
 \right];
 \label{e:u1}
 \\
 u_2&=\frac{\Gamma(c)\Gamma(a+b-c)}{\Gamma(a)\Gamma(b)} (1-z)^{c-a-b} 
 \F\left[\begin{matrix}
                                                               c-a,c-b\\ c-a-b+1
                                                              \end{matrix}; 1-z
 \right];
 \label{e:u2}
 \\
 u_3&=\frac{\Gamma(c)\Gamma(b-a)}{\Gamma(b)\Gamma(c-b)} z^{-a} 
  \F\left[\begin{matrix}
                                                              a, 1-c+a\\ 1-b+a 
                                                              \end{matrix}; z^{-1} \right]
                                                              \nonumber
\end{align}
and also for other summands in the right-hand sides of formulas Erd\'elyi, et al. \cite{HTF1},
(2.10.1)--(2.10.4). \hfill $\square$

\sm

{\sc Proof of Lemma \ref{l:for-difference}.}
First, let $a$, $b$,  $c$ be in a general position.
By Erd\'elyi, et al. \cite{HTF1}, (2.10.1), (2.10.5),
\begin{equation}
F(a,b;c;z)=u_1+u_2
,
\label{eq:u1u2}
\end{equation}
where $u_1$, $u_2$ are given by (\ref{e:u1})--(\ref{e:u2}).
The function $u_2$ is ramified at $z=1$. Passing around this point we get
a function
$$
\wt F:=
u_1+ e^{2\pi i(c-a-b)} u_2.
$$
By  analytic continuation, $\cL\wt F=0$.
The factor $e^{2\pi i(c-a-b)} $  does not change under the shifts
 $T_a$, $T_b$, $T_c$.
Therefore the summands  $u_1$, $u_2$ satisfy the same equation, $\cL u_1=0$, 
$\cL u_2=0$.
We apply the same transformation (\ref{eq:u1u2}) to the summand $u_1$ and repeat the same reasoning.
We observe that 
$$
\frac{\pi \Gamma(c)\Gamma(c-1)}{\Gamma(a)\Gamma(b)\Gamma(c-a)\Gamma(c-b)
\sin \pi(a+b-c)}
z^{1-c}\F[a+1-c,b+1-c;2-c,z]
$$
satisfies the same equation. This expression differs from (\ref{eq:u1u2}) by
the factor
$ e^{i\pi (a+b-c)} \sin \pi(a+b-c)$,
which is invariant under the shifts $T_a$, $T_b$, $T_c$. 

Passing to a limit we omit restrictions to $a$, $b$, $c$.
\hfill $\square$

\sm

{\sc Proof of Proposition \ref{pr:contiguous}.} We use the expression (\ref{eq:A})
for $\FF[\bfa,\bfb;\bfc;z]$. Obviously, the first summand satisfies
the  system (\ref{eq:cLcL}). By Lemma \ref{l:for-difference}, the expression
\begin{multline*}
 \frac{e^{\pi i(c-a-b)}\Gamma(c)\Gamma(c-1)}{\Gamma(a)\Gamma(b)\Gamma(c-a)\Gamma(c-b)}
z^{1-c}\F\left[\begin{matrix}a+1-c,b+1-c\\2-c\end{matrix};z\right]
\times\\\times
 \frac{e^{\pi i(c'-a'-b')}\Gamma(c')\Gamma(c'-1)}{\Gamma(a')\Gamma(b')\Gamma(c'-a')\Gamma(c'-b')}
\ov z^{1-c'}\F\left[\begin{matrix}a'+1-c',b'+1'-c\\2-c'\end{matrix};\ov z\right].
\end{multline*}
satisfies the system (\ref{eq:cLcL}).
It differs from the second summand in (\ref{eq:A})
by a  factor
$$
\frac {i^0 \sin \pi a' \sin \pi b' \sin \pi (c'-a') \sin \pi (c'-b')}
{\sin \pi c' \sin \pi (c'-1)}.
$$
This expression is invariant with respect to shifts $T_a$, $T_{a'}$, \dots. 
Therefore the second summand
in (\ref{eq:A}) also satisfies the system.
\hfill $\square$

\sm

{\bf\punct One difference operator.}
By \cite{Ner-index}, formula (2.3),
the
Gauss hypergeometric function
$
\F[p,q;r;z]
$
satisfies
the following difference equation
\begin{multline}
-z\cdot  \F(p,q;r;z)=\\
=\frac{q(r-p)}{(q-p)(1+q-p)}\,\F(p-1,q+1;r;z)-\\
-\Bigl[\frac{q(r-p)}{(q-p)(1+q-p)}+
 \frac{p(r-q)}{(p-q)(1+p-q)}\Bigr]\F(p,q;r;z)+\\
+\frac{p(r-q)}{(p-q)(1+p-q)}\,\F(p+1,q-1;r;z).
\label{eq:difference-old}
\end{multline}

Define the  difference  operators $L$, $L'$
acting on functions of the variables $\bfa$, $\bfb$, $\bfc$, $z$ by
\begin{equation}
L =
\frac{b(c-a)}{(b-a)(1+b-a)}
\bigl(T_a^{-1}\,T_b- 1\bigr)
-\frac{a(c-b)}{(a-b)(1+a-b)}
\bigl(T_a T_b^{-1}- 1\bigr);
\end{equation}
\begin{equation}
L'=
\frac{b'(c'-a')}{(b'-a')(1+b'-a')}
\bigl(T_{a'}^{-1} T_{b'}- 1\bigr)
-\frac{a'(c'-b')}{(a'-b')(1+a'-b')}
\bigl(T_{a'} T_{b'}^{-1}- 1\bigr).
\end{equation}


\begin{corollary}
\label{cor:difference-eq}
The complex hypergeometric function $\FF[\bfa,\bfb;\bfc;z]$
 satisfies the following system of difference equations
 \begin{align}
L\,\,\FF[\bfa,\bfb;\bfc;z]=z\,\, \FF[\bfa,\bfb;\bfc;z];
\label{eq:L}
\\
L'\,\,\FF[\bfa,\bfb;\bfc;z]=\ov z\,\, \FF[\bfa,\bfb;\bfc;z].
\label{eq:L-prime}
 \end{align}
\end{corollary}

{\bf \punct Some properties of the kernel $\cK$.} 
We have the following corollaries from our previous considerations.

\sm

1) By (\ref{eq:ab-ba}) {\it  $\cK_{a,b}$ is even},
\begin{equation}
 \cK_{a,b}(z;-k,-\sigma)=\cK_{a,b}(z;k,\sigma).
\end{equation}

2) By (\ref{eq:f-conjugate}),
\begin{equation}
 \ov{\cK_{a,b}(z;k,-\ov \sigma)}=\cK_{a,b}(z;k,\sigma).
\end{equation}
In particular, {\it $\cK_{a,b}(z;k,\sigma)$ is real on $\Lambda$}.

\sm

3) By Proposition \ref{pr:system-hypergeometric}, {\it $\cK_{a,b}(z;k,\sigma)$ satisfies the following differential equations:}
\begin{align}
 \frD \, \cK_{a,b}(z;k,\sigma)&=\tfrac 14(k+\sigma)^2 \, \cK_{a,b}(z;k,\sigma);
 \\
  \ov\frD\, \cK_{a,b}(z;k,\sigma)&=\tfrac 14(k-\sigma)^2 \, \cK_{a,b}(z;k,\sigma).
\end{align}

4) By Corollary \ref{cor:difference-eq}, {\it $\cK_{a,b}(z;k,\sigma)$ satisfies the following difference equations:}
\begin{align}
\frL\, \cK_{a,b}(z;k,\sigma)&= z\, \cK_{a,b}(z;k,\sigma);
\\
\ov\frL \,\cK_{a,b}(z;k,\sigma)&= \ov z \, \cK_{a,b}(z;k,\sigma).
\end{align}

\section{Nonexistence of commuting self-adjoint extensions\label{s:differential}}

\COUNTERS

Here we prove that for $(a,b)\notin \Pi$ 
the operators $\tfrac12(\frD+\ov\frD)$, $\tfrac1{2i}(\frD-\ov\frD)$
defined on $\cD(\dot\C)$ do not admit commuting self-adjoint extensions.
We analyze the set of possible generalized eigenfunctions and show that this set is 
too small.

\sm

{\bf \punct Generalized eigenfunctions.}
Denote by $\cD'(\dot\C)$ the  space of distributions on $\dot\C$.
We have a nuclear rigging (see \cite{BSU}, Section 14.2)
$$
\cD(\dot\C)\subset L^2( \C,\mu_{a,b})\subset \cD'(\dot\C),
$$
and  apply the usual formalism of generalized eigenfunctions,
see  \cite{BSU}, Chapter 15.

Recall that we have formally symmetric and formally commuting operators
$$
D_+:=
\tfrac 12(\frD+\ov\frD),\qquad
D_-:=\tfrac 1{2i}(\frD-\ov\frD)
$$
in $L^2(\C,\mu)$ (defined on the domain $\cD(\dot\C)$) and 
the spectral problem
\begin{equation}
\frD \Phi=\zeta \Phi,
\qquad
\ov \frD \Phi=\ov \zeta \Phi.
\label{eq:DDPhi}
\end{equation}

{\it Suppose that the operators $D_+$,
$D_-$ admit commuting self-adjoint extensions.}
Then the operator $U$ of  spectral decomposition can be written in terms of
generalized eigenfunctions. Precisely, there exist a  space
$R$ equipped with a measure $\rho$ and 
an injective  measurable map $r\mapsto \phi_r$
from $R$ to 
$\cD'(\dot\C)$ such that  
$$
D_+\phi_{r}=a(r) \phi_{r}, \qquad D_-\phi_r=b(r) \phi_r,
$$
where $a(r)$, $b(r)$ are real-valued functions,
and the pairing
$$
U f(r)=\{f, \phi_{r}\}
$$
 of $f\in \cD(\dot\C)$ and $\phi_r$
 determines a unitary operator
 $L^2(\C, \mu)\to L^2(R,\rho)$, see textbook \cite{BSU},
 Subsect. 15.2.3%
 \footnote{Basically, this is a result of Kostyuchenko and Mityagin 
 \cite{KM1}--\cite{KM2} with weaker conditions for a rigging.}%
).

Since the operator $\frD$ is elliptic, generalized eigenfunctions
are smooth  on $\dot\C$, see e.g.,
\cite{BSU}, Theorem 16.2.1. 
 Therefore 
in our case
 generalized eigenfunctions $\phi_{r}$ are usual smooth
solutions of the system of differential equations.

We also can identify the measure space $R$ with its image, and so we can think that {\it the measure $\rho$
is sitting on the space $\Omega$ of smooth solutions of the
systems} (\ref{eq:DDPhi}),
{\it where $\zeta$ ranges in $\C$}.  We intend to show that for any
measure $\rho$ on $\Omega$ the operator
$J:L^2(\Omega,\rho)\to L^2(\C,\mu_{a,b})$
defined by
$$
U h (z)=\int_{\Omega} h(r) \phi_r(z)\,d\rho(r)  
$$
is not unitary. Precisely:

\begin{lemma}
\label{l:atomic}
Let $(a,b)\notin \Pi$.
 Let $\rho$ be a measure on $\Omega$, and let the corresponding operator $U$ be bounded.
 Then $\rho$ is an atomic measure supported by a finite set.
\end{lemma}

The idea of a proof is simple, it is explained in the next subsection,
a formal proof is completed in Subsect. \ref{ss:non-self}.

\sm

{\it Lemma {\rm \ref{l:atomic}} implies that
for $(a,b)\notin \Pi$ the operators $D_+$, $D_-$ have no commuting
self-adjoint extensions}.

\sm

{\bf\punct Almost proof of Lemma \ref{l:atomic}.} For $\zeta$ being in  a general position,
the system (\ref{eq:DDPhi}) has a unique solution, and  it has the form $\FF[\cdot;z]$.
Denote by $\Omega_{hyp}$ the subset of $\Omega$ consisting of the functions 
$\FF[\cdot;z]$.
We wish to prove the following statement:

\begin{lemma} 
\label{l:with-exceptions}
Let $(a,b)\notin \Pi$ and $a+b$, $a-b$, $a$, $b\notin\Z$. 
Let $\rho$ be a measure on $\Omega$, and let the corresponding operator $U$ be bounded.
Then $\rho$ is atomic on $\Omega_{hyp}$.
\end{lemma}

{\sc Proof.}
Set $\zeta=\lambda^2$. Then a hypergeometric solution of the system
(\ref{eq:DDPhi}) has one of the two forms:
$$
\FF\left[\begin{matrix}
    a+\lambda|a-\ov \lambda, a-\lambda|a+\ov \lambda
    \\ a+b|a+b
   \end{matrix},z
\right],\qquad
\FF\left[\begin{matrix}
    a+\lambda|a+\ov \lambda, a-\lambda|a-\ov \lambda
    \\ a+b|a+b
   \end{matrix},z
\right].
$$
In the first case we have $(a+\lambda)-(a-\ov\lambda)=2\Re \lambda\in \Z$,
hence $\lambda=\frac12(k+is)$, where $k\in\Z$, $s\in\R$.
We come to the functions $\cK_{a,b}(z;k,is)$. 

In the second case we have $\lambda-\ov\lambda\in\Z$, i.e., $\lambda=\tau\in\R$.
We come  to the functions
\begin{equation}
\cK(z;0,\tau)=
\FF\left[\begin{matrix}
    a+\tau|a+\tau, a-\tau|a-\tau
    \\ a+b|a+b
   \end{matrix},z
\right].
\label{eq:complementary}
\end{equation}

Next, we will show that 
\begin{equation}
\text{\it
the measure 
$\rho$ is zero on the set of all $\lambda=\tfrac12(k+is)$
with $s\ne 0$.}
\label{eq:text}
\end{equation}

Our kernel has the following asymptotics at $z=0$ and $z=1$:
\begin{align}
\cK(z;k,is)=& \bigl(1+O(z)\bigr)+B(k,is) |z|^{2-2a-2b} \bigl(1+O(z)\bigr)\qquad \text{as $z\to 0$},
\label{eq:asz0}
\\
\cK(z;k,is)=& C(k,is) \bigl(1+O(1-z)\bigr)+
\nonumber
\\
\quad\quad &+
 D(k,is) |1-z|^{2b-2a}\bigl(1+O(1-z)\bigr)\qquad \text{as $z\to 1$},
\label{eq:asz1}
\end{align}
where the coefficients $B$, $C$, $D$ are continuous non-vanishing functions on $\Lambda$
and all $O(\cdot)$ are uniform on compact subsets of $\Lambda$
(see formulas (\ref{eq:A})--(\ref{eq:B11})).

For definiteness, assume that $a+b>2$. Consider a point $(k_0,is_0)\in \Lambda$, $s_0\ne 0$
and a neighborhood $\cN$ of $(k_0,is_0)$. Assume that $\rho(\cN)>0$. Denote by $I_\cN$ the indicator function of the set
$\cN$. The function $U I_\cN$ has the following asymptotics at
$z=0$:
$$
U I_\cN(z)= \alpha\,(1+O(z))+\beta\, |z|^{2-2a-2b} \bigl(1+O(z)\bigr) \qquad \text{as $z\to 0$}.
$$
Due to uniformity $O(\cdot)$, for a sufficiently small neighborhood $\cN$ we have
$\alpha\ne 0$, $\beta\ne 0$. Since $a+b>2$, 
the actual asymptotics is
$$
U I_\cN(z)= \beta\, |z|^{2-2a-2b} (1+O(z)).
$$
Therefore
$$U I_\cN\notin L^2\bigl(\C, |z|^{2a+2b-2}|1-z|^{2a-2b}\dd z\bigr).$$
This contradicts to boundedness of $U$. Thus any point has a neighborhood of zero measure,
and this implies claim (\ref{eq:text}) in the case $a+b>1$.

In domains $a+b<0$, $a-b<-1$, $a-b>1$ we get the same effect.

\sm

Next, examine the complementary series
$\cK(z;0,\tau)$
of eigenfunctions, see (\ref{eq:complementary}).
We have the same asymptotics (\ref{eq:asz0})--(\ref{eq:asz1}),
we only must write the coefficients of the form $A(0,\tau)$, $B(0,\tau)$, $C(0,\tau)$, $D(0,\tau)$
in (\ref{eq:asz0})--(\ref{eq:asz1}).
These functions have zeros and poles on the axis $\tau\in\R$.
The same argument as above shows that if $\tau_0$ is not a zero and not a pole
of all our coefficients, then the measure $\rho$ is zero on a sufficiently small neighborhood
of $\tau_0$. The set of zeros and poles is countable. This completes the proof of the lemma.
\hfill $\square$

\sm

{\bf\punct Proof of non-self-adjointness.%
\label{ss:non-self}} However, our system of differential equations
(\ref{eq:DDPhi})
has solutions that have not the form $\FF$, and enumeration of all possible degenerations is tedious.
So we continue the proof of Lemma \ref{l:atomic} without constrains of Lemma \ref{l:with-exceptions}.
Due to the homographic transformations, without loss of generality we can set
  \begin{equation}
  a+b>2.
  \label{eq:a+b2}
  \end{equation}
First, we examine asymptotics in a neighborhood of $z=0$.

\sm

{\sc Asymptotics at $z=0$. Non-logarithmic case.}
If $a+b\ne 2$, $3$, \dots, then the equation $\frD \Phi=\lambda^2\Phi$ has two holomorphic solutions,
$$
\Psi_1(z):=\F\left[\begin{matrix}a+\lambda,a-\lambda\\a+b\end{matrix};z\right],\quad 
\Psi_2(z):=z^{1-a-b}\F\left[\begin{matrix}1-b+\lambda,1-b-\lambda\\2-a-b\end{matrix};z\right].
$$
The equation $\ov\frD \Phi=\ov\lambda^2\Phi$ has two antiholomorphic solutions 
$$
\wt \Psi_1(\ov z):=\F\left[\begin{matrix}a+\ov\lambda,a-\ov\lambda\\a+b\end{matrix};\ov z\right],\qquad 
\wt \Psi_2(\ov z):=\ov z^{\,1-a-b}\F\left[\begin{matrix}1-b+\ov\lambda,1-b-\ov\lambda\\2-a-b\end{matrix};\ov z\right].
$$
Therefore a single-valued solution of the system must have the form
$$
A\Psi_1(z) \wt \Psi_1(\ov z)+ B \Psi_2(z) \wt\Psi_2(\ov z).
$$
The first term has $L^2(\C,\mu_{a,b})$-asymptotics at $z=0$, by (\ref{eq:a+b2}) the second term has non-$L^2$-asymptotics.
Thus  the spectral measure $\rho$ is supported by the set of functions of the form  $\Psi_1(z) \wt\Psi_1(\ov z)$.

\smallskip

{\sc  Asymptotics at $z=0$. Logarithmic case.} Now let $a+b=n=2$, 3,\dots.
Then the equation  $\frD \Phi=\lambda^2\Phi$ has two holomorphic solutions,
$$
\Psi_1(z)=\F[a+\lambda, a-\lambda;n,z],\qquad \Psi_2(z),
$$
where $\Psi_2(z)$ is a logarithmic solution, which has the form 
(\ref{eq:logarithmic-solution}).
The equation  $\frD \Phi=\ov\lambda^2\Phi$ has two antiholomorphic solutions,
$$
\wt\Psi_1(\ov z)=\F[a+\ov\lambda, a-\ov\lambda;n,\ov z],\qquad \wt\Psi_2(\ov z).
$$
A single valued solution must have the form
$$
A\Psi_1(z)\,\wt\Psi_1(\ov z)+B \Bigl(
\Psi_1(z)\wt\Psi_2(\ov z)+ \Psi_2( z)\wt\Psi_1(\ov z)
\Bigr).
$$
The asymptotics of the second summand is $(z^{-n+1}+\ov z^{\,-n+1})+O(z^{-n+2})$ if 
$n\ge 3$. If $n=2$ we have $(z^{-1}+\ov z^{\,-1})+O(z^{-\epsilon})$. We get a
non-$L^2$ asymptotics.

\sm

Thus, {\it for $a+b>2$ the spectral measure is supported by set of functions
of the form
$\Psi_1(z)\wt\Psi_1(\ov z)$}.

\sm

{\sc Single-valuedness near $z=1$. Non-logarithmic case.}
Assume that $a-b\notin\Z$.
We apply formulas Erd\'elyi, et al., \cite{HTF1}, (2.10.1), (2.10.5) and write explicit  expansions of
$\Psi_1$, $\ov\Psi_1$ at $z=1$.
\begin{align*}
\Psi_1(z)=A_1 G_1(1-z)+A_2 (1-z)^{b-a} G_2(1-z);
\\
\wt\Psi_1(\ov z)=\wt A_1 \wt G_1(1-\ov z)+\wt A_2 (1-\ov z)^{b-a} \wt G_2(1-\ov z),
\end{align*}
where $G_1$, $G_2$ are certain series $\F$ and the coefficients $A_1$, $A_2$ are products of gamma functions,
see the explicit formulas (\ref{eq:analytic-cont})--(\ref{eq:analytic-cont1})  above.
Clearly, the product $\Psi_1(z)\wt\Psi_1(\ov z)$ can be single-valued only if 
$A_2=\wt A_2=0$, or
$A_1=\wt A_1=0$. 
Looking to the explicit expressions for the gamma-coefficients, we observe that the first case happens
if both  hypergeometric series $G_1(z)$, $G_2(z)$ are terminating
(i.e., $a-\lambda=0$, $-1$, \dots or  $a+\lambda=0$, $-1$, \dots, in particular, $\lambda$
is real). The second variant holds if and only if both  series $G_2(1-z)$, $\wt G_2(1-\ov z)$
are terminating (i.e., $b-\lambda=0$, $-1$, \dots or  $b+\lambda=0$, $-1$, \dots).

\sm

{\sc Single-valuedness near $z=1$. Logarithmic case.}
Now let $b-a\in\Z$.
The transposition $a\leftrightarrow b$ corresponds to a homographic transformation
of differential operators,
it preserves the condition $a+b\ge 2$. Therefore we can assume $m:=b-a\ge 0$. 
Represent
$\Psi_1(z)$, $\wt\Psi_1(\ov z)$ 
as combinations of  basic solutions of the hypergeometric equations
at the point $z=1$,
\begin{align*}
\Psi_1(z)=
A\,\, \F[a+\lambda,a-\lambda; b-a+1;z]+B\, \Theta(1-z);
\\
\wt \Psi_1(\ov z)=
\wt A\,\, \F[a+\ov\lambda,a-\ov\lambda; b-a+1;\ov z]+\wt B\, \wt\Theta(1-\ov z),
\end{align*}
where $\Theta(1-z)$ is a logarithmic series of the type (\ref{eq:logarithmic-solution}), see
 Erd\'elyi, et al., \cite{HTF1}, (2.10.12).
A straightforward calculation shows that the  product $\Psi_1(z)\,\wt\Psi_1(\ov z)$  can be single valued near $z=1$
only if $B=\wt B=0$. Therefore $\Psi_1(z)$ is  single valued near $z=1$, and therefore
it is a single valued solution of
a hypergeometric equation on the whole plane $\dot\C$.
Hence (see Erd\'elyi, et al., \cite{HTF1}, Subset. 2.2.1) $\Psi_1(z)$ is a polynomial.

\sm

{\sc Behavior at infinity.} Thus the spectral measure $\rho$
is supported by generalized  eigenfunctions of the  following types
$$
p_1(z)\, p_2(\ov z),\qquad (1-z)^{b-a|b-a}\, q_1(z)\, q_2(\ov z),
$$
where $p_j$, $q_j$ are polynomials. However, our density
$\mu_{a,b}(z)$ has a behavior $\sim |z|^{4a-2}$ at infinity and therefore the space $L^2$ 
can contain only a finite number  orthogonal functions of such a type.
\hfill $\square$

\section{Symmetry of differential operators\label{s:symmetry1}}

\COUNTERS

Here we show that $J^*_{a,b}$ sends $\cD_\even(\dot{\Lambda})$ to $\cR_{a,b}$
and verify that  $\frD$ and $\ov\frD$ are  adjoint one to another  on $\cR_{a,b}$.

\sm

In this section we denote by $D_r(u)\subset \C$ (resp. $\ov D_r(u)$) the open (resp. closed) disc in $\C$ of radius $r$ with center at $u$. By $S_r(u)$ we denote 
the circle $|z-u|=r$.

\sm

{\bf \punct The map $J^*_{a,b}$ on the space $\cD_\even(\dot\Lambda)$.}
Introduce  a natural topology in the space $\cR_{a,b}(\dot\C)$  defined 
in Subsect. \ref{ss:1.1}.  
Consider the space $\cR(0)$ of functions in $\ov D_{1/3} (0)$ having the form
$\alpha(z)+\beta(z)|z|^{2a+2b-2}$, where $\alpha(z)$, $\beta(z)$ are smooth in
$\ov D_{1/3} (0)$ up to the boundary.
Let $C^\infty_{\mathrm{flat}}(\ov D_{1/3}(0))\subset C^\infty (\ov D_{1/3}(0))$ be the subspace consisting of all
functions
that are flat
at 0. 
The space $\cR(0)$ is a quotient space
$$
\cR(0)\simeq
\Bigl[
C^\infty(\ov D_{1/3}(0))\oplus |z|^{2a+2b-2}C^\infty(\ov D_{1/3}(0))\Bigr]/ C^\infty_{\mathrm{flat}}(\ov D_{1/3}(0)).
$$
We equip $\cR(0)$  with the topology of a quotient space. In the same way
we define a topology in the space $\cR(1)$ of smooth functions in $\ov D_{1/3}(1)$ having the form
$\gamma(z)+\delta(z)|1-z|^{2a-2b}$.

We define a topology in $\cR_{a,b}$ as a weakest topology
such that:

\sm

a) The restriction operators
$$
\cR_{a,b}\to \cR(0), \quad \cR_{a,b}\to \cR(1),\quad
\cR_{a,b}\to C^{\infty}\Bigl(\ov D_2(0)\setminus (D_{1/3}(0)\cap D_{1/3}(1))\Bigr)
$$
 are continuous.

\sm 

b) For all $\alpha$, $\beta$, $N$ the following seminorms are continuous
\begin{equation}
p_{\alpha,\beta,N}(f)=
\sup_{\C\setminus D_2(0)} |z|^{2+\alpha+\beta} (\ln|z|)^N
\biggl|\frac{\partial^{\alpha+\beta} f(z)}{\partial z^\alpha\partial\ov z^\beta} \biggr|.
\label{eq:seminorm}
\end{equation}

Recall that
$
\dot\Lambda:=\Lambda\setminus\{(0,0)\}
$.

\begin{lemma}
For $|z|>2$, $(k,s)\in \dot\Lambda$ 
we have the following expansion
\begin{multline}
 \cK(z;k,is)=
 z^{-a-\frac{k+is}2\bigl|-a-\frac{-k+is}2} B(k,s; z^{-1})+\\+
  z^{-a+\frac{k+is}2\bigl|-a+\frac{-k+is}2} B(-k,-s;  z^{-1}),
  \label{eq:nadoelo}
\end{multline}
where the expression $B(k,s;u)$ for fixed $k$ is smooth $s$ except the point
$(k,s)=(0,0)$.
\end{lemma}

{\sc Proof.}
We refer to  expansion (\ref{eq:C})--(\ref{eq:C11}). Notice that
for $k=0$, $s=0$ we have a singularity in this expansion (but the kernel 
itself is analytic at this point).
\hfill $\square$

\begin{proposition}
\label{l:JRab}
{\rm a)}
 Let $\Phi\in \cD_\even(\dot\Lambda)$. Then
 $J^*_{a,b} \Phi\in \cR_{a,b}$.
 
 \sm
 
{\rm b)} Moreover, the operator $J_{a,b}^*$ is a continuous operator from  $\cD_\even(\dot \Lambda)$
to $\cR_{a,b}$.
\end{proposition}

{\sc Proof.} Forms of asymptotics of  $J^*_{a,b}\Phi$ at 0 and 1 follow from the expressions
(\ref{eq:A}), (\ref{eq:B}). Let us examine the asymptotics at $z\to\infty$.
Without loss of generality we can assume that $|k|$ is  fixed.
We write
\begin{equation*}
 J^*_{a,b} \Phi(z)
 =
 z^{-a-\frac{k}2|-a+\frac{k}2}  \int_\R  z^{-\frac{is}2|-\frac{is}2}  B(k,s;  z^{-1})\, \Phi(k,s)\,ds
 +
 \Bigl\{ \begin{matrix}
          \text{similar}\\\text{term}
         \end{matrix}.
\Bigr\}
\end{equation*}
Differentiating the first summand by $\frac {\partial^{\alpha+\beta}}{\partial z^\alpha\partial \ov z^\beta}$ and
keeping in  mind (\ref{eq:C}) and Lemma \ref{l:differentsiruem}, we get an expression
of the form
\begin{multline*}
 z^{-a-\frac{k}2-\alpha|-a+\frac{k}2-\beta} 
\sum_{0\le p\le \alpha, 0\le q \le \beta}
 \int_\R  z^{-\frac{is}2|-\frac{is}2}
 U^{\alpha,\beta}_{p,q}(a,b,k,s)
 \times\\\times 
 \frac{\Gamma^\C(-k-is|k+is)\,(a+\frac{k+is}2)_p\, (a+\frac{-k+is}2)_p\,
(a+\frac{k-is}2)_q\, (a+\frac{-k-is}2)_q
}
{\Gamma^\C(b-\frac{k+is}2|b-\frac{-k+is}2)\Gamma^\C( a+\frac{k+is}2|a+\frac{-k+is}2)
\,(a+b)_p\,(a+b)_q}
 \times\\\times 
 \F\left[  
 \begin{matrix}
 a+\frac{k+is}2+p, a+\frac{-k+is}2+p
 \\
 a+b+p
 \end{matrix}; z^{-1}
 \right] 
   \F\left[  
 \begin{matrix}
 a+\frac{-k-is}2+q, a+\frac{k-is}2+q
 \\
 a+b+q
 \end{matrix}; \ov z^{\,-1}
 \right] 
  \times\\\times 
  \Phi(k,s)\,ds,
\end{multline*}
where $ U^{\alpha,\beta}_{p,q}(a,b,k,s)$ are polynomials.
It is easy to verify that the integrand is a smooth compactly supported function on $\dot\Lambda$.
Next, we write
$$
|z|^{-is} =\frac{i}{\ln |z|} \frac\partial{\partial s} |z|^{-is},
$$
integrate our   expansion by parts $N$ times and observe that
$p_{\alpha,\beta,N}(J^*_{a,b}\Phi)<\infty$.

The continuity follows from the same considerations.
\hfill $\square$

\sm

As a corollary, we obtain the following lemma.

\begin{lemma}
\label{l:J*C}
 The operator $J^*_{a,b}$ is continuous as an operator from $\cD_\even(\dot\Lambda)$ to the space $L^2(\C, \mu_{a,b})$.
\end{lemma}

{\sc Proof.} Indeed, for $(a,b)\in\Pi$ the identical
embedding $f\mapsto f$ of $\cR_{a,b}$ to $L^2(\C, \mu_{a,b})$
is continuous. \hfill $\square$

\begin{lemma}
 If $f\in \cR_{a,b}$, then $\frD f\in \cR_{a,b}$.
\end{lemma}

{\sc Proof.} Let us check the behavior of $\frD f$ at 0, For definiteness 
assume that $a+b\ne1$. Then near zero we have
\begin{multline*}
\frD f=\frD\bigl(\alpha(z)+\beta(z) |z|^{1-a-b}\bigr)=
\\=
\Bigl\{\Bigl(z(1-z) \frac{d^2}{dz^2}+ (a+b)\frac d{dz}\Bigr)z^{1-a-b}
\Bigr\}\cdot
\ov z^{1-a-b} \beta(z) + \Bigl\{\text{the rest}\Bigr\}.
\end{multline*}
Obviously, the rest has the form $\wt\alpha(z)+\wt\beta(z)|z|^{2-a-b}$
with smooth $\wt\alpha$, $\wt\beta$. The expression in the curly brackets%
\footnote{Cf. \cite{Fed2}, Sect.I.2.}
is $-(a+b)(a+b-1)z^{1-a-b}$. \hfill $\square$

\sm

{\bf \punct Symmetry of differential operators.}

\begin{proposition}
	\label{pr:D-symmetry}
	For any $f$, $g\in \cR_{a,b}(\dot\C)$
$$
\la \frD f, g\ra=\la f,\ov \frD g\ra.
$$
\end{proposition}

{\sc Proof.} Let $f$, $g\in \cR_{a,b}$.
We wish to show that
$$
\int_{\dot\C}\bigl( \frD f(z)\cdot \ov {g(z)}- f(z)\cdot \frD \ov {g(z)}\bigr)\,\mu_{a,b} \dd z=0.
$$
By Lemma \ref{l:J*C}, $\frD f$, $\frD g\in L^2(\C,\mu_{a,b})$.
Therefore our improper integral absolutely converges,
 we write  it as 
$$
\lim_{\epsilon\to 0} \int_{D_{1/\epsilon} (0)\setminus \bigl(D_\epsilon(0)\cup D_\epsilon(1)\bigr)} \bigl(\dots\bigr)\,
\frac{dz\wedge d \ov z}{2i}.
$$
Next, we integrate two times by parts in $z$ (with the Green formula) and
after a simple calculation come to
\begin{equation}
\lim_{\epsilon\to 0}
\Biggl\{
\int_{S_{1/\epsilon}(0)} V(z) d\,\ov z- \int_{S_{\epsilon}(0)} V(z)\, d\ov  z- \int_{S_{\epsilon}(1)} V(z)\, d\ov z
\biggr\}
\label{eq:limSSS}
,
\end{equation}
where 
$$
V(z)=\Bigl(\frac{\partial f}{\partial z}\cdot \ov{g(z)}-
f(z)\cdot \frac{\partial \ov{g(z)}}{\partial z} 
\Bigr) z(1-z) \,\mu_{a,b}(z).
$$
We claim that all summands in (\ref{eq:limSSS}) tend to 0. For the first summand this is clear. For the second summand
we represent $f$, $\ov g$ as 
$$
f(z)=\alpha(z)+\beta(z)\, z^{1-a-b|1-a-b},\qquad
\ov g(z)=\gamma(z)+\delta(z)\, z^{1-a-b|1-a-b}.
$$
Then $V(z)$ transforms to an expression of the following type:
\begin{multline*}
\Bigl(A(z)+ B(z) z^{-a-b|1-a-b} + C(z) z^{2-2a-2b|2-2a-2b} \Bigr)
\times\\\times
z(1-z)\cdot z^{a+b-1|a+b-1}(1-z)^{a-b|a-b},
\end{multline*}
where $A(z)$, $B(z)$, $C(z)$ are smooth near 0.
We emphasize that the term with $ z^{1-2a-2b|2-2a-2b}$
in the bracket
appears with the coefficient
$$
(2-2a-2b)\bigl(\beta(z)\delta(z)-\beta(z)\delta(z)\bigr)=0.
$$
 Thus we get summands with the following behavior at 0:
$$
\sim A(0)\, z^{a+b|a+b-1},\qquad \sim B(0)\,  z^0, \qquad \sim C(0)\,  z^{2-a-b|1-a-b}.
$$
Since $0<a+b<2$ all powers are $>-1$ and therefore%
\footnote{See a discussion of a parallel situation for ordinary differential operators in
\cite{Ner-jacobi}, Section 1. However, in the one-dimensional case we must impose boundary conditions in
such points.}
$\int_{|z|=\epsilon}(\dots)\,d\ov z$ tends to 0.\,
 $\square$

\section{The operator $J_{a,b}^*$ is an isometry\label{s:isometry1}}

\COUNTERS

Here we prove half of Theorem \ref{th:spectral1}.

\sm

{\bf \punct The statement.} 
First, denote by $\Lambda_+$ the subset of $\Lambda$ consisting of 
$(k+is)/2$ such that $k>0$ or $k=0$ and $s>0$. We have an
obvious identification $\cD_{\even}(\dot\Lambda)\simeq \cD(\Lambda_+)$. 

\begin{lemma}
 \label{l:J*-embed}
 Let $u(\lambda)$, $v(\lambda)$ be smooth compactly supported 
 function on $\Lambda_+$. Then 
  $$\la J_{a,b}^* u, J_{a,b}^*v \ra_{L^2( \dot \C,\mu_{a,b})}=
 2 \la u,v\ra_{L^2(\Lambda_+,\kappa_{a,b})}.
 $$ 
\end{lemma}

Our proof is based on heuristic arguments outlined in Berezin, Shubin \cite{BSh}, Section 2.6,
for ordinary differential operators.
However, this way is tiresome.

\sm

{\bf\punct Preliminary remarks.}
Recall that
$$
J^*_{a,b} u(z)=2\int_{\Lambda_+} u(\lambda)\, \cK(z,\lambda)\,\kappa_{a,b}(\lambda)
\wt d\lambda.
$$
By Lemma \ref{l:J*C}, this operator is continuous
 as an operator $\cD(\dot\Lambda_+)\to L^2( \dot \C,\mu_{a,b})$.
Therefore the sesquilinear form
\begin{equation}
T(u,v):=
\la J_{a,b}^* u, J_{a,b}^*v \ra_{L^2( \dot \C,\mu_{a,b})}
\label{eq:lara}
\end{equation}
is continuous as a form $\cD(\Lambda_+)\times \cD(\Lambda_+) \to \C$.
By the kernel theorem (see, e.g., \cite{Hor}, Sect. 5.2)  it is determined by a distribution.
Formally, we transform (\ref{eq:lara}) as
\begin{align}
& 
 \int\limits_{\dot\C} \Bigl( \int\limits_{\Lambda_+} u(\lambda) \cK(z,\lambda)\kappa_{a,b}(\lambda)\,\wt d\lambda\Bigr)
 \cdot \Bigl(\int\limits_{\Lambda_+} \ov{v(\nu)} 
 \cK(z,\nu)\kappa_{a,b}(\nu)\,\wt d\nu\Bigr)\,\, \mu_{a,b}(z)\,\dd z
 \label{eq:line-1}
 =\\
 &=
 \int_{\Lambda_+} \int_{\Lambda_+} u(\lambda)\ov {v(\nu)}
 H(\lambda,\nu)\kappa_{a,b}(\lambda)\kappa_{a,b}(\nu)\,\wt d\lambda\,\wt d\nu,
 \label{eq:line-2}
\end{align}
where
\begin{equation}
  H(\lambda,\nu)=
  \int_{\dot\C}
 \cK(z,\lambda)\, \cK(z,\nu)\,\mu_{a,b}(z)\,\dd z. 
 \label{eq:line-3}
\end{equation}

Notice that all integrals in line (\ref{eq:line-1}) converge absolutely.
However, the triple integral $\int_{\Lambda_+}\int_{\Lambda_+}\int_{\dot \C}$
is not absolutely convergent. The integrand in (\ref{eq:line-3})
decreases as $|z|^{-2}$ and the integral diverges.

 However, we  regard $H(\lambda,\nu)$ as a distribution, then
 Lemma \ref{l:J*-embed} can be reformulated in the form:
 
 \begin{lemma} 
 \label{l:J-delta}
 We have the following identity of distributions on $\cD(\Lambda_+)\times \cD(\Lambda_+)$:
\begin{equation}
 H(\lambda,\nu)=\delta(\lambda-\nu).
 \label{eq:J-delta}
\end{equation}
\end{lemma}

{\bf \punct Orthogonality of packets.}

\begin{lemma}
\label{l:orthogonal}
 Let $u$, $v\in \cD(\Lambda_+)$
 and  supports  $\supp(u)$, $\supp(v)$ have empty intersection.
 Then
 $$
 \la J_{a,b}^* u, J_{a,b}^*v \ra_{L^2( \dot \C,\mu_{a,b})}=0
 .
 $$
\end{lemma}

\sm

{\sc Proof.} 
Denote $D_+:=\frac12(\frD+\ov\frD)$, $D_-=\frac1{2i}( (\frD-\ov \frD)$.
By Proposition
\ref{l:JRab}, $J^*_{a,b}u$ is contained in the space $\cR_{a,b}$.
By
Proposition \ref{pr:D-symmetry}, the operators $D_+$, $D_-$ are formally symmetric on $\cR_{a,b}$. Since they formally commute, for any real polynomial
$p(D_+,D_-)$ we have
$$
\la p(D_+,D_-) J_{a,b}^* u, J_{a,b}^*v \ra= \la  J_{a,b}^* u, p(D_+,D_-) J_{a,b}^*v \ra,
$$
or 
\begin{equation}
\la  J_{a,b}^*\, p(\Re \lambda, \Im\lambda) \cdot u, J_{a,b}^*v \ra=
\la  J_{a,b}^* u, J_{a,b}^*\, p(\Re \lambda, \Im\lambda)\cdot v \ra
,
\label{eq:orthogon}
\end{equation}
where $\cdot$ denotes the operator of multiplication by a function.
We choose a sequence $p_N$ of polynomials such that $p_N$ uniformly
converges to 1 on $\supp(u)$ with all derivatives and converges to
$0$ on $\supp(v)$. By Lemma \ref{l:J*C}  the map $J^*_{a,b}$ is continuous as a map $\cD(\Lambda_+)\to
L^2(\C,\mu_{a,b})$.
Replacing $p$ by $p_N$
in (\ref{eq:orthogon}) and passing to a limit,
we come to the desired statement.
\hfill $\square$

\sm

{\bf\punct Next reduction of our statement.} 
%
Let $S(u,v)$ be an Hermitian form on $\cD(\Lambda_+)$.
We say that $S$ is {\it $C^\omega$-smooth} if it has the form
$$
S(u,v)=\int_{\Lambda_+} \int_{\Lambda_+}
M(\lambda,\nu)\, u(\lambda)\ov{v(\nu)}
 \,\wt d\lambda\,\wt d\nu,
$$
where $M$ is a real analytic function on $\Lambda_+\times \Lambda_+$.

\begin{lemma}
\label{l:orthogonal2}
We have
\begin{equation}
\la  J_{a,b}^* u, J_{a,b}^*v \ra_{L^2(\C, \mu_{a,b})}
=
\la u, v\ra_{L^2(\Lambda,\kappa_{a,b})}+ S(u,v),
\label{eq:o1}
\end{equation}
where $S(u,v)$ is $C^\omega$-smooth.
\end{lemma}

This lemma together with Lemma \ref{l:orthogonal} imply the desired statement, i.e., the identity
(\ref{eq:J-delta}). Indeed, for any $u$, $v$ with disjoint support,
we have
$$
\int_{\Lambda_+} \int_{\Lambda_+}  M(\lambda,\nu) \,u(\lambda)\, \ov v(\nu)\, \kappa_{a,b}(\lambda)\kappa_{a,b}(\nu)\,\wt d\lambda\,\wt d\nu=0,
$$
 therefore $M(\lambda,\nu)=0$.

\sm

The rest of this section is occupied by the proof of Lemma \ref{l:orthogonal2}.

\sm

{\bf \punct Beginning of the proof of Lemma \ref{l:orthogonal2}. Cleaning of the problem.}
{\it Step} 1.
 Represent
$$
u=\sum_k u_k\,\delta(\Re \lambda-k/2), \qquad
v=\sum_l v_l\,\delta(\Re \lambda-l/2),
$$
in fact the sums are finite and $u_k$, $v_l$ depend on a real variable $s$.
By Lemma \ref{l:orthogonal}, we have
$$\la  J_{a,b}^* u_k, J_{a,b}^*v_l \ra=0\qquad \text{for $k\ne l$}.
$$
Therefore it is sufficient to examine only inner products 
$$
\la  J_{a,b}^* u_k, J_{a,b}^*v_k \ra=
\int_{\dot\C} R(z)\,\dd z,
$$
where
\begin{multline*}
R(z):= \int_{\Lambda_+} u_k(is) \,
 \cK\bigl(z,\tfrac12(k+is)\bigr)\kappa_{a,b}\bigl(\tfrac12 (k+is)\bigr)
 \,ds
 \times\\\times
 \int_{\Lambda_+} \ov{v_k(it)}\, \cK\bigl(z,\tfrac12(k+it)\bigr)
 \kappa_{a,b}\bigl(\tfrac12(k+it)\bigr)\,dt\,\, \mu_{a,b}(z).
 \end{multline*}
  
{\it Step} 2. Represent the integral as
$\int_{|z|\le 2} R+ \int_{|z|\ge 2} R$. 

Let us show that the first summand is $C^\omega$-smooth.
In this case the triple integral absolutely converges and can be written as
$$
\int_{|z|\le 2} R\,\dd z=
\int_{\R} \int_{\R} u_k(is) \ov{v_k(it)}
L(s,t) \,ds\,dt,
$$
where
$$
L(s,t)=
\int_{|z|\le 2} \cK\bigl(z,\tfrac12(k+is)\bigr)\,\cK\bigl(z,\tfrac12(k+it)\bigr)\,\mu_{a,b}(z)\,\dd z.
$$
Integrand makes sense for complex $s$, $t$ that are sufficiently
close to $\R$ 
and the integral absolutely converges
(singularities at $z=0$ and $1$ have the forms (\ref{eq:asz0}), (\ref{eq:asz1})).
Therefore $L(s,t)$ is a holomorphic function in $s$, $t$ near $\R\times\R$.

\sm

Therefore our question is reduced to an examination the integral
$$
\int_{|z|>2} R(z)\,\dd z
$$

{\it Step} 3.
{\it A decomposition of the kernel.}
Applying Theorem \ref{th:long}.c,
we write $\cK(z,\lambda)$ in the domain
$|z|\ge 2$
as
\begin{multline}
\cK(z,\lambda)=
W_1+W_2+W_3:=\\=
A(\lambda) (-z)^{-a-\lambda|-a+\ov\lambda}+
A(-\lambda) (-z)^{-a+\lambda|-a-\ov\lambda}+\Psi(z,\lambda)
,
\label{eq:cK-exp1}
\end{multline}
where 
$$
A(\lambda)=\frac{ \Gamma^\C(2\lambda|-2\ov\lambda)}
{\Gamma^\C(b-\lambda|b+\ov\lambda) \Gamma^\C(a-\lambda|a+\ov\lambda)}
$$
and
$$
\Psi(z,\lambda)=O(|z|^{-2a-1})\qquad \text{as $z\to\infty$.}
$$
Notice that
\begin{equation}
|A(\lambda)|^2=A(\lambda)\, A(-\lambda)=\kappa_{a,b}^{-1}(\lambda).
\label{eq:|A|}
\end{equation}

 Therefore the integral $\int_{|z|>2} R(z)\,\dd z$
splits into a sum of 9 summands $V_{\alpha\beta}$,
where $\alpha$, $\beta=1$, 2, or 3,
\begin{multline*}
V_{\alpha\beta}:=\int_{|z|>2} 
\int_{\R} W_\alpha(z;k,s) u_k(is)\,\kappa_{a,b}\bigl(\tfrac12(k+is)\bigr)\,ds
\times\\\times
\int_{\R} W_\beta(z;k,t) v_k(it)\,\kappa_{a,b}\bigl(\tfrac12(k+it)\bigr)\,dt
\cdot \mu_{a,b}(z)\,\dd z.
\end{multline*}

\sm

{\it Step} 4.
 For five summands $V_{13}$,
$V_{23}$, $V_{31}$, $V_{32}$, $V_{33}$ 
we  immediately get
absolute convergence of triple integrals and 
$C^\omega$-smoothness.
For instance, 
\begin{multline*}
 V_{13}=\int_{\R} \int_{\R} u_k(is) \ov{u_k(it)}\,
  A\bigl(\tfrac12(k-is)\bigr)^{-1} 
  \times\\\times
  \biggl[\int_{|z|\ge 2} \Bigl(\frac z{\ov z}\Bigr)^{-k}|z|^{-2a+is}\,
\ov{\Psi\bigl(z,\tfrac12 (k+is)\bigr)}\,\mu_{a,b}(z)\dd z\biggr]
 \,ds\,dt.
\end{multline*}
(we simplified the integrand using (\ref{eq:|A|})).
The expression in  the square brackets is real analytic (the integrand decreases as $|z|^{-3}$).

\sm

{\it Step} 5.
Non-obvious summands are $V_{11}$, $V_{12}$, $V_{21}$, $V_{22}$.
We start with $V_{11}$,
\begin{multline*}
V_{11}:=\int_{|z|\ge 2}
  \int_\R u_k(is) 
 A\bigl(\tfrac12(k+is)\bigr) \Bigl(\frac z{\ov z}\Bigr)^{-k/2}|z|^{-2a-is}
 \kappa_{a,b}\bigl(\tfrac12 (k+is)\bigr)
 \,ds
 \times\\\times
 \int_{\R} \ov{ A\bigl(\tfrac12(k+it))\,v_k(it)} 
 \Bigl(\frac z{\ov z}\Bigr)^{k/2}|z|^{-2a+it}
 \kappa_{a,b}\bigl(\tfrac12(k+it)\bigr)\,dt\,\, \mu_{a,b}(z)\,\dd z.
\end{multline*}
For $k=0$ we must keep in mind that the integration  $\int_\R$ 
actually is taken over a ray $[\epsilon,\infty)$ for some $\epsilon>0$.
Applying (\ref{eq:|A|}), we come to
\begin{multline}
V_{11}:= \int\limits_{|z|\ge 2}
  \int\limits_\R \int\limits_\R u_k(is)\ov{v_k(it)}  A\bigl(\tfrac12(k-is)\bigr)^{-1}
  A\bigl(\tfrac12(k+it)\bigr)^{-1} |z|^{-4a-is+it} \,ds \,dt
  \times\\\times
  \mu_{a,b}(z)\,\dd z.
  \label{eq:V11}
\end{multline}
Next, we notice that
$$
\mu_{a,b}(z)=|z|^{2a+2b-2}|1-z|^{2a-2b}=|z|^{4a-2}+O(|z|^{4a-3})\qquad 
\text{as $z\to\infty$.}
$$
We write
\begin{equation}
\mu_{a,b}(z)=|z|^{4a-2}+(\mu_{a,b}(z)-|z|^{4a-2}),
\label{eq:mu-sum}
\end{equation}
substitute this to (\ref{eq:V11}) and decompose (\ref{eq:V11}) as a sum of two integrals.
The second summand immediately gives a $C^\omega$-smooth term.
The first summand  is the topic of our interest. It equals 
the following expression: 
\begin{multline}
I(u,v):= 
\int\limits_{|z|\ge 2}
  \int\limits_\R \int\limits_\R u_k(is)\ov{v_k(it)}  A\bigl(\tfrac12(k-is)\bigr)^{-1}
  A\bigl(\tfrac12(k+it)\bigr)^{-1}
  \times\\\times
  \boxed{ |z|^{-2-is+it}} \,ds \,dt
  \,\dd z.
  \label{eq:Iuv}
\end{multline}

{\bf\punct Application of the Sokhotski formula and disappearance of a singular term.}

\sm

{\it Step} 6. {\it Extension to the complex domain.}
Now consider a function $I(u,v,\epsilon)$  obtained by replacing $s\mapsto s-i\epsilon$ 
in the boxed term, $\epsilon>0$. The new triple integral absolutely converges,
we can change the order of integrations and explicitly integrate in $z$.
We get
$$
I(u,v,\epsilon)= 
  \int\limits_\R \int\limits_\R u_k(is)\ov{v_k(it)}  A\bigl(\tfrac12(k-is)\bigr)^{-1}
  A\bigl(\tfrac12(k+it)\bigr)^{-1} 
  \frac{2^{-is-\epsilon+it}}{-is-\epsilon+it} \,ds \,dt.
$$
Next, we claim that 
$$
I(u,v)=\lim_{\epsilon\to+0}
I(u,v,\epsilon).
$$
 Indeed, we 
integrate $I(u,v,\epsilon)$ two times by parts in $s$ and come to
\begin{multline*}
I(u,v,\epsilon)=
 \int\limits_{|z|\ge 2}
  \int\limits_\R \ov{v_k(it)}  A\bigl(\tfrac12(k+it)\bigr)^{-1}  \int\limits_\R 
  \frac {\partial^2}{\partial s^2}
  \biggl[u_k(is) A\bigl(\tfrac12(k-is)\bigr)^{-1}
 \biggr] 
  \times\\\times
  \frac{|z|^{-2-\epsilon-is+it}}{i^2 \ln^2|z|} \,ds \,dt
  \,\dd z
\end{multline*}
The new triple integral absolutely converges and is  continuous at $\epsilon=+0$. 

Thus we come to the so-called distribution  $\frac1{x-i\epsilon}$, see,
e.g., \cite{GSh}. Recall  the Sokhotski formula
\begin{equation}
\lim_{\epsilon\to +0}\int_\alpha^\beta \frac{f(y)\,dy}{x-y-i\epsilon}
= \mathrm{p.v.}\int_\alpha^\beta \frac{f(y)\,dy}{x-y}+\pi i f(x),
\label{eq:Sokh}
\end{equation}
where $\mathrm{p.v.}$ denotes the principal value of an integral.

Applying this formula and keeping in  mind (\ref{eq:|A|}),
we come to
\begin{multline}
 I(u,v)=  \mathrm{p.v.} \int\limits_\R \int\limits_\R u_k(is)\ov{v_k(it)} \,
 A\bigl(\tfrac12(k-is)\bigr)^{-1}
  A\bigl(\tfrac12(k+it)\bigr)^{-1} 
  \frac{2^{-is+it}}{-is+it} \,ds \,dt
  +\\+
  \pi\int_\R u_k(is)\ov{v_k(is)}\, \kappa_{a,b}\bigl(\tfrac12(k+is)\bigr)\,ds.
  \label{eq:vpdelta1}
\end{multline}

 {\it Step} 8. We deal with $V_{22}$ in the same way and come to
 \begin{multline}
 V_{22}=  \mathrm{p.v.} \int\limits_\R \int\limits_\R u_k(is)\ov{v_k(it)} 
 A\bigl(\tfrac12(k-it)\bigr)^{-1}
  A\bigl(\tfrac12(k+is)\bigr)^{-1} 
  \frac{2^{-it+is}}{-it+is} \,ds \,dt
  +\\+
  \pi\int_\R u_k(is)\ov{v_k(is)} \kappa_{a,b}\bigl(\tfrac12(k+is)\bigr)\,ds
  + \Bigl\{\text{a $C^\omega$-smooth term}\Bigr\}.
  \label{eq:vpdelta2}
\end{multline}
Next, we take the sum $V_{11}+V_{22}$ modulo $C^\omega$-smooth terms. The expression
$$
\frac{ A\bigl(\tfrac12(k-is)\bigr)^{-1}
  A\bigl(\tfrac12(k+it)\bigr)^{-1} 2^{-is+it}-
  A\bigl(\tfrac12(k-it)\bigr)^{-1}
  A\bigl(\tfrac12(k+is)\bigr)^{-1} 2^{-it+is}
  }
  {-i(s-t)} 
$$
has the form
$$
\frac{L(t,s)-L(s,t)}{s-t}
$$
with analytic $L(t,s)$.
It 
has a removable singularity on the line $t=s$. 
Thus the first summands in (\ref{eq:vpdelta1}) and (\ref{eq:vpdelta2})  give
us a $C^\omega$-smooth term, the second summands give us 
the first term in (\ref{eq:o1}), i.e., the desired delta-function.

\sm

{\bf \punct End of the proof of Lemma \ref{l:orthogonal2}.}

\sm

{\it Step} 9. Next, we examine the term $V_{12}$. 
We write the integral and apply the transformation
(\ref{eq:mu-sum}). We get a sum of a $C^\omega$-smooth term and
the integral
\begin{multline*}
 J(u,v)=\int_{|z|\ge 2}\int_\R\int_\R
 u_k(is)\ov{v_k(it)} A\bigl(-\tfrac 12 (k+is)\bigr)^{-1} A\bigl(-\tfrac 12 (k+it)\bigr)^{-1}
 \times\\\times
 \Bigl(\frac z{\ov z}\Bigr)^k\, \boxed{|z|^{-2-is-it}}\, \dd z \,dt\,ds.
\end{multline*}
As above, we change $s\mapsto s-i\epsilon$ in the box and get integrals
$J(u,v,\epsilon)$ with $\epsilon>0$. As above,
\begin{multline*}
J(u,v;\epsilon)=
 \int_\R\int_\R
 u_k(is)\ov{v_k(it)} A\bigl(-\tfrac 12 (k+is)\bigr)^{-1} A\bigl(-\tfrac 12 (k+it)\bigr)^{-1}
 \times\\\times
\biggl[ \int_{|z|\ge 2}
\Bigl(\frac z{\ov z}\Bigr)^{-k} |z|^{-2-\epsilon-is-it}\, \dd z \biggr]
 \,dt\,ds.
\end{multline*}
If $k>0$, then the term in square brackets is zero (we pass to polar coordinates
and get 0 after the integration with respect to the angle coordinate).
If $k=0$, then
we get
$$
\frac{2^{-\epsilon-is-it}}
{\epsilon+ i(s+t)}.
$$
However,
$\supp(u_0)$, $\supp(v_0)$ are contained in domains $s>0$, $t>0$,
and 
actually we have no singularity.
Thus $V_{12}$ is $C^\omega$-smooth.

The same examination shows $C^\omega$-smoothness of $V_{21}$.
This completes the proof of Lemma \ref{l:orthogonal2}.
\hfill $\square$

\sm



\section{Asymptotics of the kernel in the parameters\label{s:asymptotic}}

\COUNTERS

{\bf\punct The statement.}
Let us modify a notation for the kernel $\cK$. Set
\begin{multline*}
\cK^\circ(z;\lambda; \sigma):=
\\:=
\frac1{\Gamma^\C(a+b|a+b)}\,
\FF\left[\begin{matrix}
    a+\lambda+\frac\sigma 2|a-\ov\lambda+\frac\sigma 2,\, a-\lambda-\frac\sigma 2|a+\ov\lambda-\frac\sigma 2
    \\ a+b|a+b
         \end{matrix}; z
  \right]
  =\\=
  \frac1{\Gamma^\C(a+b|a+b)}\,
 \FF\left[\begin{matrix}
    a+\frac{k+\sigma+is} 2|a+\frac{-k+\sigma+is} 2,\, a+\frac{-k-\sigma-is} 2|a+\frac{k-\sigma-is} 2
    \\ a+b|a+b
         \end{matrix}; z
  \right] ,
\end{multline*}
where $\lambda\in\Lambda$, $\sigma\in\R$.
In fact,
$$
\cK^\circ\bigl(z;\tfrac{k+is}2;\sigma\bigr)
=\cK(z; k, \sigma+is).
$$
 However, in calculations of this section the variables $\sigma$ and $is$ have  different meanings.

Denote
$$
t_\pm(z)=1\pm \sqrt{1-1/z}.
$$

\begin{theorem}
\label{th:asymptotic}
Then for a fixed $z$ we have the following asymptotic expansion
  \begin{multline}
 \cK^\circ(z;\lambda;\sigma)= \frac1 {\Gamma^\C\bigl(a-\lambda-\frac\sigma 2  |a+\ov \lambda-\frac\sigma 2\bigr)\,
 \Gamma^\C\bigl(b+\lambda+\frac\sigma 2 |b-\ov \lambda+\frac\sigma 2\bigr) \cdot |\lambda|}
\times \\\times|1-1/z|^{-1/2}\cdot |1-z|^{b-a}\cdot |z|^{-a-b}
\times\\\times
\biggl[ \Bigl(\frac{t_-(z)}{t_+(z)}\Bigr)^{\frac\sigma 2 +\lambda|\frac\sigma 2 -\ov\lambda}
\sum_{k\ge 0,\,l\ge 0,\, k+l< N}
\frac{\ov\lambda^{-k} \lambda^{\,-l}}
{k!\, l!} A_k(\sigma,\sqrt{1-z})\, {A_l(\sigma,\sqrt{1-\ov z})}
+\\+
\Bigl(\frac{t_+(z)}{t_-(z)}\Bigr)^{\frac\sigma 2 +\lambda|\frac\sigma 2 -\ov\lambda}
\!\!\!\!\!\!
\sum_{k\ge 0,\,l\ge 0,\, k+l< N}
\frac{\ov \lambda^{-k} \lambda^{\,-l}}
{k!\, l!} A_k(\sigma,-\sqrt{1-z})\, {A_l(\sigma,-\sqrt{1-\ov z})}
\biggr]+ \\+ R_N(z, \sigma,\lambda), 
\label{eq:}
 \end{multline}
 where $A_k(\xi)$ are rational expressions in $\xi$ {\rm(}depending on the parameters $a$, $b${\rm )} having poles at $\xi=0$, $\pm 1$
 and
 $A_0=1$.
The reminder $R_N(z)$ satisfies 
 $$
 R_N(z,\sigma,\lambda)= O(|\lambda|^{-N})
, \qquad \text{\rm as $\lambda\to\infty$},
 $$
moreover $O(\cdot)$ is uniform in $z$ and $\sigma$ on compact subsets in $\dot\C\times\R$. 
\end{theorem}

The proof occupies the rest of this section.

\sm

{\sc Remark.}
This formula is a counterpart of Watson's \cite{Wat} formula for asymptotics of the
Gauss hypergeometric functions $_2F_1[a-\lambda,b+\lambda;c;z]$ in the parameter $\lambda$
(see an exposition of Watson's results in \cite{Luke}, Sect. 7.2, see also a remark
in \cite{Olv}, p.162, 
on  typos in \cite{Wat}). We do not see
a way to reduce our statement to Watson's work.
\hfill $\boxtimes$

{\sc Remark.}
Lemma \ref{l:beta-as} gives us an asymptotics of the gamma-factor in (\ref{eq:}).
\hfill $\boxtimes$

\sm

{\bf\punct Stationary phase approximation.}
We transform  $\cK^\circ(z,\lambda,\sigma)$
as 
\begin{equation}
\tfrac1 {\Gamma^\C\bigl(a-\lambda-\frac\sigma 2  |a+\ov \lambda-\frac\sigma 2\bigr)\,
 \Gamma^\C\bigl(b+\lambda+\frac\sigma 2 |b-\ov \lambda+\frac\sigma 2\bigr) \cdot |\lambda|}
\int_\C
R(t ,z,\sigma) \exp\Bigl\{ Q(t,z,\lambda,\sigma)   \Bigr\}
\dd t,                  
\label{eq:st1}
\end{equation}
where
\begin{equation}
R(t ,z):=t^{a-\frac\sigma2-1|a-\frac\sigma2-1} (1-t)^{b+\frac\sigma2-1|b+\frac\sigma2-1} (1-tz)^{-a-\frac\sigma2}
\end{equation}
and
\begin{multline}
Q(t,z,\lambda):= 
\lambda \ln \Bigl(\frac{t(1-zt)}{1-t} \Bigr)- \ov\lambda\, \ov{\ln \Bigl(\frac{t(1-zt)}{1-t} \Bigr)}
=\\=
ik \Im \ln \Bigl(\frac{t(1-zt)}{1-t} \Bigr)+i s \ln \Bigl(\frac{t(1-zt)}{1-t} \Bigr).
\label{eq:St2}
\end{multline}
The function $\Im \ln (\dots)$ is ramified, however the exponent is well-defined
and formulas below contain only partial derivatives of $\ln(\dots)$, which are independent of the choice
of a branch.

We apply the stationary phase approximation, see, e.g., Fedoryuk \cite{Fed}, H\"ormander \cite{Hor}.
Singular points are 
$0$, $1$, $\infty$.
Stationary 
points are
$$
t_\pm=1\pm \sqrt{1-1/z},
$$
they are the same for  both summands in (\ref{eq:St2}). This could be a fatal obstacle for an evaluation
of a uniform asymptotics, however this does not happen.  
Also the domain of convergence of the integral (\ref{eq:st1}) is smaller than it is necessary for 
our purposes. 

Consider a partition of unity
$$
1=\rho_0+\rho_1+ \rho_{z^{-1}} + \rho_\infty+ \rho_{t_+}+\rho_{t_-}+  \tau,
$$
where $\rho_\alpha$ is zero outside a small neighborhood of $\alpha$, and 
$\tau$ is zero in neighborhoods of $0$, $1$, $z^{-1}$, $\infty$, $t_\pm$.
According to this partition we expand (\ref{eq:st1}) into a sum of 7 integrals,
$$
I=I_0+I_1+I_{z^{-1}} +I_\infty+ I_{t_+}+I_{t_-}+ J.
$$

Obviously (see \cite{Fed}, Lemma III.2.1), for each $N$ we have
$$
J=O(k^2+s^2)^{-N}\qquad \text{as $n+is\to\infty$}.
$$

{\bf\punct  Preparatory statement.}

\begin{theorem}
\label{th:parametric}
Let $\Omega$ be a domain in $\C$, $f(t)$, $\phi(t)$ be holomorphic in $\Omega$. Let $t_0$ be
a unique zero of
$\phi'(t)$ in $\Omega$ and $\phi''(t_0)\ne 0$.
Let $\rho(t)$ be a $C^\infty$-smooth function compactly supported by $\Omega$ such that
$\rho=1$ in a neighborhood of $t_0$. 
Consider the integral 
\begin{equation}
\label{eq:I-lambda}
I(\lambda)=\int_\Omega\rho(t)\, f(t)\,\ov {f(t)} \exp\Bigl\{ i \Re (\lambda \phi(t) \Bigr\}\,\dd t
,
\end{equation}
where $\lambda\in\C$ is a parameter. Then

\sm

{\rm a)} For $|\lambda|>1$ we have the following  expansion
\begin{multline}
 I(\lambda)= \frac1{|f''(t_0)|\, |\lambda|}\exp\Bigl\{ i  \Re (\lambda \phi(t_0) \Bigr\}
 \times\\\times
 \Bigl(\sum_{k\ge 0,\,j\ge 0,\, k+l< N} \frac{\lambda^{-k}\ov\lambda^{\,-l}}{k!\,l!}
 a_k (f,\phi)\, a_l(\ov f,\ov \phi)+R_N(\lambda)
 \Bigr),
\end{multline}
where $a_k$ are rational expressions 
$$a_k=a_k \bigl(\phi(t_0), \phi'(t_0), \dots;  f(t_0), f'(t_0), \dots;
\phi''(t_0)^{-1}\bigr)
$$
and  $a_0=1$. The reminder $R_N$
satisfies
\begin{equation}
R_N(\lambda)= O(|\lambda|^N) \qquad \text{as $\lambda\to\infty$.}
\label{eq:RN}
\end{equation}

{\rm b)} The asymptotic expansion  
$$
I(\lambda)\sim |\lambda|^{-1} \sum_{k\ge 0,\,l\ge 0} \frac{c_{kl}}{\lambda^k \ov\lambda^{\,l}}
\qquad \text{as $\lambda\to\infty$}
$$
can be written as
\begin{multline}
I(\lambda)\sim 
 \frac1{|f''(t_0)|\, |\lambda|}\exp\Bigl\{ i  \Re (\lambda \phi(t_0) \Bigr\}
 \times\\\times 
 \exp\Bigl\{ \frac i{2\ov \lambda \,\phi''(t_0)}\frac{\partial^2} {\partial t^2 }\Bigr\} 
 \Bigl(f(t)
 \exp\Bigl\{\ov\lambda\bigl(\phi(t)-\phi'(t_0)-\tfrac 12 \phi''(t_0) (t-t_0)^2\bigr)\Bigr\}\Bigr)\biggr|_{t=t_0}
 \times\\\times
  \exp\Bigl\{ \frac i{2 \lambda\, \ov{\phi''(t_0)}}\frac{\partial^2} {\partial \ov t^2 }\Bigr\} 
 \Bigl(f(t)
 \exp\Bigl\{\lambda\bigl(\ov{\phi(t)-\phi'(t_0)-\tfrac 12 \phi''(t_0) (t-t_0)^2}\bigr)\Bigr\}\Bigr)\biggr|_{t=t_0}.
\end{multline}

{\rm c)} Let $\phi=\phi_\alpha$, $f=f_\alpha$ smoothly depend on a parameter $\alpha$,
where $\alpha$ ranges in a compact domain  $K\subset \C$ and the conditions of the preamble
of the theorem are satisfied for all $\alpha$.
 Then $O(\cdot)$ in {\rm (\ref{eq:RN})} is uniform
in $\alpha\in K$.
\end{theorem}

{\sc Proof.} b)
We use Fedoryuk \cite{Fed}, Proposition III.2.2 or H\"ormander \cite{Hor}, Theorem 7.7.5.
Let $f$ be a smooth compactly supported function on $\R^n$, let $S$ be smooth.
Consider an $n$-dimensional integral 
$$
I(\sigma):=\int_{\R^n} f(x) \exp\{i \sigma S(x)\} \,d x, \qquad t\ge 1.
$$
Let $x_0$ be a unique critical point of $S$ on the support of $f$, let it be nondegenerate.
Let $H(x_0)$ be the Hessian of $S$ at $x_0$ (i.e., the matrix composed of second partial derivatives),
let $\sgn H(x_0)$ denote the signature of the Hessian
(the number of positive eigenvalues minus the number of negative eigenvalues). Consider
the second order differential operator 
$$
L:=\frac i2\la H(x_0)^{-1} \nabla_x,\nabla_x\ra,
$$
where $\nabla_x$ denotes the vector column composed of $\frac\partial {\partial x_1}$, \dots, 
$\frac\partial {\partial x_n}$. Denote
\begin{equation}
S(x,x_0):=S(x)- S(x_0)-\frac 12 \la H(x_0) (x-x_0), (x-x_0)\ra,
\label{eq:Sxx}
\end{equation}
this expression is  the part of the Taylor expansion of $S(x)$ at $x_0$ starting cubic terms.
Then the following 
 expansion take holds:
\begin{multline}
I(\sigma)=
\Bigl(\frac{2\pi} \sigma \Bigr)^{n/2} |\det H(x_0)|^{-1/2} \exp\Bigl[\frac {i\pi} 4\sgn H(x_0) \Bigr]
\times\\\times
\biggl(
\sum_{k=0}^{N-1} \frac{\sigma^{-k}}{k!} \,L^k\bigl(f(x) \exp\{ i\sigma S(x,x_0) \}\bigr)\Bigr|_{x=x_0}
+
\sigma^{-N+[ 2N/3 ]}
V(\sigma)
\biggr),
\label{eq:fedor}
\end{multline}
where $V(\sigma)$ is bounded.

\sm

Let us return to our integral (\ref{eq:I-lambda}).
Without loss of generality, we can set $t_0=0$, $\phi''(t_0)=1$, i.e.,
$$
\phi(t)=\tfrac12 t^2+ r(t), \qquad \text{where $r(0)=r'(0)=r''(0)=0$.}
$$
Set $\lambda=s e^{i\theta}$, $s>0$. Set $z=x+iy$, then
$$
\phi(x,y)=\tfrac12(x^2-y^2+2ixy) +r(x,y).
$$
Thus we come to an oscillating integral in $s$ with  the parameter $\theta$, 
$$
I(s,\theta)=
\int \rho(x,y) f(x,y) \ov{f(x,y)} \exp\Bigl\{ i s(\cos\theta \Re \phi(x,y)+ \sin\theta \Im \phi(x,y)) \Bigr\}\,dx\,dy.
$$
We wish to apply the general statement formulated above. The Hessian 
is given by
$$
H= 2
\begin{pmatrix}
 \cos\theta & \sin\theta\\
 \sin\theta& -\cos\theta
\end{pmatrix},\qquad 
H^{-1}= \frac 12
\begin{pmatrix}
 \cos\theta & \sin\theta\\
 \sin\theta& -\cos\theta
\end{pmatrix}.
$$
The signature is 0.
The differential operator $L$ is
$$
L=\frac i4\Bigl(\cos \theta\Bigl(\frac{\partial^2}{\partial x^2}-\frac{\partial^2}{\partial y^2}\Bigr)+2\sin \theta
\frac{\partial^2}{\partial x\partial y}\Bigr)
= \frac i2\Bigl(e^{i\theta} \frac{\partial^2}{\partial t^2}+ e^{-i\theta} \frac{\partial^2}{\partial \ov t^2} \Bigr).
$$
Next, we rewrite our phase function $S(\cdot)$
as
$$
e^{-i\theta} \phi(t)+e^{i\theta}\ov{\phi(t)}.
$$
Therefore the expression (\ref{eq:Sxx}) is
$$
e^{-i\theta} r(t)+e^{i\theta}\ov{r(t)}.
$$
Applying (\ref{eq:fedor}), we
get
\begin{multline*}
I(s,\theta):=
\frac {2\pi} s
\exp\Bigl\{ \frac i{2s} \Bigl(e^{i\theta} \frac{\partial^2}{\partial t^2}+
e^{-i\theta} \frac{\partial^2}{\partial \ov t^2} \Bigr) \Bigr\}
\times\\\times
\Bigl(f(t)\ov{f(t)} \exp\bigl\{is\bigl(e^{-i\theta} r(t)+e^{i\theta}\ov{r(t)}\bigr)\bigr\}\Bigr)
\biggr|_{t=0}
=\\=
\frac {2\pi} s
\exp\Bigl\{ \frac i{2s e^{-i\theta}}  \frac{\partial^2}{\partial t^2}\Bigr\}
\Bigl(f(t) \exp\bigl\{is e^{-i\theta} r(t)\bigr\}\Bigr)\biggr|_{t=0}
\times\\\times
\exp\Bigl\{ \frac i{2s e^{i\theta}}  \frac{\partial^2}{\partial \ov t^2}\Bigr\}
\Bigl(\ov{f(t)} \exp\bigl\{is e^{i\theta} \ov{r(t)}\bigr\}\Bigr)\biggr|_{t=0}.
\end{multline*}
We obtained  asymptotics in $s$ for fixed $\theta$.
However, $\theta$ ranges in a compact set, by \cite{Fed}, Theorem III.2.2,
we get that the term $V(\cdot)$ in (\ref{eq:fedor})
is bounded uniformly in $\theta$. 

\sm

a) follows from b).

\sm

c) We again refer to the parametric version of the stationary phase approximation, see
\cite{Fed}, Theorem III.2.2.
\hfill $\square$

\sm

{\bf \punct Contribution of the stationary points.} Let us apply 
Theorem \ref{th:parametric} to our integral (\ref{eq:st1}).
 We have
\begin{align*}
f(t)&=R(t,z)=
\Bigl(\frac{t}{1-zt}\Bigr)^{a} 
\, (1-t)^{b}\, \Bigl(\frac{1-t}{t(1-zt)}\Bigr)^{\frac\sigma2}\, \bigl(t (t-1)\bigr)^{-1};
\\
\phi(t)&= 2 \ln \Bigl(\frac{1-t}{t(1-zt)} \Bigr).
\end{align*}

Denote 
$$
\zeta=\sqrt{1-1/z}.
$$
We substitute $t=t_+$ and transform the factors of $R(t,z)=f(t)\ov{f(t)}$:
\begin{align}
\Bigl(\frac{t}{1-zt}\Bigr)^{a|a}\biggr|_{t=t_+}& = \Bigl(\frac{1-\zeta^2}{-\zeta}\Bigr)^{a|a}
=\bigl((z-1)z\bigr)^{-a/2|-a/2};\\
(1-t)^{b|b} \Bigr|_{t=t_+}
&= \Bigl(\frac{1-z}{z}\Bigr)^{b/2|b/2};
\\
\Bigl(\frac{1-t}{t(1-zt)}\Bigr)^{\frac\sigma2}\Bigr|_{t=t_+}&=\Bigl(\frac{1-\zeta}{1+\zeta}\Bigr)^{\frac \sigma2|\frac \sigma2}
\\
\bigl(t (t-1)\bigr)^{-1|-1}\Bigr|_{t=t_+}
&=\Bigl(\frac{-1}{\zeta(1+\zeta)} \Bigr)^{1|1}.
\end{align}
Next, 
\begin{equation*}
 \phi(t_+)=2\ln\Bigl(\frac{1-\zeta}{1+\zeta}\Bigr),
\end{equation*}
therefore
\begin{equation*}
 \exp\Bigl\{ i\Re \bigl(\phi(t_+) \tfrac 12(k+is)\bigr)\Bigr\}
=
 \Bigl(\frac{1-\zeta}{1+\zeta}\Bigr)^{\lambda|-\ov\lambda}=
 \Bigl(\frac{t_-}{t_+}\Bigr)^{\lambda|-\ov\lambda}.
\end{equation*}
Finally,
$$
\phi''(t)=
\frac{-2}{(1 - t)^2} + \frac 2{t^2} + \frac{2 z^2}{(1 - t z)^2},
$$
and 
$$
\phi''(t_+)= \frac{-4}{\zeta(1 + \zeta)^2}.
$$
Uniting these data we get that the leading term at the point $t_+$ is
\begin{equation}
- |\zeta|\, |1-z|^{b-a}\, |z|^{-a-b}
\Bigl(\frac{t_-}{t_+}\Bigr)^{\lambda|-\ov\lambda} \cdot \frac 1{(k^2+s^2)^{1/2}}.
\end{equation}

Th general form of the asymptotic expansion at $t=t_+$ follows from Theorem \ref{th:parametric}.

\sm

{\bf \punct Contributions of the singular points.}
\begin{lemma}
 Contributions at the singular points $0$, $1$, $\infty$ are 
 $O(|\lambda|^{-N})$ for any $N$.
\end{lemma}

{\sc Proof.}
For definiteness  examine
the point 0.
We have the integral 
$$
I_0(\lambda)=\int_\C \rho_0(t)\, t^{a-1|a-1} (1-t)^{c-a-1} (1-zt)^{-a}
\Bigl(\frac{t(1-zt)}{1-t}\Bigr)^{\lambda|-\ov\lambda}\dd t,
$$
defined as an analytic continuation.
Keeping in  mind that a support of $\rho_0$ can be chosen sufficiently small,
we pass to a new variable in a neighborhood of $0$,
$$
u=\frac{t(1-zt)}{1-t},
$$
and come to an integral of the form 
$$
I_0(\lambda)=\int_\C u^{a-\lambda-1|a+\ov \lambda-1} \Phi(u)\dd u,
$$
where $\Phi$ is a smooth compactly supported function. It remains to apply Theorem \ref{th:mellin-meromorphic}.

\sm

Argumentation for other singular points is the same. 
\hfill $\square$

\section{Symmetry of difference operators\label{s:symmetry2}}

\COUNTERS

Here we prove Theorem \ref{th:paley}, i.e., show that  if $f\in\cD(\dot\C)$,
then
$J_{a,b}f$ is contained in the space $\cW_{a,b}$ of meromorphic functions 
on $\Lambda_\C$.
Also we show that 
$\frL$ and $\ov \frL$ are formally adjoint one to another on $\cW_{a,b}$,
see Theorem \ref{th:symmetry}.

\sm

{\bf\punct Beginning of the proof of Theorem \ref{th:paley}.}
We follow  the list of properties in the definition of $\cW_{a,b}$,
see Subsect. \ref{ss:complex-domain},

\sm

 a)  is a corollary of the  symmetry
$
\cK_{a,b}(z; -k, -\sigma)=\cK_{a,b}(z; k, \sigma)
$.

\sm

b) We must examine poles of $\cK_{a,b}(z; k, \sigma)$ as a function of the variable $\sigma$ for a fixed
$z\in\dot\C$,
$k\in\Z$.
Let $a+b\ne1$. We look to the expansion (\ref{eq:A}) of $\FF[\cdot]$
at $z=0$. The only source of  poles
of $\cK$
 are zeros of the denominators
in (\ref{eq:A11}), i.e., zeros of the expression
\begin{multline}
 R(k,\sigma):=\Gamma^\C  \bigl(a+\tfrac{k+\sigma}2\bigl| a+\tfrac{-k+\sigma}2\bigr)
 \,
 \Gamma^\C  \bigl(a+\tfrac{-k-\sigma}2\bigl| a+\tfrac{-k-\sigma}2\bigr)
 \times\\
 \Gamma^\C  \bigl(b+\tfrac{k+\sigma}2\bigl| b+\tfrac{-k+\sigma}2\bigr)
 \,
 \Gamma^\C  \bigl(b+\tfrac{-k-\sigma}2\bigl| b+\tfrac{-k-\sigma}2\bigr).
 \label{eq:Rks}
\end{multline}
This gives us the desired list of possible poles.

Let us examine  the case $a+b=1$. The decomposition of the
hypergeometric functions  (\ref{eq:A}) at $z=0$ 
produces an expression
of the type
\begin{equation}
\cK_{a,b}(z; k, \sigma)=
\frac{u_{a,b}(z,k,\sigma)-v_{z,a,b}(z,k,\sigma)}{a+b-1}
\label{eq:a+b-1}
\end{equation}
with $u_{a,b}$, $v_{a,b}$ having poles at
zeros of $R(k,\sigma)$. A decomposition at $z=1$
gives 
$$
\cK_{a,b}(z; k, \sigma)=
\frac{U_{a,b}(z,k,\sigma)-V_{a,b}(z,k,\sigma)}{a-b},
$$
therefore  the singularity in (\ref{eq:a+b-1})
at $a+b=1$ is removable.

\sm

d) Indeed, we have $ \cK_{a,b} (p,q)=  \cK_{a,b} (q,p)$, i.e.,
\begin{multline*}
\FF\left[\begin{matrix}
a+\frac{p+q}2\bigl| a+\frac{-p+q}2,\, a+\frac{-p-q}2\bigl| a+\frac{p-q}2 
\\
a+b\bigl|a+b
         \end{matrix}; z
\right]
=\\=
\FF\left[\begin{matrix}
a+\frac{p+q}2\bigl| a+\frac{p-q}2,\, a+\frac{-p-q}2\bigl| a+\frac{-p+q}2 
\\
a+b\bigl|a+b
         \end{matrix}; z
\right].
\end{multline*}
This is a special case of the symmetry (\ref{eq:additional}).

We also mention the following similar identity for (\ref{eq:Rks}):
\begin{equation}
R(p,q)=R(q,p),
 \label{eq:RR}
\end{equation}
it is a special case of (\ref{eq:R-additional}). 

\sm

The statement c) about the behavior at infinity is a corollary of the expansion (\ref{eq:}) and the following lemma 

\begin{lemma}
\label{l:asy}
Let $t_\pm(z)$ be as in Theorem {\rm \ref{th:asymptotic}}. Let 
 $\Phi\in \cD(\dot \C)$ be a function with a  simply connected support.
 Then for any $A>0$ for any $N>0$ in the strip $|\Re \sigma|<A$ we have 
$$
\int_{\dot\C}
\Phi(z) \Bigl(\frac{t_-(z)}{t_+(z)}\Bigr)^{(k+\sigma)/2|(-k+\sigma)/2}\,\dd{z}
=
O\bigl(k^2+(\Im \sigma)^2\bigr)^{-N}
$$
as $(k^2+(\Im \sigma)^2)\to \infty$.
\end{lemma}

We need  a simply connected support since the integrand is ramified at the points $z=0$, $z=1$.
A proof of the lemma requires some preparations.

\sm

{\bf \punct A change of variable.}
We define a new variable
\begin{equation}
p:= \frac{t_+(z)}{t_-(z)},
\label{eq:p}
\end{equation}
The inverse map 
is done by
\begin{equation}
z=
\zeta(p):=\frac{(p+1)^2}{4p}.
\label{eq:zeta}
\end{equation}
The map $\zeta(p)$ determines a two-sheet covering map
from
\begin{equation}
\ddot\C:=\C\setminus \{0,1,-1\}
 \label{eq:ddot}
\end{equation}
to $\dot\C$.
Notice that
\begin{align}
1-z&=-\frac{(p-1)^2}{4p}, & \sqrt{1-1/z}&=\frac{p-1}{p+1},&
\zeta'(p)&=\frac{p^2-1}{4p^2},
\label{eq:zeta1}
\\
t_+&=\frac{2p}{p+1},& t_-&=\frac{2}{p+1},&
\frac{t_+}{t_-}&=p. 
\label{eq:zeta2}
\end{align}
Also,
\begin{equation}
\zeta(p^{-1})=\zeta(p),\qquad \zeta'(p^{-1})\, p^{-1}=\zeta(p)\, p.
\label{eq:chet}
\end{equation}

{\bf\punct Proof of Theorem \ref{th:paley}.c.}

\sm

{\sc Proof of Lemma \ref{l:asy}}
We substitute $z=\zeta(p)$ to the integral and get
$$
\frac 1{16}
\int_{\ddot\C}  p^{(k+\sigma)/2|(-k+\sigma)/2} \Bigl( \Phi(\zeta(p))
\,|p^2-1|^2\, p^{-2}\Bigr)\,\dd{p}.
$$
This is a Mellin transform of a function compactly supported by $\ddot \C$.
In virtue of Theorem
\ref{th:mellin-meromorphic} the integral rapidly decreases in the union of strips $|\Re \sigma|<A$.

\sm

{\sc Proof of the statement c) of Theorem \ref{th:paley}.}
We represent $\phi(z)$ as a sum of functions in $\cD(\dot\C)$ with simply connected supports.
Next, we decompose the kernel according to Theorem \ref{th:asymptotic} and apply the lemma to each summand.

\sm

{\bf\punct Continuity.}

\begin{corollary}
\label{cor:JCL2}
 The map $J_{a,b}$ is a continuous map from $\cD(\dot\C)$
 to $L^2_\even(\Lambda, \kappa_{a,b})$.
\end{corollary}

{\sc Proof.} Define the following seminorms on the space of smooth functions on $\Lambda$:
$$
p_{\alpha,N} (F)=\sup_{\lambda\in\Lambda}
 \Bigl| \frac{\partial^N F}{\partial \sigma^N} (1+|\lambda|)^\alpha\Bigr|,
$$
and the space $\cY$ defined by these seminorms. Clearly, our proof provides a continuity
of $J_{a,b}$ as a map $\cD(\dot\C)$ to $Y$. 
It remains to notice that the identical embedding $f\mapsto f$ of $Y$ to $L^2$ is continuous.

If $k=0$ and  $a=1$ (or $b=1$), then elements of $\cW_{a,b}$ have a pole of order two%
\footnote{At the same point the spectral density has a zero of order 4.}
at
$k=0$, $s=0$. In this case we write $\lambda^2 F$ instead of $F$ in the definition of the seminorms.
\hfill $\square$

\sm

{\bf \punct Invariance of $\cW_{a,b}$.}
Consider the difference operators $\frL$, $\ov\frL$   defined
above (\ref{eq:frL}), 
\begin{multline}
\frL F(k,\sigma)= \frac{(a+\frac{k+\sigma}2) (b+\frac{k+\sigma}2)}{(k+\sigma)(1+k+\sigma)} \bigl(F(k+1,\sigma+1)-F(k,\sigma)\bigr)
+\\+
\frac{(a+\frac{-k-\sigma}2) (b+\frac{-k-\sigma}2)}{(-k-\sigma)(1-k-\sigma)} \bigl(F(k-1,\sigma-1)-F(k,\sigma)\bigr);
\label{eq:frL-bis}
\end{multline}
\begin{multline}
\ov\frL F(k,\sigma)= \frac{(a+\frac{-k+\sigma}2) (b+\frac{-k+\sigma}2)}{(-k+\sigma)(1-k+\sigma)} \bigl(F(k-1,\sigma+1)-F(k,\sigma)\bigr)
+\\+
\frac{(a+\frac{k-\sigma}2) (b+\frac{k-\sigma}2)}{(k-\sigma)(1+k-\sigma)} \bigl(F(k+1,\sigma-1)-F(k,\sigma)\bigr).
\end{multline}

\begin{lemma}
The space $\cW_{a,b}$ is invariant with respect to 
 the operators $\frL$, $\ov\frL$.
\end{lemma}

{\sc Proof.}
Since $F(0,-1)=F(1,0)= F(-1,0)=F(0,1)$,  the expressions
$$
\frac{F(k+1,\sigma+1)-F(k,\sigma)}{1+k+\sigma},\qquad \frac{F(k-1,\sigma-1)-F(k,\sigma)}{1-k-\sigma}
$$
have no poles at $k=-1$, $\sigma=0$ and $k=1$, $\sigma=0$ respectively.

Since a function $F(k,\sigma)$ is even, it can not have a pole of order 1 at $k=0$, $\sigma=0$.

New poles of $F(k+1,\sigma+1)$ that are not poles of $F(k,\sigma)$ are annihilated by the rational factor in (\ref{eq:frL-bis}).

The condition $\frL F(p,q)=\frL F(q,p)$ follows from a straightforward calculation.
\hfill $\square$

\sm

{\bf \punct Symmetry.}

\begin{theorem}
\label{th:symmetry}
For $(a,b)\in\Pi$,
 for  $F$, $G\in \cW_{a,b}$ we have
 \begin{equation}
 \label{eq:Lsym}
 \la \frL F, G\ra_{L^2(\Lambda_\C,dK_{a,b})}=
  \la F, \ov \frL  G\ra_{L^2(\Lambda,dK_{a,b})}.
  \end{equation}
\end{theorem}

\begin{corollary}
 Operators $\frac 12(\frL+\ov\frL)$, $\frac 1{2i}(\frL-\ov\frL)$
 are symmetric on the $J_{a,b}$-image of $\cD(\dot\C)$.
\end{corollary}

{\sc Remark.} In fact, the proof uses only properties of $F\in \cW_{a,b}$ in
strips $|\Re\sigma|<1+\epsilon$. So we can define operators $\frL$, $\ov\frL$ on a space of meromorphic functions
 in the strip satisfying   an obvious list of conditions.
\hfill $\boxtimes$

\sm

{\bf\punct Proof of Theorem \ref{th:symmetry} for the case $(a,b)\in\Pii$.}
First, we 
notice that for pure imaginary $\sigma$ we have 
$\ov {G(k,\sigma)}=\ov{G(k,-\ov \sigma)}$, the last function is meromorphic and also is contained in $\cW_{a,b}$.
Let $R(k,\sigma)$ be given by (\ref{eq:Rks}). Then 
\begin{multline}
4\pi^2 i\,
  \la \frL F, G\ra=
  \sum_k\int_{i\R} 
 \biggl\{ \frac{(a+\frac{k+\sigma}2) (b+\frac{k+\sigma}2)}{(k+\sigma)(1+k+\sigma)} \bigl(F(k+1,\sigma+1)-F(k,\sigma)\bigr)
+\\+
\frac{(a+\frac{-k-\sigma}2) (b+\frac{-k-\sigma}2)}{(-k-\sigma)(1-k-\sigma)} \bigl(F(k-1,\sigma-1)-F(k,\sigma)\bigr)\biggr\}
\times\\\times
\ov{G(k,-\ov \sigma)}\,\,
(k-\sigma)(k+\sigma) R(k,\sigma)\, d\sigma.
\label{eq:4summands}
\end{multline}

Let us expand the expression in the curly brackets $\{\dots\}$ as a sum of 4 summands that include
$F(k+1,\sigma+1)$, $F(k,\sigma)$, $F(k-1,\sigma-1)$, $F(k,\sigma)$. The whole expression $\{\dots\}$ is holomorphic near the contour of integration.
The summands have simple poles on the contour, and we  pass to an integration in the sense of principal values.

Let us examine the summand corresponding $F(k+1,\sigma+1)$. We get
\begin{equation}
  \sum_k \mathrm{v.p.}
  \int_{i\R} \frac{k-\sigma}{1+k+\sigma}
 F(k+1,\sigma+1)\,
\ov{G(k,-\ov \sigma)}\,
\wt R(k,\sigma)\, d\sigma,
\label{eq:intFG}
\end{equation}
where
\begin{multline}
 \wt R(k,\sigma):= \bigl(a+1+\tfrac{k+\sigma}2\bigr) \bigl(b+1+\tfrac{k+\sigma}2\bigr) R(k,\sigma)
 =\\=
 \Gamma^\C  \bigl(a+1+\tfrac{k+\sigma}2\bigl| a+\tfrac{-k+\sigma}2\bigr)
 \,
 \Gamma^\C  \bigl(a+\tfrac{-k-\sigma}2\bigl| a+\tfrac{-k-\sigma}2\bigr)
 \time\sigma\\
 \Gamma^\C  \bigl(b+1+\tfrac{k+\sigma}2\bigl| b+\tfrac{-k+\sigma}2\bigr)
 \,
 \Gamma^\C  \bigl(b+\tfrac{-k-\sigma}2\bigl| b+\tfrac{-k-\sigma}2\bigr).
 \label{eq:wtRks}
\end{multline}

\begin{lemma}
 For $0<a<1$, $0<b<1$ the integrand in {\rm(\ref{eq:intFG})}
 has no poles in the strip $-1<\Im \sigma<0$.
\end{lemma}

{\sc Proof.} We enumerate possible (simple) poles of the factors.

\sm

a) Factor $\ov {G(k,-\ov \sigma)}$. In this case we can have poles if $k=0$.
Since $a<1$, $b<1$ the poles $2-2a$, $2-2b$ are outside our strip.
On other hand  the pole $2a-2$ (resp. $2b-2$) is contained in the strip if $1/2<a<1$ (resp. if $1/2<b<1$). 

\sm

b) Factor $F(k+1,\sigma+1)$ has a pole in our strip for $k=-1$ at $\sigma=2a-1$ (resp. $\sigma=2b-1$) if $0<a<1/2$, 
(resp. $0<b<1/2$).

\sm

c) Since $a>0$, $b>0$ the expression $\wt R(k,\sigma)$ has no poles in our strip.

\sm

However, the poles of  $\ov {G(k,-\ov \sigma)}$ and
of $F(k+1,\sigma+1)$ are zeros of $\wt R(k,\sigma)$. Therefore
the product is holomorphic.
\hfill $\square$

\sm

\begin{lemma}
In {\rm(\ref{eq:intFG})}, we can change the integration contour to $1+i\R$.
\end{lemma}

{\sc Proof.} The integrand has no poles between contours $i\R$ and $1+i\R$, but has poles on contours,
the
integral is taken in the sense of principal values. We have only two such poles,
 $\sigma=0$ on the contour $i\R$ for $k=-1$ and $\sigma=-1$ for $k=0$.
 Thus the difference between the two integrals is $2\pi$ by half of the sum of residues, i.e.,
 \begin{multline*}
 \frac {2\pi} 2
 \Biggl\{ (-1-\sigma)
 F(0,\sigma+1) 
\ov{G(-1,-\ov \sigma)}
\bigl(a+1+\tfrac{-1+\sigma}2\bigr) \bigl(b+1+\tfrac{-1+\sigma}2\bigr)   R(-1,\sigma)\Bigr|_{\sigma=0}
 +\\+
 (0-\sigma)
 F(1,\sigma+1) 
\ov{G(0,-\ov \sigma)}
\bigl(a+1+\tfrac{\sigma}2\bigr) \bigl(b+1+\tfrac{\sigma}2\bigr)   R(0,\sigma)\Bigr|_{\sigma=-1}
\biggr\}.
 \end{multline*}
 Let us show that the sum is zero.
 Since $F$, $G$ are even and satisfy (\ref{eq:kk}), we have
 $$
  F(0,1)=F(1,0), \qquad \ov{G(-1,0)}= \ov{G(0,-1)}
 $$
 By (\ref{eq:RR}), we have
 $$
 R(-1,0)=R(0,-1).
 $$
 The remaining factors give
 $$
 -(a+\tfrac 12) (b+\tfrac 12) \quad\text{and}
 \quad
  (a+\tfrac 12) (b+\tfrac 12),
 $$
 i.e., the same expressions with different signs.
 \hfill $\square$
 
 \sm

{\sc End of the proof of Theorem \ref{th:symmetry} for the case $(a,b)\in\Pii$.}
Thus we can replace the integration in (\ref{eq:intFG})
by the integration over the contour $-1+i\R$.
We change the variables $l=k+1$, $t=\sigma+1$ and get
\begin{equation*}
 \sum_l 
   \mathrm{v.p.}
  \int_{i\R} \frac{l-t}{-1+l+t}
 F(l,t)\,
\ov{G(l-1,-\ov t+1)}\,
\wt R(l-1,t-1)\, dt.
\end{equation*}
Next,
$$
\wt R(l-1,t-1)=R(l,t) \bigl(a+\tfrac{-l-t}2\bigr) \bigl(b+\tfrac{-l-t}2\bigr),
$$
and we come to
\begin{multline*}
 \sum_l 
   \mathrm{v.p.}
  \int_{i\R}
 F(l,t)
\biggl[ \frac{\bigl(a+\tfrac{-l-t}2\bigr)\bigl(b+\tfrac{-l-t}2\bigr)}{(-l-t)(1-l-t)} \,\ov{G(l-1,-\ov t+1)}\biggr]
 \times\\\times
(l-t)(l+t)\, R(l,t)\, dt.
\end{multline*}
We transform the expression in the big brackets to the form
$\ov {U(l,-\ov t)}$ , where 
$$
U(l,t)=\frac{\bigl(a+\tfrac{-l+t}2\bigr)\bigl(b+\tfrac{-l+t}2\bigr)}{(-l+t)(1-l+t)}
\, G(l-1,t+1).
$$
Thus we finished the transformation of the summand of the (\ref{eq:4summands}) corresponding to
$F(k+1,\sigma+1)$. The transformation of the summand corresponding to $F(k-1,\sigma-1)$
is similar. The case of the summands $F(k,\sigma)$ is obvious. We come 
to the desired expression.
\hfill $\square$

\sm

{\bf\punct End of the proof of Theorem \ref{th:symmetry}.}
Due to the homographic transformations, it is sufficient to examine the case $a<0$.
Let $\Phi$, $\Psi\in \cW_{a,b}$.
Denote
\begin{equation}
U(a,b;k,\sigma):= 
\Phi(k,\sigma)\,\ov{\Psi(k,-\sigma)} \, \kappa_{a,b}(k,\sigma). 
\end{equation}
For $(a,b)\in \Pii$ we have
\begin{equation} 4\pi^2 i
\la \Phi,\Psi\ra_{L^2(\Lambda,\kappa_{a,b})}
=\sum_k\int_{i\R} U(a,b;k,\sigma)\,d\sigma.
\label{eq:PhiPsi}
\end{equation}
We wish to write the analytic continuation of this expression to the
domain $(a,b)\in \Pi$, $a\le 0$.

Possible singularities of $U$ as a function in $\sigma$
in the strip $|\Re \sigma|<1$ are the following:

\sm

--- if $b>1/2$, then  both functions $\Phi$, $\Psi$ have  poles at $(k,\sigma)=(0,\pm(2-2b)$;

\sm

---  $\kappa_{a,b}(k,\sigma)$ has poles at $(k,\sigma)=(0,\pm 2a)$.

\sm

Due to our restrictions $2b-2<2a<-2a<2-2b$.

Thus all summands of (\ref{eq:PhiPsi}) except
$0$-th are holomorphic in $|a|<1-b$.

\begin{lemma}
	\label{l:contour}
	 Fix $b$.
 Assume that $\Phi$, $\Psi$ be even rapidly decreasing meromorphic functions 
in the strip $|\Re \sigma|<1$  satisfying the condition {\rm(\ref{eq:kk})} and having poles only 
 at the points $(0,\pm(2-2b)$. Then the following expression
 is holomorphic in the domain $|a|<1-b$:
 \begin{equation}
 \gamma^b(a):=\begin{cases}
          \int_{i\R} U(a,b;0, \sigma)\, d\sigma ,\qquad &\text{if $a\ge 0$},
          \\
          \int_{i\R} U(a,b;0, \sigma)\, d\sigma +4\pi i\, \mathrm{res}_{\sigma=2a}
          U(a,b;0, \sigma) ,\qquad &\text{if $a\le 0$}.
           \end{cases}
           \end{equation}
\end{lemma}

{\sc Proof of Lemma \ref{l:contour}.} Denote 
$$
\gamma_{i\R} (a) =\int_{i\R} U(a,b;0, \sigma)\, d\sigma 
,\quad
\Xi_\pm(a):=2\pi i\, \mathrm{res}_{\sigma=\pm 2a}  U(a,b;0, \sigma).
$$
Since $U$ is even in $\sigma$, we have $\Xi_-(a)=-\Xi_+(a)$.
Due to the factor $(k+\sigma)(k-\sigma)$ in the Plancherel density,
we have $\Xi_\pm(0)=0$. Therefore $\Xi_\pm(a)$ are holomorphic
in the disk $|a|<1-b$.

Consider a contour $L$ on the plane $\sigma\in\C$ composed of the ray
$(-\infty,b-1+\epsilon]$, the upper half of the circle
$|\sigma|=1-b-\epsilon$ and the ray $[1-b-\epsilon, +\infty]$.
The function
$$\gamma_L(a):=\int_{L} U(a,b;0, \sigma)\, d\sigma$$
is holomorphic in $a$ for $|a|<1-b$. For $\Re a>0$ we have
$\gamma_L(a)=\gamma_{i\R}(a)-\Xi_+(a)$. For $\Re a<0$ we have
$\gamma_L(a)=\gamma_{i\R}(a)-\Xi_-(a)$. This gives us 
the analytic continuation.
\hfill $\square$

\sm

{\sc Proof of Theorem \ref{th:symmetry} for
$a<0$.} Thus the analytic continuation of 
$ 4\pi^2 i\la \Phi,\Psi\ra_{L^2(\Lambda,\kappa_{a,b})}$ to the domain $a<0$
is given by
$$
\gamma^b(a)+\sum_{k\ne 0}\int_{i\R} U(a,b;k,\sigma)\,d\sigma,
$$
i.e., for $a<0$ we get $4\pi^2 i\la \Phi,\Psi\ra_{L^2(\Lambda_\C,dK_{a,b})}$.

Now we see that  both sides of (\ref{eq:Lsym}) are real analytic
in the parameter $a$ and coincide for $a>0$. Therefore they coincide for
$a<0$. \hfill $\square$



\section{The operator $J_{a,b}$ is an isometry\label{s:isometry2}}

\COUNTERS

Here we prove the second part of Theorem \ref{th:spectral1}.

\sm

{\bf \punct Statement.}

\begin{lemma}
 \label{l:J-embed}
 Let $f$, $g$ be smooth compactly supported functions on $\dot \C$. Then
 $$
 \la J_{a,b} f, J_{a,b} g\ra_{L^2(\Lambda,dK_{a,b})}=\la f, g\ra_{L^2( \dot \C,\mu_{a,b})}
 .$$
\end{lemma}

Here a way of a proof is simpler than in  Section \ref{s:isometry1}.
We show that $J_{a,b}$ is a perturbation of a version of the Mellin transform.

\sm

{\bf \punct Orthogonality of packets.}

\begin{lemma}
\label{l:packet2}
 Let $f$, $g\in \cD(\dot\C)$. Let $\supp(f)\cap \supp(g)=\varnothing$. Then
 $$\la J_{a,b}f,J_{a,b}g\ra_{L^2_\even(\Lambda_\C,dK_{a,b})}=0.$$
\end{lemma}

{\sc Proof.} By  Corollary \ref{cor:JCL2} the operator $J_{a,b}$ is continuous as an operator 
$\cD(\dot\C)\to L^2(\Lambda_\C,dK_{a,b})$, by Theorem \ref{th:symmetry} it is symmetric on the image of  $\cD(\dot \C,\mu_{a,b})$.
We consider the difference operators 
$$
\tfrac12(\frL+\ov\frL),\qquad \tfrac1{2i}(\frL-\ov\frL)
$$
and repeat the proof
of Lemma \ref{l:orthogonal}.
\hfill $\square$

\sm

{\bf \punct Decomposition of the kernel.}
Starting from this place we examine the restriction of $J_{a,b} f$
to $\Lambda$.
Recall that the operator $J_{a,b}$ is defined by the formula
\begin{equation}
J_{a,b} f(\lambda)=\int_{\dot\C} f(z)\,\cK(z,\lambda)\mu_{a,b}(z)\,\dd z  
.
\label{eq:Jf}
\end{equation}

Decompose the kernel $\cK(z,\lambda)$ according to
(\ref{eq:}) with $N=3$. We consider $\lambda\in\Lambda$, and therefore we set $\sigma=0$.
Denote by $\omega$ the factor depending on $\lambda$ in the front of the expansion.
We have
\begin{equation}
\omega(\lambda)\,\ov{\omega(\lambda)}=\kappa_{a,b}^{-1}(\lambda)
.
\label{eq:omega-omega}
\end{equation}
Notice also that the expression in
brackets in (\ref{eq:}) has a singularity at $\lambda=0$. 
Denote by $\Theta(\lambda)$ a smooth function, which 
equals 0 for $|\lambda|\le1/3$ and $1$ for $|\lambda|\ge 1/2$.
Represent the kernel as
\begin{multline*}
 \cK(z,\lambda)=\omega(\lambda)
 |1-z|^{b-a} |z|^{-a-b}  \times\\\times
\biggl\{
  \biggl[\Bigl(\frac{t_+(z)}{t_-(z)}\Bigr)^{\lambda|-\ov\lambda}
 +
 \Bigl(\frac{t_-(z)}{t_+(z)}\Bigr)^{\lambda|-\ov\lambda}
 \biggr]
 +\\+
\Theta(\lambda) \biggl[ \Bigl(\frac{t_+(z)}{t_-(z)}\Bigr)^{\lambda|-\ov\lambda}
\sum_{k\ge 0,\,l\ge 0,\, 1\le k+l\le 2}
\frac{\ov\lambda^{-k} \lambda^{\,-l}}
{k!\, l!} A_k(\sqrt{1-z})\, A_l(\sqrt{1-\ov z})
+\\+
\Bigl(\frac{t_-(z)}{t_+(z)}\Bigr)^{\lambda|-\ov\lambda}
\!\!\!\!\!\!
\sum_{k\ge 0,\,l\ge 0,\, 1\le k+l\le 2}
\frac{\ov \lambda^{-k} \lambda^{\,-l}}
{k!\, l!} A_k(-\sqrt{1-z})\, A_l(-\sqrt{1-\ov z})\biggr]
+ 
R_3(z,\lambda)
\biggr\},
\end{multline*}
where $R_3(z,\lambda)$ is a smooth function in $z\in \dot \C$
and $\lambda$,
$$R_3(z,\lambda)=O(|\lambda|^{-3})
\qquad \text{as $\lambda\to\infty$.}
$$
uniformly on compact subsets of $\dot\C$. The summands corresponding to
$k=0$, $l=0$ are smooth at $\lambda=0$, so we do not multiply them by 
the patch function $\Theta(\lambda)$.

Next, we change the variable as in (\ref{eq:p})--(\ref{eq:chet}):
$$
\zeta(p):=\frac{(p+1)^2}{4p}
$$
and represent the operator
$J_{a,b}$ in the form
\begin{multline*}
 J_{a,b} f(\lambda)=\omega(\lambda)
 \int_{\ddot\C} f\bigl(\zeta(p)\bigr)
 |1-\zeta(p)|^{a-b-1/2} |\zeta(p)|^{-1/2} 
 |\xi'(p)|^2
 \times\\\times
 \biggl\{
 \biggl[p^{\lambda|-\ov\lambda}
 +
 (p^{-1})^{\lambda|-\ov\lambda}
 \biggr] 
 +\\+
\Theta(\lambda) \biggl[ p^{\lambda|-\ov\lambda}
\sum_{k\ge 0,\,l\ge 0,\, 1\le k+l\le 2}
\frac{\ov\lambda^{-k} \lambda^{\,-l}}
{k!\, l!} A_k\Bigl(\frac{p-1}{p+1}\Bigr) {A_l\Bigl(\frac{\ov p-1}{\ov p+1}\Bigr)}
+\\+
(p^{-1})^{\lambda|-\ov\lambda}
\!\!\!\!\!\!
\sum_{k\ge 0,\,l\ge 0,\, 1\le k+l\le 2}
\frac{\ov \lambda^{-k} \lambda^{\,-l}}
{k!\, l!} A_k\Bigl(-\frac{p-1}{p+1}\Bigr) {A_l\Bigl(-\frac{\ov p-1}{\ov p+1}\Bigr)}\biggr]
+ R\bigl(\zeta(p),\lambda\bigr)
\biggr\}\,\dd p,
\end{multline*}
where $\ddot\C$ denotes
$
\ddot\C:=\C\setminus\{0,1,\, -1\}
$ as above.

It is convenient to split the operator $J_{a,b}$  into a sum of operators, 
\begin{equation}
J_{a,b}= [V_{0,0}^+ +V_{0,0}^-]  +\sum V_{k,l}^+ +\sum V_{k,l}^-+ V_{\mathrm{rem}},
\label{eq:summands}
\end{equation}
where 
the summands correspond to the summands of the previous formula.
We also denote
$$
\gamma(p):= |1-\zeta(p)|^{a-b-1/2} |\zeta(p)|^{-1/2} 
 |\zeta'(p)|^2.
$$



\sm

{\bf \punct The main term.}

\begin{lemma}
\label{l:mellin}
 The operator $\frac 1{2\pi}(V_{0,0}^++V_{0,0}^-)$ is a unitary operator
 from $L^2(\C,\mu_{a,b})$ to $L^2_\even(\Lambda,\kappa_{a,b})$.
\end{lemma}

{\sc Proof.} 
\begin{multline*}
 \bigl\la (V_{0,0}^++V_{0,0}^-)f,(V_{0,0}^++V_{0,0}^-)g\bigr\ra_{L^2(\Lambda, \kappa)}
 =\\=
 \int\limits_\Lambda \Bigl(\int\limits_{\ddot \C} f\bigl(\zeta(p)\bigr)\,\gamma(p) 
 \bigl(p^{\lambda|-\ov\lambda}+p^{-\lambda|\ov\lambda}\bigr)\,\dd p \Bigr)
\Bigl(\int\limits_{\ddot \C} \ov{g\bigl(\zeta(q)\bigr)}\gamma(q) 
 \bigl(q^{-\lambda|+\ov\lambda}+q^{\lambda|-\ov\lambda}\bigr)\,\dd q \Bigr)\, \wt d\lambda
\end{multline*}
(we also applied (\ref{eq:omega-omega}).
Transform this expression as
\begin{multline}
 \int\limits_\Lambda \Bigl(\int\limits_{\ddot \C} f\bigl(\zeta(p)\bigr)\,\gamma(p)\, |p|^2 \,
 \bigl(p^{\lambda-1|-\ov\lambda-1}+p^{-\lambda-1|\ov\lambda-1}\bigr)\,\dd p \Bigr)
 \times\\\times
\Bigl(\int\limits_{\ddot \C} \ov{g\bigl(\zeta(q)\bigr)}\,\gamma(q)\, |q|^2
 \bigl(q^{-\lambda-1|+\ov\lambda-1}+q^{\lambda-1|-\ov\lambda-1}\bigr)\,\dd q \Bigr)\, \wt d\lambda.
 \label{eq:even}
\end{multline}
Now we apply the remark about  Mellin transforms of even functions from Subsect.
\ref{ss:Mellin-even}.
Keeping in mind (\ref{eq:chet}), we observe that
functions $f(\zeta(p))\gamma(p) |p|^2$ are $\times$-even. Therefore  both integrals
over $\ddot \C$ in (\ref{eq:even})
are Mellin transforms of even functions, and we can apply the Plancherel formula
for the Mellin transform.
We come to
\begin{multline*}
\int_{\ddot \C} f\bigl(\zeta(p)\bigr)\, \ov{g\bigl(\zeta(p)\bigr)}
\,|\gamma(p)|^2\, |p|^4\,\frac{\dd p}{|p|^2}=\\=
\int_{\ddot \C} f\bigl(\zeta(p)\bigr) \,\ov{g\bigl(\zeta(p)\bigr)}
\,|1-\zeta(p)|^{2a-2b}\,|\zeta(p)|^{2a+2b-2} |\zeta'(p)|^2
\times\\\times
\Bigl( |1-\zeta(p)|^{-1} |\zeta(p)^{-1} |p|^2 |\zeta'(p)|^2 \Bigr)\,\dd p.
\end{multline*}
By (\ref{eq:zeta1})--(\ref{eq:zeta2}) the expression in the big brackets is 1. 
Now we return to the variable $z=\zeta(p)$ and get the desired expression
$$
\int_{\dot \C} f(z)\ov{g(z)}\, \mu_{a,b}(z)\, \dd z.
$$

{\bf \punct Other terms.}

\begin{lemma}
\label{l:zero}
The Hermitian form
\begin{equation}
\label{eq:zero}
\{f,g\}:=
\bigl\la J_{a,b} f, J_{a,b}g\bigr\ra_{L^2(\Lambda,\kappa_{a,b})}-
\bigl\la (V_{0,0}^++V_{0,0}^-)f, (V_{0,0}^++V_{0,0}^-)g\bigr\ra_{L^2(\Lambda,\kappa_{a,b})}
\end{equation}
on  $\cD(\dot \C)$ 
 can be written as
 \begin{equation}
 \{f,g\}=\int_{\ddot \C} \int_{\ddot \C} K(p,q)\, f(\zeta(p))\,\ov{g(\zeta(q))} \,\dd p\,\dd q,
 \label{eq:Kpq}
 \end{equation}
where $K$ is a locally integrable function on $\ddot\C\times\ddot \C$
smooth outside the sets $p=q$, $p=q^{-1}$.
\end{lemma}

{\sc Proof.} Expanding $J_{a,b}$ according to (\ref{eq:summands}),
we get many summands in (\ref{eq:zero}).  We wish to show that
 each summand can be written as
(\ref{eq:Kpq}) with its own $K$. Let us start the discussion with the summand
\begin{multline}
\la V_{0,0}^+ f, V_{0,1}^-g\ra_{L^2(\Lambda,\kappa_{a,b})}
=
\int_\Lambda
\boxed{\frac{1-\Theta(\lambda)}\lambda}
\Bigl(\int\limits_{\ddot \C} f\bigl(\zeta(p)\bigr)\,\gamma(p)\, |p|^2\,
 p^{\lambda-1|-\ov\lambda-1}\,\dd p \Bigr)
 \times\\ \times
 \Bigl(\int\limits_{\ddot \C} \ov{g\bigl(\zeta(q)\bigr)}\,\gamma(q)\, |q|^2\,
 A_1\Bigl(-\frac{\ov q-1}{\ov q+1}\Bigr)
 \ov{(q^{-1})^{\lambda-1|-\ov\lambda-1}}\,\dd q \Bigr)\,\wt d\lambda.
 \label{eq:boxed}
\end{multline}
The integral in the first big bracket is the Mellin transform
of the function 
$$
F(p):=f\bigl(\zeta(p)\bigr)\,\gamma(p)\, |p|^2.
$$
The integral in the second big bracket is a complex conjugate to the Mellin transform
of
$$
G(q)=
g\bigl(\zeta(q^{-1})\bigr)\,\gamma(q^{-1})\, |q|^{-2}\,
 A_1\Bigl(-\frac{\ov q^{-1}-1}{\ov q^{-1}+1}\Bigr).
$$
Thus we get
$$
\la V_{0,0}^+ f, V_{0,1}^-g\ra_{L^2(\Lambda,\kappa_{a,b})}=
\int_\Lambda \cM F(\lambda)\,\ov{\cM G(\lambda)}\,\, \frac{1-\Theta(\lambda)}\lambda\,\wt d\lambda
.
$$
Denote by $L(p)$ the inverse Mellin transform of $\tfrac{1-\Theta(\lambda)}\lambda$.
It is easy to see that $L(p)$ is an integrable function with a unique singularity of the type $1/(1-\ov p)$ at
$p=1$. We rewrite our integral as
$$
\int_{\ddot\C} \int_{\ddot\C} L(p  q)\, F(p)\,\ov{G(q)}\,\dd p\,\dd q,
$$
and it has the desired form.

\sm

For other pairs  $V_{k,l}^\epsilon$,  $V_{k',l'}^{\epsilon'}$,
where $\epsilon$, $\epsilon'=\pm 1$,
we have similar calculations. Instead of the boxed factor  
in (\ref{eq:boxed}), we get 
\begin{equation}
\frac{1-\Theta(\lambda)}{\ov \lambda^{k+l'} \lambda^{l+k'}}.
\label{eq:Thetaa}
\end{equation}
For $k+l+k'+l'\le 2$ we repeat the same considerations, in these cases
inverse Mellin transforms of the functions (\ref{eq:Thetaa})
 have (integrable) singularities%
 \footnote{We can refer to corresponding formulas for the Fourier transform,
 see \cite{GGV}, Addendum, Sect. 1.7 (Russian edition) or \cite{GSh}, Sect. B.1.3
(English translation).}
  at $p=1$
of types
$$
\frac 1{1-p},\quad \frac{1-p}{1-\ov p},\qquad \ln|1-p|, \quad \frac{1-\ov p}{1-p}
.
$$

If $k+l+k'+l'\ge 3$, then this expression is integrable in $\lambda$, the triple integral 
is convergent, we can change the order of  integrations
and we immediately get an expression of the form (\ref{eq:Kpq}) with real analytic $K(p,q)$.

\sm 

For the pairs including $V_\rem$
we get  absolutely convergent triple integrals and analytic kernels $K(p,q)$. 

\sm

{\bf \punct Proof of Lemma \ref{l:J-embed}.}
Now let $f$, $g\in \cD(\dot\C)$ have disjoint supports.
Then  both terms in (\ref{eq:zero}) are zero (see Lemma \ref{l:packet2}).
Therefore the kernel $K(p,q)$ satisfy the following property:
$$
\{f,g\}=
\int_{\ddot\C} \int_{\ddot\C}K(p,q)\, \phi(p)\,\ov {\psi(q)}\,dp\,dq=0
$$
if $\phi$, $\psi$ are $\times$-even elements $\cD(\ddot\C)$ with disjoint supports. 

We claim that $\{f,g\}=0$ for any $\times$-even functions $f$, $g\in \cD(\ddot\C)$.
To observe this, we  take a $\times$-even partition of unity $\tau_j$ with small supports,
and decompose
$$
\{f,g\}=\sum_{k,l} \{\tau_k f,\tau_l g\}.
$$
Clearly, we can make this sum as close to zero as we want by refinement
of a partition of unity.
We omit trivial details.

\section{Domains of self-adjointness}

Thus $J_{a,b}$ is unitary. Clearly the multiplication operators
$$
f(z)\mapsto \tfrac12(z+\ov z) f(z), \qquad f(z)\mapsto \tfrac1{2i}(z-\ov z) f(z)
$$
defined on $\cD(\dot\C)$ are essentially self-adjoint in $L^2(\C,\mu_{a,b})$ and commute.
Therefore the operators $\frac12(\frL+\ov\frL)$, $\frac1{2i}(\frL-\ov\frL)$ are essentially self-adjoint and commute on the subspace
$J_{a,b} \cD(\dot\C)\subset L^2_\even(\Lambda_\C, dK_{a,b})$. But $\cW_{a,b}$ contains this image.
This establishes Theorem \ref{th:conditions2}.a. 

Theorem \ref{th:spectral2}.a follows from the same argumentation.

 {\tt
 \noindent
 Vladimir Molchanov\\
 Institute for Mathematics, Natural sciences,
 \\
                     and informational technologies,
 \\                   
 Tambov State University, Tambov, Russia.}
 
 \sm
 
 \tt
 \noindent
 Yury Neretin\\
  Pauli Institute/c.o. Math. Dept., University of Vienna \\
 \&Institute for Theoretical and Experimental Physics (Moscow); \\
 \&MechMath Dept., Moscow State University;\\
 \&Institute for Information Transmission Problems;\\
 URL: http://mat.univie.ac.at/$\sim$neretin/

\end{document}